\documentclass{amsart} 

\newtheorem{theorem}{Theorem}[section]
\newtheorem{lemma}[theorem]{Lemma}
\newtheorem{corollary}[theorem]{Corollary}

\theoremstyle{definition}

\theoremstyle{remark}
\newtheorem{remark}[theorem]{Remark}

\numberwithin{equation}{section}

\newcommand*\xedge[1]{$\overline{\text{#1}}$}

\usepackage[margin=1in]{geometry}

\usepackage{amsmath ,amsthm ,amssymb} 
\usepackage{graphicx} 
\usepackage[normalem]{ulem}
\usepackage{times}
\usepackage{txfonts}
\usepackage{url,hyperref}
\usepackage{enumitem}
\usepackage{pdfpages}
\usepackage
{appendix}
\usepackage{afterpage}
\usepackage{booktabs}

\newcounter{v8LemmaCounter}

\usepackage{lineno}
\usepackage{setspace} 

\title{Another Proof of the Four Colour Theorem - Part 2 - Discharging a minimal 5-Chromatic Planar Graph}
\author{Frank Allaire}
\email{frank.allaire@lakeheadu.ca}
\date{Sept. 16, 2022}		

\begin{document} 
\begin{abstract}
In RSST \cite{RSST}, they "replace the mammoth hand-checking of unavoidability that A\&H required, by another mammoth hand-checkable proof" (page 18). Here, the proof of unavoidability is accomplished in a lengthy structured hand-checkable proof whose entirety is presented in this document.
\end{abstract}
\maketitle

\section{Summary and preview}
To disprove the Four Colour Theorem (4CT) , one needs only describe a counterexample. To prove 4CT,  we assume a counterexample exists, and by proving properties of such a counterexample, eventually show that it cannot exist.
Every loopless planar graph has a vertex of degree 5 or less and 
is easily shown to be
vertex 5-colourable. A graph is k-chromatic if it admits a k-colouring but no (k-1)-colouring.

Assuming a 5-chromatic planar graph exists, then there must be (at least) one with an minimum number of vertices, a minimum 5-chromatic planar graph (m5CPG). Well known properties of a m5CPG are:\\
\indent	- it is a triangulation\\
\indent	- there are no vertices of degree less than 5\\
\indent	- every separating 5-circuit is the rim of a 5-wheel (the graph is almost 6-connected).\\
Also, a m5CPG cannot contain a reducible configuration: a cluster of vertices which can be replaced by a smaller cluster, and every 4-colouring of the result can be manipulated, typically proved possible by a Kempe chain argument,  into a 4-colouring of the assumed m5CPG \cite{FAPart1}.

The remaining discussion relates only to an assumed m5CPG.

\subsection{Unavoidable sets for a m5CPG}
From Euler's formula for non-null planar graphs, [V] - [E] + [F] = 2, if we apply a "charge" of +1 unit to each vertex, -1 to each edge, and +1 to each face, then the total of these charges over the entire planar graph must be +2. Given a planar triangulation, move all the charges to the vertices: -1/2 to the ends of each edge, and +1/3 to the vertices of each triangular face. Each vertex of degree n ends up with 1 - n/2 + n/3 = (6-n)/6 units of charge. To develop and maintain integer values later, multiply these values by 60 so each vertex of degree n has a charge of 10(6-n) and the total over the entire graph is +120 units.


If for every planar triangulation that avoids a specified set of configurations, we can redistribute the charges and prove that under this redistribution, every vertex has zero or a negative charge, then that set of configurations must be unavoidable. If all the configurations in the set are reducible, then the Four Colour Conjecture (4CC) is elevated into 4CT.

A discharging rule redistributes units of charge from a vertex in a specified cluster, the source vertex, to an indicated nearby vertex, the sink vertex. 
The sink vertex for one discharging rule could be the source vertex for a different discharging rule.

In this paper, an unavoidable set S of reducible configurations is grown by adding individual reducible configurations or sets of reducible configurations to S. Each addition validates a particular property of the assumed m5CPG and/or resulting properties of the discharging scheme as it is applied to a m5CPG. When a set of reducible configurations is added, it is possible that some elements of the set are already in S, or that not every reducible configuration in the set is necessary to validate the associated property. In particular, the unavoidable set described here is not minimal with respect to size. To begin, S is the empty set.


5-valent vertices start with a charge of +10 units and 6-valent vertices start with a charge of zero units. A desirable goal of any discharging of a m5CPG should be to send exactly 10 units from each 5-valent vertex, and the final targets of these transfers should not be 5- or 6-valent. The final charge on every 5- and 6-valent vertex will then be zero.
 \\
 \\
 \\
Notation:\\
 A vertex of degree 7 or greater is called a Major vertex, possibly indicated by M,  while vertices of degree 5 or 6 are minor vertices, possibly indicated by m.\\
(from \cite{FAPart1} Section 3, slightly modified) 
... A cluster configuration is composed of a boundary circuit Q whose inside is triangulated by a specified structure of vertices and edges, forming a near-triangulation T = (Q;T). No inside edge is a diagonal of Q. Since a cluster is a near triangulation without diagonals, its structure can be described by indicating the degrees  and adjacencies of the interior vertices. We use the following notation. The string 5-6[5665] -9 indicates a cluster with five (the leading 5) inside vertices: a central 6-valent vertex (before the square brackets) adjacent in order to vertices of degree 5, 6, 6 and 5 (inside the brackets). These vertices are called the first neighbours of the central vertex. These inner vertices triangulate a boundary circuit of nine (the last descriptor) vertices. As a place holder, we use an x to indicate a boundary vertex whose degree is unspecified, and X or B if it is unspecified but restricted to being Major. The configuration 6-8[56(5)75]-12 is similar to the one above except it has a '5cap' (in parentheses): a 5-valent vertex adjacent to the
 preceding and following
 specified vertices on the side away from the central vertex,  8-valent in this case. To reduce the number and levels of bracketing, a cap 
  may be abbreviated 0 for (5), 1 for (6), 2 for (7) and 3 for (8).
Occasionally a 5- or 6-valent vertex at distance 2 or 3 from the central vertex is attached to the preceding
 vertex but not always to the following vertex of the description. This is indicated by ((5)) or ((6)). Once the remaining vertices of specified degree are drawn, the remaining connection or connections are forced by the cluster being a triangulation. See Figure ~\ref{fig-RedUsedSoFar} for examples.

For large configurations, it may be desirable to use two adjacent central vertices with the vertex adjacent to both central vertices delimited by hyphens. Thus, 8-5[6078(5)05]-12 is also 78[56-5-555]. Since the cluster is a near triangulation, the prefix and suffix can be derived from the remainder of the description and are usually omitted.


In figures and text, the degree or range of degrees of a vertex is indicated 
by the symbols below.
Using this key it is convenient and usually possible to describe a cluster of vertices by drawing only the interior vertices of a configuration and the edges joining these. If necessary to avoid ambiguity, or to emphasize a feature, some of the remaining edges or boundary vertices can be indicated.\\
 \\
\indent \indent
\begin{tabular}{ r |  l   c c   r |  l }
\hline
\label{tabl-symbols}
5-valent: 	&$\bullet$& & &at least 5-valent:&x or (null)\\
6-valent: &$\triangle$& & &at least 6-valent: &* \ or equivalent\\
7-valent: &$\medcirc$& &	&at least 7-valent:&$\subset$ or M or B\\
8-valent: &	$\square$& & &at least 8-valent:& $\sqsubset$\\
9-valent: &$\bigtriangledown$& & & & \\
minor: &m& & & &  \\
Major: &M& & & & \\
\hline
\end{tabular}\\
\indent \indent \indent \indent \indent \indent Degree codes for vertices.
\section{The Discharging Rules - Part I: Transfers across (minor,minor) edges}
\subsection{Discharging 5-valent vertices}\label{sec-Disch5}
\begin{figure}[th!b]
  \begin{center}
    \includegraphics[width=6in]{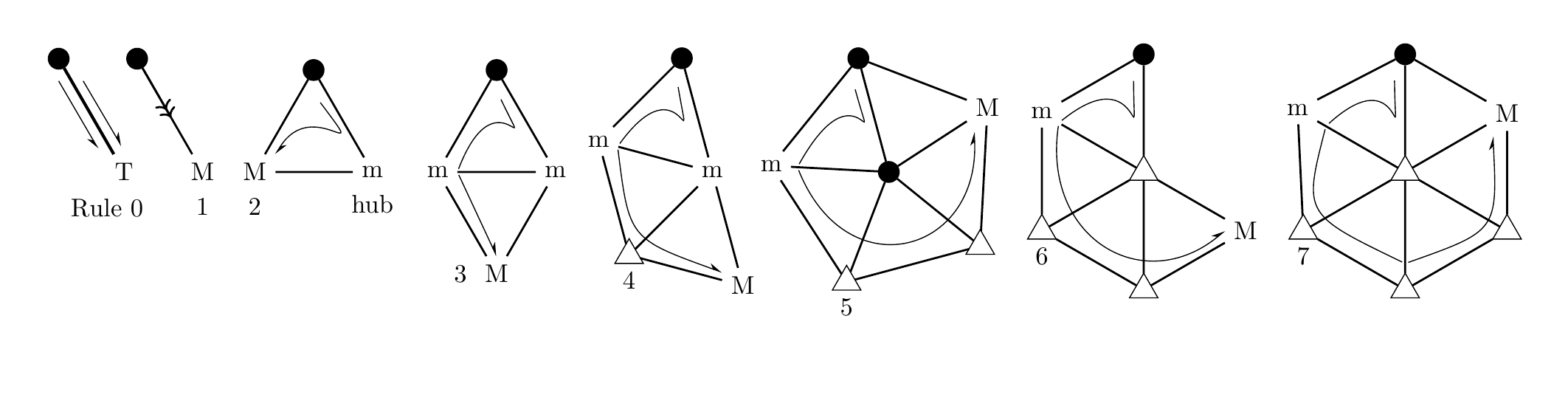}
  \end{center}
  \caption{ Rules 0 to 7:  Discharging 5-valent vertices.}
  \label{fig-Rules0to10}
\end{figure}
A very slick discharging for 5-valent vertices was described by J. Mayer  \cite{Mayer}, \cite{Mayer2} and extended by RSST. Rules 1 to 7 presented in Figure ~\ref{fig-Rules0to10} have exactly the same charge transfers as RSST from source vertex to sink Major vertex but 
the trajectories are slightly different and even these have an equivalence mapping except for Rule 2.

Rule 0: Send 1 unit of charge out along each side of each spoke of the 5-wheel.

Rule 1: If the first target of this transfer is Major, that Major vertex is the sink for both units of  charge. 
It is convenient to consider these two units as flows "inside" the (5,Major) edge to distinguish these from 
other trajectories that remain "along" or "beside" this and later edges, ending at a Major vertex labelled M, the sink for each of these flows.

Rule 2: Otherwise, this first target is minor and is used as the hub for this charge transfer. The charge flow reflects away from the minor hub vertex and arcs around the hub vertex toward a second target, the third vertex of its surrounding triangular face.

If this second target of this transfer is Major, that Major vertex is the sink for this charge transfer. 

Rule 3: Otherwise, this second target is also minor and the charge crosses the (hub, minor vertex) edge and proceeds beside the rim edge of the hub vertex toward the next neighbour of the hub, the third target.
5[555], 5[565], and 6[565] are all reducible, so adding these to the set to be avoided ensures that this third target is not 5-valent.

If this third target of this transfer is Major, that Major vertex is the sink for this charge transfer. 

Rule 4: Otherwise, this third target is 6-valent and the charge crosses the (hub, 6) edge toward the next neighbour of the hub, the fourth target.
Adding reducible 5[5665], and 6[5665] to the set to be avoided ensures that this fourth target is not 5-valent.

Here and in Rules 5 and 6, if the second target is 5-valent, then the reducible configurations needed to force the new target to be at least 6-valent have already been added to the set to be avoided. For the purpose of determining the reducible configurations sufficient to force the new target to be at least 6-valent, we can assume that the second target is 6-valent.

If this fourth target of this transfer is Major, that Major vertex is the sink for this charge transfer.

Rules 5 and 6: Otherwise, this fourth target is 6-valent and the charge crosses the (hub, 6) edge toward the next neighbour of the hub, the fifth target.

Rule 5. If the hub is 5-valent, this fifth target is the final neighbour of the hub. Adding reducible 5[56665] and 5[56666] to the set to be avoided assures that this fifth and final target is Major and that Major vertex is the the sink for this charge transfer.

Rule 6: Otherwise the hub is 6-valent. Adding reducible 6[56665] to the set to be avoided assures that this fifth target is not 5-valent.\\
If this fifth target of this transfer is Major, that Major vertex is the sink for this charge transfer.

Rule 7: Otherwise the hub is 6-valent and this fifth target is 6-valent.  The charge crosses the (6-hub, 6) edge toward the final neighbour of the hub.

\noindent Adding reducible 6[566665]  and 6[566666] to the set to be avoided assures that this sixth and final target is Major and that Major vertex is the sink for this charge transfer.

\subsubsection{Mayer's discharging and the minor maxim}
Combining Rules 1 and 2, a 5-valent vertex sends 2, 3,or 4 units to each adjacent Major neighbour, depending on whether respectively  0, 1, or 2 of their common neighbours, also called shoulder vertices, is/are minor. Rule 3 sends 2 units from a 5-valent vertex across each minor-minor edge of its rim. Those were Mayer's rules, and his proof that it discharged the 5-valent vertex enumerated the 8 necklace patterns of a 5-cycle with beads of two types: minor and Major. In every case, the outflow from a 5-valent vertex is exactly 10 units. 
Mayer's result was empirical, almost serendipitous. Here, it is explained as a consequence of Rule 0 together with the maxim that a charge does not stop at minor vertices as it travels around the hub vertex to its Major destination. 
Although the third target was known to be at least 6-valent, Mayer did not specify where the 2 units that flow across the minor-minor edge should end up. The continuation around the minor shoulder vertices was the obvious choice as soon as 6[566666]-11  and the smaller configurations were known to be reducible
.
By avoiding reducible configurations, the Major destination must exist.

\subsubsection{Preview of Unavoidability}
So far, the set to be avoided contains only reducible configurations. Adding two non-reducible configurations makes the set unavoidable and its proof demonstrates one of several approaches taken in this paper.

\begin{lemma}
Let S be the set of reducible configurations mentioned so far and let U be the set S to which the (5,5) edge and the 7-valent vertex are added. This set U is unavoidable for a m5CPG.
\label{lem-U-is-Unavoidable}
\end{lemma}
\begin{proof} Eliminating the (5,5) edge limits a hub vertex to being 6-valent and also eliminates the co-terminal (co-sink?) transfers which occur when the second target is 5-valent.
After the discharging, final charge on each vertex is as follows:
\begin{enumerate}
\item  5-valent vertex: zero units. The initial charge on each 5-valent vertex is ten units and ten units are distributed to Major vertices by the discharging.
\item  6-valent vertex: zero units.  The initial charge on each 6-valent vertex is zero and is unchanged by the discharging.
\item For a Major vertex M, consider the inflow along each spoke of its wheel.\\
If the end of a spoke is also a Major vertex,  that spoke contributes zero charge: no rule has flow along a (Major, Major) edge.\\
Each 6-valent neighbour may support an inflow of at most one unit on each side of its spoke and the vertex opposite any 1 unit flow beside a spoke must be its 6-valent hub. A 6-valent rim vertex with 2 units flowing along its spoke (1 unit on each side) must be the middle vertex of a 666 triple of consecutive rim vertices.\\
Each spoke from a 5-valent neighbour delivers 2, 3 or 4 units from rules 1 and 2 but no more.
With no (5,5)-edge, even Rule 3's transfer ends along a (6,M) edge.
In any case, aggregate the contribution of each 5-valent rim vertex together with its non-5-valent clockwise neighbour on the rim for a total of at most 5 units from the pair of vertices. Non-paired rim vertices are 6-valent or Major and contribute at most 2 units each.
The final charge on a Major vertex of degree $n$ is at most
 $10(6-n) + 2.5n = 60-7.5n \le\,0\text{ for } n\ge8$.
\end{enumerate}
\end{proof}
In summary, applying the discharging rules to a m5CPG which avoids the configurations in $U$ results in every vertex having non-positive charge, so such a m5CPG cannot exist.
Another consequence of U being unavoidable is that a m5CPG must contain either a (5,5) edge or a 7-valent vertex since it cannot contain a reducible configuration.

\subsection{More cross minor-minor edge transfers}\label{More mm Transfers}
\begin{figure}[th!b]
  \begin{center}
    \includegraphics[width=6in]{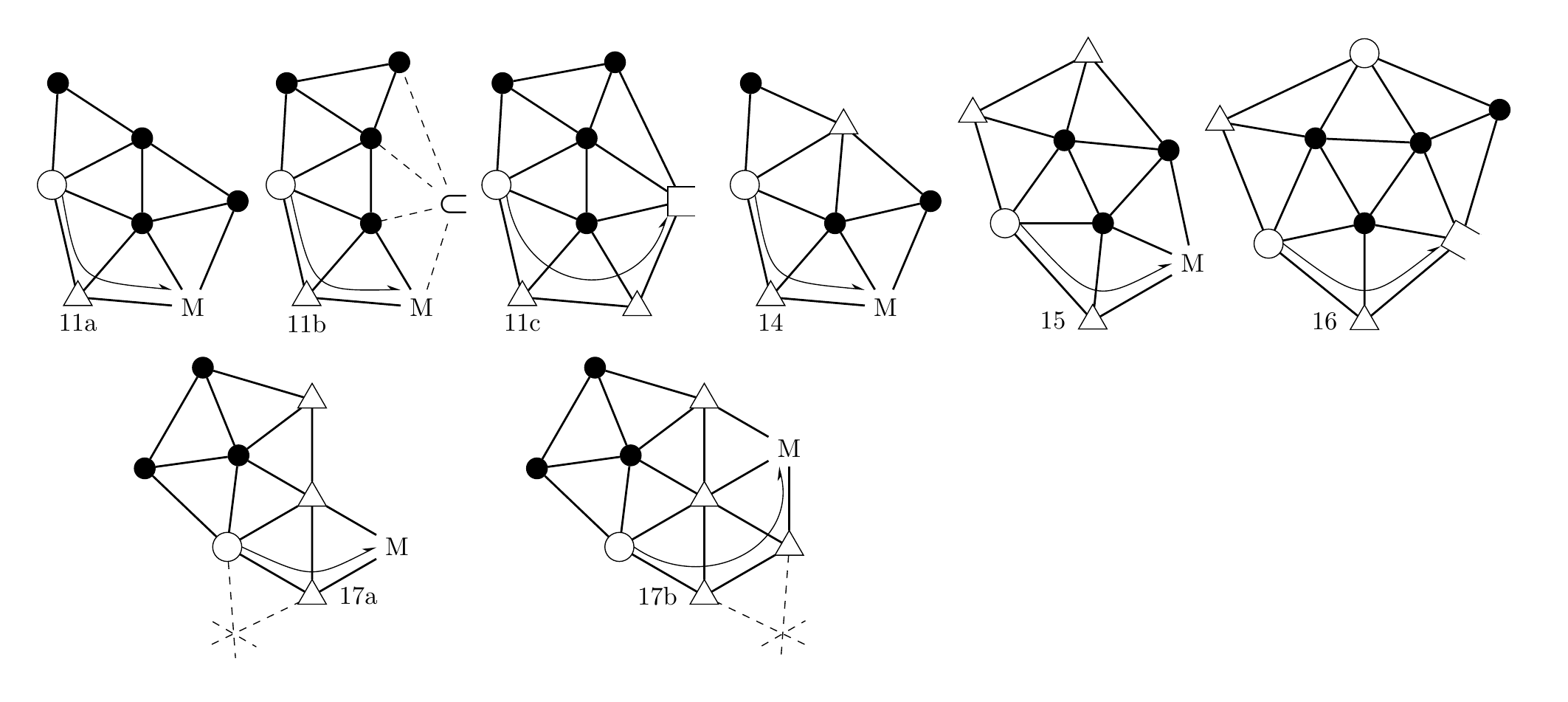}
  \end{center}
  \caption{ Rules 11 to 17:  Other cross (minor, minor) edge transfers.}
  \label{fig-Rules11to17}
\end{figure}
	
The discharging rules 2 to 7  all involve a source 5-valent vertex sending 1 unit of charge around an adjacent minor hub vertex, across zero or more (minor hub, minor) edges to the first Major vertex adjacent to the hub. 
Rules 11a to 17b of Figure ~\ref{fig-Rules11to17}
describe charge transfers from a 7-valent vertex 
around an adjacent minor hub vertex, across one  or two (minor hub, minor) edges to the first Major vertex adjacent to the hub. 

\begin{figure}[th!b]
  \begin{center}
    \includegraphics[width=6in]{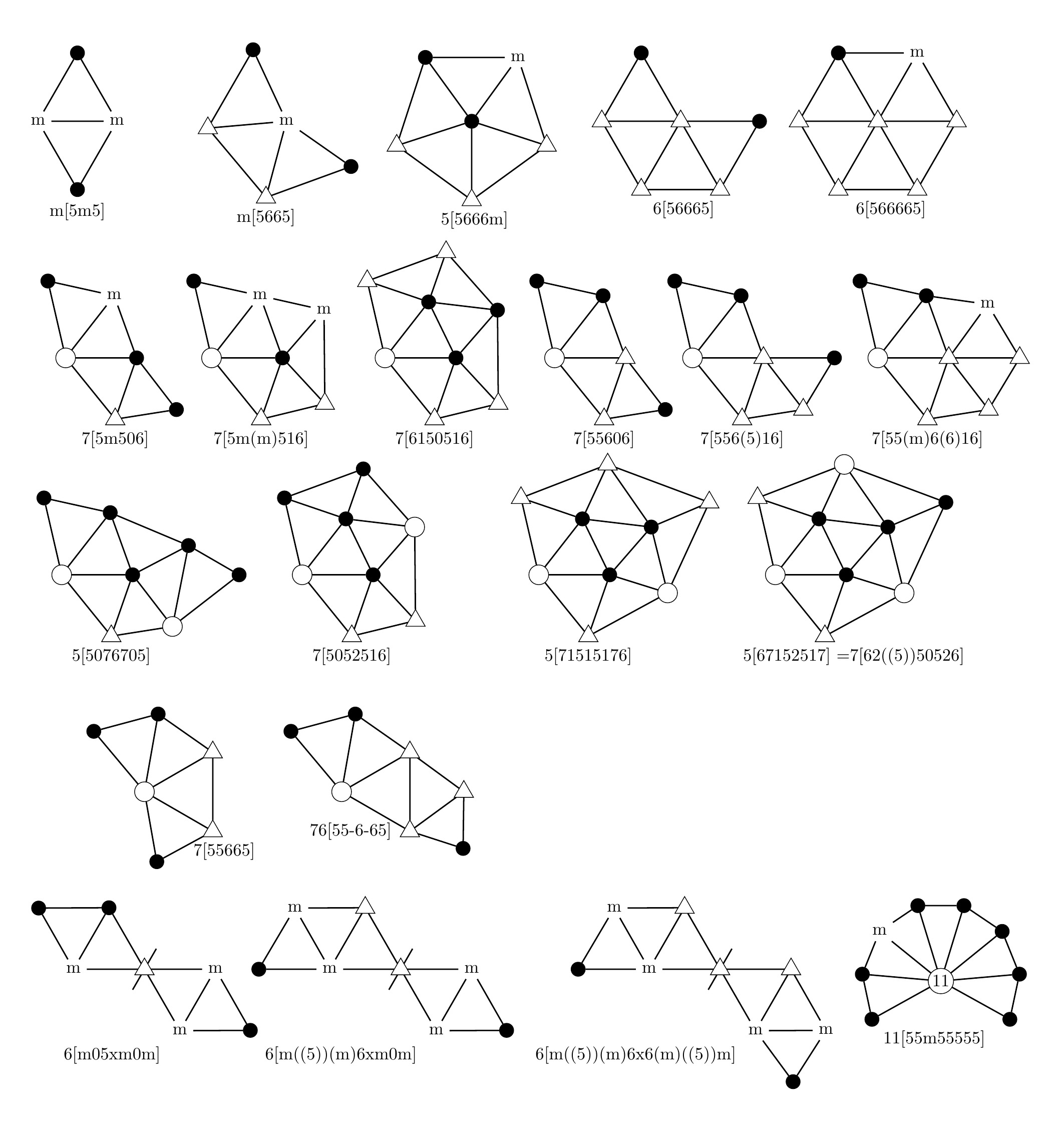}
  \end{center}
  \caption{ Reducible Configurations Used up to Section ~\ref{sec-Cross mm edge Limits}.}
  \label{fig-RedUsedSoFar}
\end{figure}

In the description of Rules 3 to 7,  reducible configurations were added to the set to be avoided to ensure that a Major target existed and hence the outflows of rule 0, the basic source structure, were assured.
Here, underlying source configurations with one unit of outflow could be 7[555(T)6], 7[565(T)6], 7[61505(T)6], 7[62((5))505(T)6] and 7[556(T)6], because reducible configurations force the existence of a Major target around the hub (Figure ~\ref{fig-RedUsedSoFar}, line 2). 
However, that generality is not needed and in fact complicates the analysis of the discharging.  

Several diagrams in Figure ~\ref{fig-RedUsedSoFar} use m to represent a vertex that may be either 5-valent or 6-valent. When fully specified, each configuration is  reducible with an appropriate reducer if necessary.  This just one of the ways  a single diagram will represent as many as over a hundred reducible  configurations.

\begin{remark}
From now on, any indication that a configuration( or set of configurations) "is reducible" will imply "and is added to (or Unioned with) $S$, the set to be avoided". Also, "is the Major target" will imply "and is the sink vertex for the charge transfer". A set of configurations is reducible iff every configuration in the set is reducible.
\end{remark}

\subsubsection{One-way Transfers}\label{ssec-78xmm1-way}
The descriptions of the discharging rules for a 5-valent vertex include reducible configurations to ensure that a Major sink exists, the rule applies, and the discharge occurs. In particular, these and future discharging rules describe a one-way transfer from the source vertex to the Major sink vertex. 

With the extra vertices, only 11a+11a, 11c+11c, 15+15, and 16+11a are compatible and without  m[5m5]. The four reducible configurations in  line 3 of Figure~\ref{fig-RedUsedSoFar} ensure that if the Major target of rules 11 to 17 
is 7-valent, then it is not part of a  
structure with that 7-valent target as a source
sending charge back to the first 7-valent source vertex.
As a bonus, the lower bounds on the target vertices for Rules 11c and 16 are now established. Two more reducible configurations, 7[55665] and 76[55-6-65] ensure that 
if the source structures 17a and 17b occur in a m5CPG, the dashed-attached vertices must have the indicated minimum degrees, as in 11b by reducible m[555].


\subsection{Preview of cross Major-minor edge transfers}  \label{sec-CrossMm}
(Figure ~\ref{fig-CrossMm})  To complete this discharging for 4CT, unlike in Lemma ~\ref{lem-U-is-Unavoidable} there are charge transfers from a 7- or 8-valent vertex to an adjacent Major vertex. Again, each flow is one unit from a specified source to a specified target along a specified trajectory.
Here, each trajectory crosses one (Major,minor) edge, one unit with the Major vertex as the hub, and occasionally a second unit crossing the same edge but with the minor vertex as the hub. If both shoulder vertices of the (Major source, Major sink) edge are minor, there is a possible transfer of four units from a 7-valent vertex to the same adjacent Major vertex, each transfer with a different trajectory. Transfers from an 8-valent vertex to an adjacent Major vertex, M, occur only when both shoulder vertices are 5-valent and this flow to M is along a (5,M) edge with the  8-valent vertex as the hub.

\subsection{Radial inflow.}\label{sec-MaxRad-Inflow}
The total inflow to a Major sink vertex inside and along both sides of a spoke is called the radial inflow from that spoke, or, from the rim vertex at the other end of the spoke, regardless of location of the actual source vertex.
\begin{lemma} (Figure ~\ref{fig-Max-Radial-Inflows}) The radial inflow to a major vertex M is:
\begin{enumerate}
\item at most 2 units from a major vertex and is positive only when that source vertex is 7-valent, 
\item at most 3 units from a 6-valent vertex.
\item 6 units from the 5-valent vertex in the cluster M[60506], otherwise at most 5 units.
\end{enumerate}
\label{lem-MaxRadial}
\end{lemma}

\begin{figure}[ht]
  \begin{center}
    \includegraphics[width=3.in]{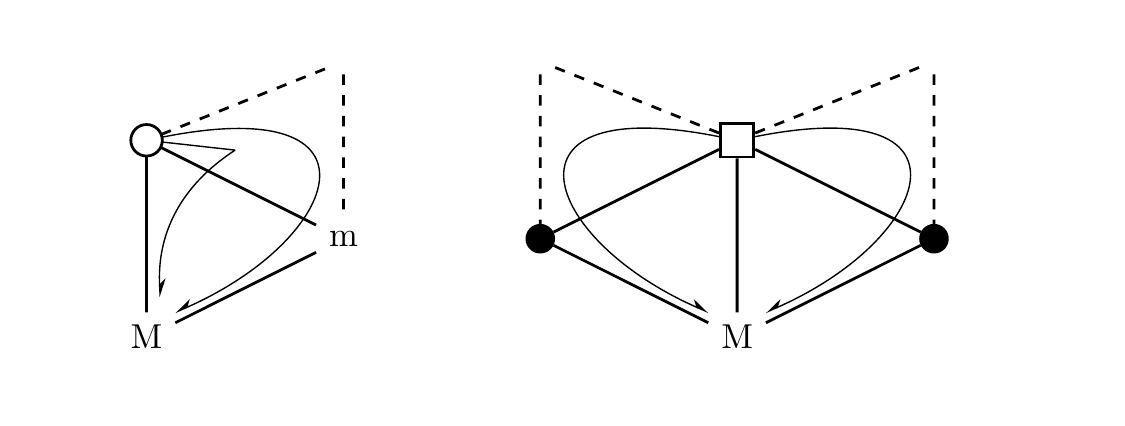}
  \end{center}
   \caption{Preview of cross (Major,minor) edge transfers.
  \label{fig-CrossMm}			}
\end{figure}
\begin{figure}[ht]
 \begin{center}
    \includegraphics[width=6in]{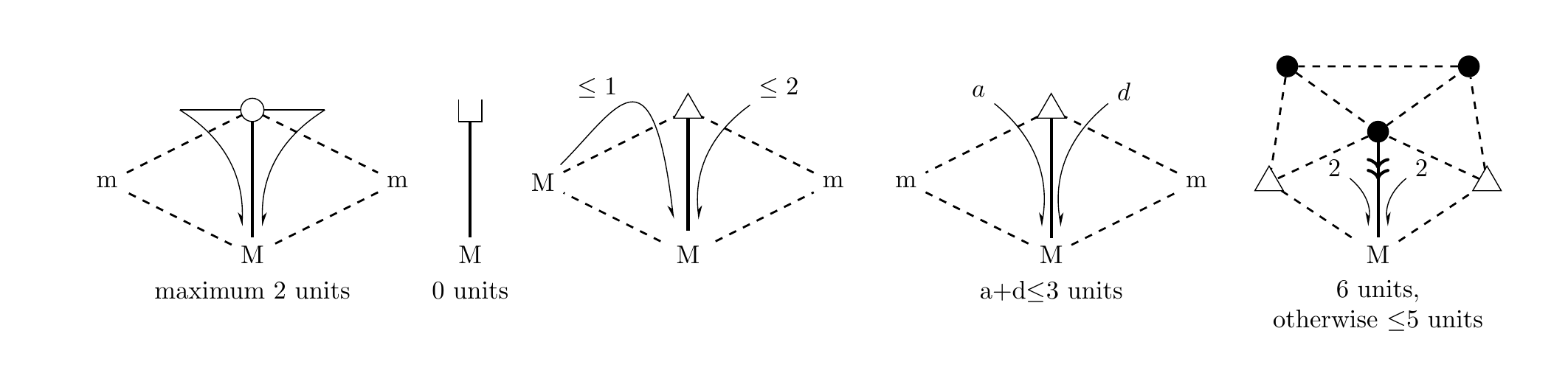}
  \end{center}
  \caption{Maximum Radial inflows from a Major, 6-valent or 5-valent vertex}
  \label{fig-Max-Radial-Inflows}
\end{figure}

\begin{proof} 

(1) follows from Figure ~\ref{fig-CrossMm} and its description in Section ~\ref{sec-CrossMm}.

(2): Inflow to M along a side of a spoke from a 6-valent vertex occurs only when the shoulder vertex of this edge is the hub of the transfer. With itself as hub, a Major shoulder vertex generates at most one such unit . A minor shoulder vertex as hub directs flow across one or more (minor, minor hub) edge(s) to a sink M only if the hub has a 5-valent or sourcing 7-valent neighbour upstream.
For two units to flow along a side of a (6,M) edge to M, the hub must be minor and the upstream neighbourhood of this hub must be one of hub[M655], hub[M6655] or  (hub=6)[M66655].
Every source structure driving 2 units of flow along both sides of a (6,M) edge to M contains a reducible configuration
shown in Figure ~\ref{fig-RedUsedSoFar}.

(3): The maximum inflow to a Major vertex along each side of a spoke from a 5-valent vertex is 2 units and this occurs only with a minor hub and a 5cap on the (5, hub) edge, Rules 2 and 3. Together with Rule 1, the flow "inside" the edge, the maximum flow to M is 6 units and requires 5caps on both  (5, hub) rim edges. Avoiding 5[555] forces both minor hubs to be 6-valent.
\end{proof}
\section{Vertices of high degree are not overcharged}
At this point, as in RSST, we can show that vertices of sufficiently high degree cannot be overcharged. Their discharging has the feature that all transfers were "inside" edges, and allows transfers to and from minor vertices.
They showed that
no more than 5 units flows inside any spoke to a Major vertex. 
From this,  they conclude that the final charge on a Major vertex of degree $n$ is at most
\indent $10(6-n)+5n = 60-5n\le0\text{ for }n\ge12$. 
\begin{theorem}\label{12+notOver}
Vertices of degree 12 or greater are not overcharged.
\end{theorem}
\begin{proof}
This discharging has one possibility of 6 units coming in along a spoke, the 5-valent  rim vertex of the cluster M[60506], but there are two adjacent 6-valent vertices, each with a radial inflow of two or perhaps three units. Furthermore, M[605060506] contains 6[505x505x]-8, a reducible configuration, so each [60506] neighbour cluster must be distinct.
Grouping the radial inflows from the three rim vertices of  a 60506 cluster together, the total inflow from the three spokes is at most 12 units, averaging at most 4 units per spoke. All other  5-valent vertices give a radial inflow of at most 5 units. 
Similar to RSST, the final charge on a vertex of degree $n$ is at most $10(6-n)+5n = 60-5n\le0\text{ for }n\ge12$.
\end{proof}
\begin{remark} In describing the neighbourhood of a central vertex, power notation is used to abbreviate consecutive vertices of a given degree or a sequence of clusters of a given shape. 11[55555555xxx] can be written 11[$\text5^\text8$xxx] and 11[55655555xxx] is 11[$\text5^\text2$6$\text5^\text5$xxx], and together, the pair of (reducible) configurations is 11[$\text5^\text2$m$\text5^\text5$].

Another way to express a family of configurations is by using braces to indicate that every permutation of the contained vertices or structures is allowed. Inside these braces, power notation no longer indicates consecutivity.
10[5\{$\text5^\text4$,m\}5xxx] denotes 10[5555555], 10[5655555], 10[5565555], and 10[5556555] (all of which are reducible), bounded by a circuit of size 11 or 12. Other families of reducible configurations are: 10[5\{$\text5^\text4,\text{m}^\text2$\}5xx]-13, and
10[5\{5,5,m\}60506xxx]-13, the trailing  -13 indicating the size of the largest boundary circuit, obtained when each m=6.
\end{remark}

\begin{theorem}\label{11notOver}
Vertices of degree 11 are not overcharged.
\end{theorem}
\begin{proof}
11[55555555xxx]-12 and 11[55555655xxx]-13 are reducible.\\ 
Recall that $\ast$ indicates a vertex of degree 6 or more.
We evaluate the final charge  via an enumeration of the neighbourhoods of an 11-valent vertex in decreasing  order of the number of 5-valent neighbours:

11[ $\text5^\text{11}$] and 11[ $\text5^\text{10}\ast$] contain reducible 11[$\text5^\text8$xxx] which also eliminates neighbourhoods from11[\{$\text5^\text9$,$\ast^\text2$\}] when the * 's are consecutive or one 5 apart. Otherwise  the * 's are at least two 5's apart in which case avoiding 11[55655555] forces the $\ast$'s to be Major. Since there are no [656] neighbour runs, the final charge is at most --50 +9$\cdot$5+2$\cdot$2=-1.

With only three neighbours of degree at least 6, 11[\{$\text5^\text8$,$\ast^\text3$\}] has at most one 60506 neighbour group, so the final charge is at most --50+1$\cdot$6+7$\cdot$5+3$\cdot$3=0.

11[\{$\text5^k$,$\ast^{(11-k)}$\}] contains at most three non-overlapping 60506 neighbour groups and the final charge is at most\\ 
--50+3$\cdot$6+(k--3)$\cdot$5+(11--k)$\cdot$3=2k--14 $\le$ 0 for $k$ $\le$ 7.
\end{proof}

\section{ Cross (minor, minor)-Edge Flow Limits, single or adjacent pairs } \label{sec-Cross mm edge Limits}

The same approach is used in showing that 10-valent vertices are not overcharged: radial inflows exclude the neighbourhoods with a small number of 5-valent vertices, and reducible first neighbourhoods exclude or limit the inflow from neighbourhoods with a large number of 5-valent vertices.
In between, and for vertices of degrees 7, 8, or 9, a more  detailed evaluation of possible inflows is needed. 
These enumerations and evaluations have some common elements but each has its own tailor-made lemmas and approach.

The Mayer values of 2, 3, or 4 units account for the inflow from 5-valent neighbours according to Rules 1 and 2. All other flows to a Major sink vertex come across the rim edges of its wheel, either (minor, minor) or (Major, minor). No discharging trajectory crosses (Major, Major) edges. After crossing this rim edge, the charge flows beside the spoke opposite its hub to its Major target.
Coupling the flows on both sides of a spoke gave the radial inflow limits.
Coupling the flows across (minor,minor) rim edges leads to Cross mm-edge flow limits.

\begin{remark}Perhaps here is a good place to give reading advice: (m,m) should be read "minor minor"  rather than "em em" and similarly, M should often be read as "Major", m0m could be read as "minor, 5cap, minor". "Star" has fewer syllables  than "asterisk". An integer "k" may be an abbreviation for  the more accurate "k-valent vertex".
\end{remark}

\begin{figure}[ht] 
  \begin{center}
    \includegraphics[
    width=6in
    ]{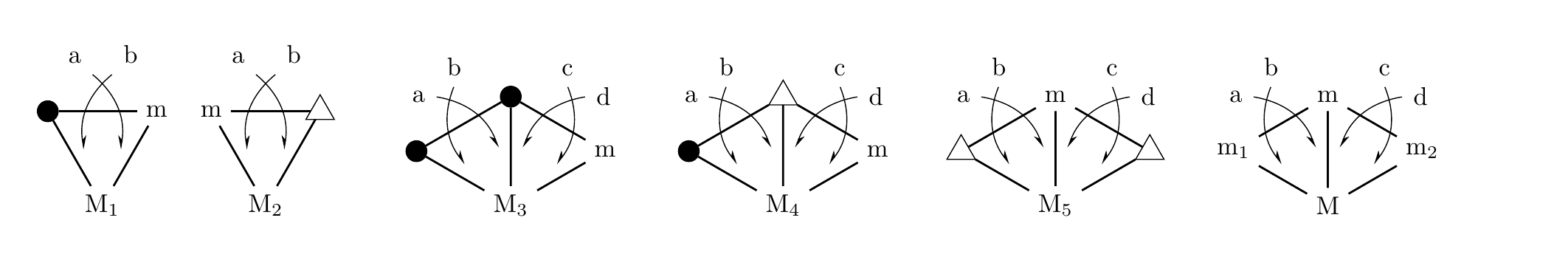}
  \end{center}
  \caption{Cross mm edge limits - Lemma ~\ref{lem-cross-mm}}
  \label{fig-Cross mm edge limits}
\end{figure}

\begin{lemma}\label{lem-cross-mm} The following limits apply to cross (m,m) edge transfers to a Major sink, rules 3 to17. 
Use Figure ~\ref{fig-Cross mm edge limits} to disambiguate the flows ( possibly summed from two adjacent upstream 5-valent source vertices ) named a, b, c and d, each of which is at most 2 units and crosses a (minor hub, minor) edge, continues along the (minor, M) edge and is sinked at M. Unlike the radial limits of Lemma ~\ref{lem-MaxRadial}, these transfers exclude Rule 2, the inside a (5,m,M) triangle transfers from 5 to M, which are instead included in the Mayer 2,3,4 values.

\begin{tabular}{ r r l}		

\phantom{0}1&$M_{1}[55]:$&$
a=b=1$$\iff$$5$-$cap; otherwise\  a=b=0.$\\

2&$M_{1}[56]:$&$
b=1$$\iff$$5$-cap; otherwise b=0; ($\implies$a+b $\leq$ 3)\\
&&$a+b=3$$\iff$$ M_{1}[5(5)0(*)6],\ the\ 55$-cap;\\
&&$a+b=2$$\iff$$ M_{1}[5(*)06](a=b=1,\ the\ lone\ 5$-cap) or M$_{1}$[50516](a=2, b=0).\\

3&$M_{2}[m6]:$&$a=2$$\iff$$ m[M_{2}655]\ or\ m[M_{2}6655]\ or\ m=6[M_{2}66655];$\\
&&$a=1$$\iff$$m[M_{2}65*]\ or\ m[M_{2}665*]\ or\ m=5[M_{2}6665]\ or\ m=6[M_{2}6665*]$\\
&&\phantom{abcde}$
\ or \ m=6[M_{2}66665] \ or\  m=5[M_{2}67055] \ or\  m=5[M_{2}670(5)5] \ or\  m=5[667050M_{2}]$\\
&&\phantom{abcde}$
 \ or\  m=5[M_{2}67065]\ or\  m=5[M_{2}671515]\ or\ m=5[M_267152(5)5]$\\
&&\phantom{abcde}$
\ or\  m=6[M_{2}670506]\ or\  m=6[M_{2}6670506].$\\

4&$M_{2}[66]:$&$a+b=2$$\iff$$M_{2}[6(*)0(*)6] (a=b=1,\ the\ lone\ 5$-$cap)\ or$\\&&\phantom{abcde}$a=2,\ b=0\ from M_{2}[m6]\ above\ or\  a=0,\ b=2\ from\ M_{2}[m6]\ above;$\\ 
&&$a+b=3$$\iff$$ M_{2}[6(5)0(*)6],\ the\ 55$-$cap;$\\
&&$a+b\le3.$\\

5&$M_{3}[555]:$&$a+b+c+d=2\iff$a 5cap on either 55, otherwise a+b+c+d=0.\\

6&$M_{3}[556]:$&
$a+b+b+d \leq4;\ a+b+c+d=4$$\iff$$
M_3[50516];$\\
&&$a+b+c+d=3$$\iff$$  M_{3} $ is the Major target vertex in discharging rule 11a, 15, or 16.\\

7&$M_{5}[656]:$&
 $a+b+c+d=6\iff$$ M_{5}[60506];$\\
&&$a+b+c+d=3 $$\iff$$M_{5}[60516];$\\
&&$a+b+c+d=2 $$\iff$$M_{5}[605(M)6];$\\
&&otherwise $M_{5}[6(*)5(*)6]$ and $a+b+c+d=0$.\\

8&$M[m_{1}6m_{2}]:$&a=2$\implies$d=0 except c=d=1 for the lone 5cap;\\
&&no 5cap $\implies a+b+c+d \leq 2$ except M[5051(6)616] with a=2, c=1.\\

9&$M_{4}[565]:$&$a+b+c+d\leq5;\ a+b+c+d=5$$\iff$$ M_{4}[5060(5)5].$\\

10&$M_{4}[566]:$
&$a+b+c+d\leq 5;$\\
&&$a+b+c+d=5\implies$ (55cap + lone 5cap, either way) or M$_{4}[5(5)0(6)616].$\\ 

11&$M_{5}[666]:$&$
a+b+c+d\leq 5;$\\
&&$a+b+c+d=5\implies$(55cap + lone 5cap) or $M_{5}[60(5)616].$\\ 

\end {tabular}
\end {lemma}
\begin{proof}
Each cross (minor,minor) edge trajectory from rules 3 to 17 crosses a final (hub, minor) edge before ending at its target Major vertex. \\

Items 1 and 2: 
Only Rule 3 has a trajectory with its final (hub,minor) crossed edge being (hub, 5), giving a 5cap with a=b=1 in item 1 for hub=5 and b=1 in item 2 for hub=6.
This 5cap establishes a$\geq$1 and a=2 only for a second adjacent 5-valent vertex upstream on the a arc, not allowed in item 1 by reducible 5[555] and giving the 55cap
 pointing away from the 6
in Item 2. A vertex 6-valent or higher upstream makes this a lone 5cap and  a=b=1 of Item 2. Otherwise, b=0 and  a+b=2 simplifies to a=2, requiring a 6-cap and two preceding 5-valent vertices, the neighbourhood given. \\

Item 3 enumerates the possible upstream neighbourhoods that lead to a=2 or a=1 with a minor hub and the last crossed edge being (hub,6)\\



Item 4, M[66]: If the cap on the 66 edge is a 5-valent vertex, then the possibilities are a=b=1 for the lone 5cap and 1$\leq$(a=3-b)$\leq$2 for the 55 cap, depending on which side has the extra 5-valent neighbour. Both sides cannot have an extra 5-valent neighbour by reducible 5[5665].
If the cap on the 66 edge is 7-valent, then a=1 only from source 17a and b=1 by a second 17a source violates the degree of the dashed-attached vertex in source 17a.

Otherwise, for positive flow to M, the cap on the 66 edge is 6-valent. Now a=1 from a neighbourhood listed in Item 2 with the shape 6[M66...] going around the left 6-valent vertex, and b=1 from a similar neighbourhood listed in item 2 going around the right 6-valent vertex. All such configurations are reducible, limiting a=b=1 to the lone 5cap.\\  
\noindent 66[5-6-5], 66[5-6-65], 66[5-6-665],
66[56-6-65], 66[56-6-665], 66[566-6-665]-13,\\
(Rule 17b:) 66[507-6-5], 66[507-6-65], 66[507-6-665]-13, 66[507-6-705]-13\\
The a=2 sources are superstructures of an a=1 source, so any a=2, b=1 structure also contains a reducible configuration limiting a=2, b=1 to the 55 cap, and a continuation for b=2 sources shows that a+b$\leq$3.\\

Item 5, M[555]: Non-zero flow on both 55 edges produces 5[5555M] containing 5[555]. \\

Item 6, M[556]: By items 1 and 2, a+b+c+d $\geq$5$\implies$c+d=3 by a 55 cap on the 5 side producing 5[555].\\ a+b+c+d=4 requires a+b=2 by a 5cap leaving c+d=2 which only occurs with a 6-cap on the 56.\\
\indent Suppose a+b+c+d=3. 3 units by a 55cap on the 56 is not possible, so again a=b=1 by a 5cap and the cap on 56 is at least 6-valent. A 6-valent vertex gives the total 4 above, so the cap on 56 must be 7-valent, i.e. neigbourhood [50526], expanding to precisely those sources listed.\\

Item 7, M[656]: a+b+c+d=6 $\implies$
a=d=1 $\implies$ 5caps for each 56 edge and M$_5$[60506].
With \{6,5\} caps we get a+b+c+d=3. None of sources 11 to 17 has a Major target with a 656 neighbour sequence, so \{M,5\} caps has a+b+c+d=2, and  \{*,*\} caps has a+b+c+d=0.\\

Item 8, M[m6m]: a=2 is from one of the five source structures listed in M[m6]. d=2 is from the same set.  d=1 is from any of the 15 source structures listed in M[m6]. Apart from d=1 due to a lone 5cap, all combinations are reducible.

Item 8, No 5cap: Without 5caps, b,c$\leq$1 and a+b,c+d$\leq$2. Assume the total exceeds 2 units. Wlog a+b=2, making the cap 6-valent, and c+d$\geq$1.  Suppose a=2; then b=0, and d=0 by the above, and the only way for c=1 is by the exception, M[5051(6)616]. Otherwise a=1, b=1, and b=1 requires m$_1$=6 and a 5cap is required by M$_2$[66]. \\

Item 9, M[565]: a+b+c+d=6 $\implies$ a+b= 3 =c+d with  55caps pointing away from the 6, i.e. reducible 6[505x505x]. For a+b+c+d=5, wlog c+d=3 and by M[56], has a 55cap pointing away from the 6 giving d=2. From M[m6m], a=0 except for the lone 5cap and by M[56], no other source has b$\geq$1.\\

Item 10 M[566]: a+b$\leq$3 by M[56] and c+d$\leq$3 by M[66], so a+b+c+d$>$5 must have c+d=3=a+b $\implies$ a 55cap on 56 pointing away from the 6 and a second 55cap on 66, giving either 6[555] or 6[505x506], both reducible.\\
a+b+c+d=5 forces a 55 cap on 56 or 66. On 56, it points away from the 6 giving a=2 from which d=c=1 from a lone 5cap, one of the possibilities or d=0 with c=2, the third possibility. On 66, a 55 cap pointing away from the 5 has d=2 giving a=b=1 from a lone 5cap, the second possibility or a=0. For a=0 we would need b=2, not possible by M[56]. If the 55 cap on 66 points upward, then the cap on 56 is Major, limiting a+b to 1 unit.\\


Item 11 M[666]: Similar to above, a+b+c+d=6 requires 55 caps on both 66 edges leading to 6[555] or, this time, 6[605x506x]. a+b+c+d=5 requires a 55 cap, say d=2, c=1, and a+b=2 from a lone 5cap is possible. Otherwise, a=0 and b$\leq$1. Alternately, the 55cap could point upward with d=1, c=2, and a 6-cap on the other 66 edge gives the indicated structure with the required flow. A 5cap would lead to 6[555] and a 7 or higher cap cannot deliver the required 2 units.
\end{proof}

\begin{figure}[ht] 
  \begin{center}
    \includegraphics[
    width=6in
    ]{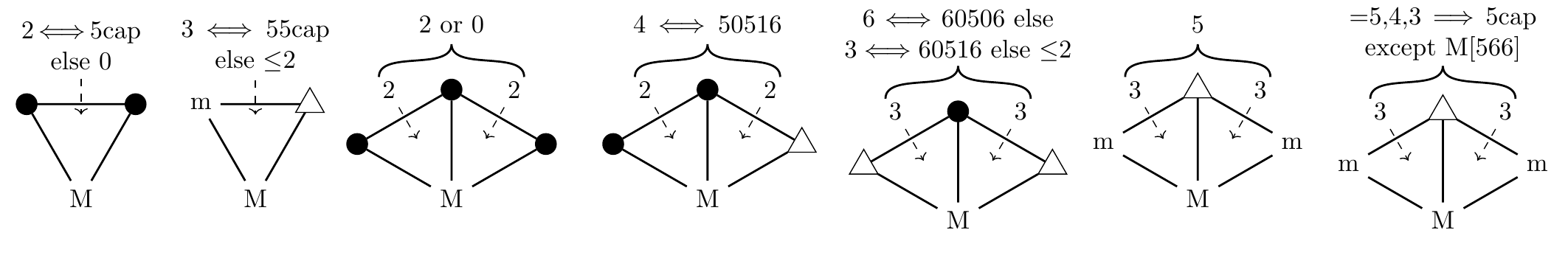}
  \end{center}
  \caption{Flow maximums across (m,m) edges or (m,m) edge pairs}
  \label{fig-Cross mm edge summary}
\end{figure}

The main results are summarized in Figure ~\ref{fig-Cross mm edge summary}.
\begin{samepage}
\begin{corollary}Across the four rim edges of M[55655], at most 6 units is delivered to M.
\label{lem-55655}
\end{corollary}
\begin{proof} From items 1 and 6, at most 2 units crosses each of the 4 edges to M. To exceed 6 units, each (5,5) must have a 5cap, forbidding a 5cap on the (5,6) edges, and requiring at least 3 units  to cross the 565 pair of edges, violating the second part of Item 8\_M[m6m].
\end{proof}
\end{samepage}

\section {Vertices of degree 10 are not overcharged}
\begin{remark} \label {Rmk10}
 Individual and families of Reducible configurations used in this section.
 
\begin{enumerate}
\item \label{two60506s} 
10[56050660506xxx]-13, 10[60506$\text5^k$60506]-13, k=1,2,3.
\item \label{xxx=mm6,m,60506} 
10[5\{5,m,60506\}5xxx]-13,10[5\{5,5,m\}60506xxx]-13,
10[5\{5,5,m,60506\}5xx]-13,10[5\{5,5,5,m\}60506xx]-13.
\item \label{xarticulatedx} 
10[x555x55555]-12,   10[x5555x5555]-12. 
\item \label{articulated60506}10[x555x\{5,5,60506\}]-13,  10[x5555x560506]-13, 10[55555x60506x]-13.
\item \label{xxx,*forcedM} 
10[5\{5,5,5,5,m\}5xxx]-12.
\item \label {10[555555,7]} 
10[5\{5,5,5,5,705\}xxx]-13, 10[5\{5,5,5,71505\}xxx]-13.
\item \label{xx,**not66}  
10[5\{5,5,5,5,6,6\}5xx]-13.
\item \label{articulated 505}
 10[505x5\{5,5,5,m\}5x]-13,  10[505x\{5,5,5,60506\}x ]-13.
\item \label{-14ring} 
10[ 5\{5,5,5,5,6,6,6\}5x]-14.

\end{enumerate}%
\end{remark}
\begin{theorem}\label{10notOver}
Vertices of degree 10 are not overcharged.
\end{theorem}
\begin{proof}
Recall the Radial inflows and their limits: 6 units for a 5-valent vertex iff it is part of a 60506 cluster, otherwise at most 5 units, 3 units for a 6-valent vertex, and 2 units for a Major vertex. (Figure ~\ref{fig-Max-Radial-Inflows}.) As with 11-valent vertices, the first subcasing is by the number of 5-valent vertices.
\begin{itemize}[leftmargin=0.25in]
\item 
A first neighbourhood with $k$ 5-valent vetices and $(10-k)$ vertices of degree at least 6, i.e. 
10[\{$\text5^k$,$\ast^{(10-k)}$\}], contains at most three non-overlapping 60506 neighbour clusters.\\ The final charge is at most
 $-40+3\cdot6+(k-3)\cdot5+(10-k)\cdot3=2k-7\le0\text{ for  }k \le3$
\item 
For k=4, the upper bound on the final charge is +1 and overcharging would occur only for all radial inflows at their maximums: three 60506 clusters and a fourth 5-valent vertex, but that contains reducible 10[56050660506xxx]-13 from group ~\ref{two60506s}.
\item 
10[\{$\text5^5$,$\ast^5$\}] has at most two 60506 clusters and an upper bound on its total radial inflow is 
2$\cdot$6+3$\cdot$5+5$\cdot$3=42, so an overcharging neighbourhood must have at least one 60506 cluster and at least four of the $\ast$'s must be 6's.
	\begin{itemize}[leftmargin=0.25in] 
	\item 
	Two 60506 clusters: The neighbourhood is 
 10[\{60506, 60506, 5, 5, 5, $\ast$\}].  If the 60506 clusters are adjacent, then there is a 5 at one end or the other yielding the reducible configuration above. Otherwise, they are separated by one, two or three 5-valent vertices, one of the other reducibles in Remark ~\ref{Rmk10} group ~\ref{two60506s}.
 	\item
One 60506 cluster: The maximum radial inflow is 41 and all limits must be at their maximum to overcharge the 10-valent vertex. In particular, the neighbourhood must be from 10[\{60506,5,5,5,5,6,6,6\}].
These 19 configurations are all bounded by 14-ring, and likely reducible, possibly with reducible subconfigurations bounded by a 15 circuit. 

Instead, examine these neighbourhoods more closely with the aim identifying those that might be overcharged and showing that these contain reducible configurations bounded a separating circuit of size at most 13.
		\begin {itemize}[leftmargin=0.25in]
		\item 
The three 6's not part of the 60506 could be consecutive, or a pair could be consecutive and the third "loose". Remark ~\ref{Rmk10} group 	~\ref{xxx=mm6,m,60506}
covers these possibilities.
		\item
Otherwise, the three unattached 6's are singletons and separate the five items: \{5,5,5,5,60506\} into groups of sizes 3-1-1- or 
 2-2-1- with the three singleton 6's in the dashed positions. 
If any  1  is a singleton 5 and therefore part of 656, then the rest of the neighbourhood begins with 5 and ends with 5 or 60506, with one singleton 6 and the remaining two items.  Again, remark ~\ref{Rmk10} group ~\ref{xxx=mm6,m,60506} covers these possibilities.

 The final case is 2-2-1- where the size 1 group is the 60506 and is bounded on both sides by a six and then a pair of 5's and the final 6:  10[556605066556]-14  (See Figure ~\ref{fig-10Targets}, row 1, \#2). This time, evaluate the total inflow to the 10-valent vertex by Mayer's 2,3,4 and the cross rim edge flows.  The 4 unit contributions from the five 5-valent vertices total 20 units. 
The two 56 edges of the 60506 cluster each contribute 3 units  and their adjacent 66 edges contribute at most 2 more each. The 56 edges before and after them are part of a 556 chain so deliver at most 2 units each, a total so far of 20+14 units. The last 4 edges are the wrap-around cluster, 10[55655], delivering at most 6 more units by Lemma ~\ref{lem-55655} and the 10-valent vertex cannot be overcharged. 
None of the nineteen 14-ring configurations need to be tested for reducibility.
 		\end{itemize} 
	\end{itemize} 
\begin{figure}[ht] 
  \begin{center}
    \includegraphics[
    width=5in
    ]{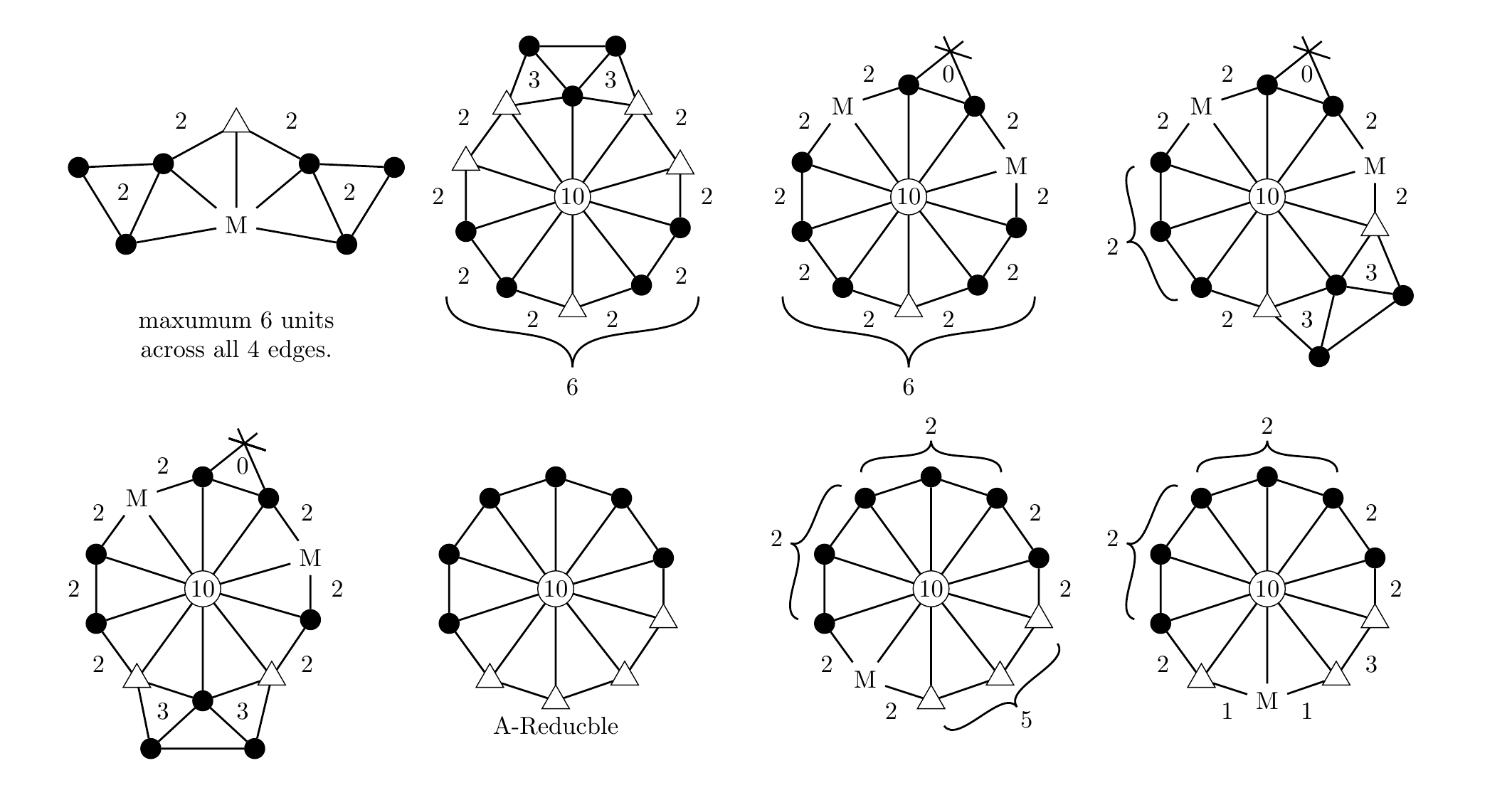}
  \end{center}
  \caption{Lemma M[55655] and some 10-vertex neighbourhoods.}
  \label{fig-10Targets}
\end{figure}
\noindent The remaining cases will be handled from ten 5-valent neighbours down to six.
\item
10[ $\text5^\text{10}$] and 10[ $\text5^\text{9}\ast$] contain reducible 10[$\text5^\text7$xxx]-11 from group ~\ref{xxx,*forcedM}  which also eliminates neighbourhoods from\\10[\{$\text5^\text8$,$\ast^\text2$\}]
when the $\ast$'s are adjacent or one 5 apart. When the two $\ast$'s are three or four 5's apart, the neighbourhood is from 
Remark ~\ref{Rmk10} groups ~\ref{xarticulatedx}.
This leaves 10[*55*555555] in which case avoiding 10[5565555xxx] of 
Remark ~\ref{Rmk10} group ~\ref{xxx,*forcedM}
forces the * s to be Major.  Since there are no [656] neighbour runs, the final charge is at most -40 +8$\cdot$5+2$\cdot$2=+4. 

Recalling Figure ~\ref{fig-CrossMm} and peeking ahead to Figures ~\ref{fig-AllSources2} and ~\ref{fig-AllSources3}, the only cross (Major,minor)-edge discharging rules that have \{5,5\} shoulders and a trajectory using the minor vertex as a hub are Rules 21 and 22. Avoiding these structures will force such a Major vertex to have a radial inflow of zero units.

Reducible configurations from Remark ~\ref{Rmk10}  group ~\ref{10[555555,7]} accomplish this task and the final charge is at most -40 +8$\cdot$5+2$\cdot$0=0.

\item
	 10[\{$\text5^7$,$\ast^3$\}]: If the three $\ast$'s are consecutive, the remainder is 10[$\text5^7$]. If they are grouped 2-1, then the singleton $\ast$ is surrounded by six 5's (actually seven) and is Major by group ~\ref{xxx,*forcedM} and delivers 0 units of radial inflow by group ~\ref{10[555555,7]}.
If the pair of $\ast$'s is 66, the neighbourhood is one of 10[5\{5,5,5,5,66\}55M]-13, each of which is a superconfiguration of one of the reducibles in Remark ~\ref{Rmk10} group ~\ref{xx,**not66}. 
With one $\ast$, one unrestricted M and the flow from the other M limited to zero, the total radial inflow is at most 7$\cdot$5 +1$\cdot$3 +1$\cdot$2 +1$\cdot$0 = 40.

Otherwise, the three $\ast$'s are singletons and separate the seven 5's into three clusters of sizes 5-1-1-, 4-2-1-, 3-3-1-, or 3-2-2-. For the first three, the $\ast$ in the first spot is again surrounded by six 5's and must be Major with zero radial inflow, same as above. The other two $\ast$'s are part of 10[5\{5,5,5,5,5,$\ast$,$\ast$\}5M] and cannot be 6,6 by Remark ~\ref{Rmk10} groups ~\ref{xx,**not66} and ~\ref{xxx,*forcedM}. 
 Again, the total radial inflow is at most 7$\cdot$5 +1$\cdot$3 +1$\cdot$2 +1$\cdot$0 = 40.

The 3-2-2- case is 10[555$\ast$55$\ast$55$\ast$]where again no pair of $\ast$'s can be 6,6, so they are $\ast$,M,M and there are three possible patterns of first neighbours:
10[555655M55M], 10[555M55655M] and 10[555M55M55M]. The radial limits give possible inflows of 42, 42, and 41 units so consider the Mayer rules and cross rim edge contributions instead.

The 2,3,4 contributions from 5-valent vertices are 24, 24, and 22. The 5M5 contributions are at most 8, 8 and 12 for subtotals of 32, 32 and 34. For the third case, the cross 55 and 555 limits are 2+2+2, so for that pattern, the 10-valent vetex is not overcharged. The second case has the cluster 55655 with the maximum of not 8 but 6 units shown earlier in Lemma ~\ref{lem-55655}, and an extra 2 units across  555 for a total of 8 units, so it is not overcharged. 

For the first one, (see Figure ~\ref{fig-10Targets}, row 1, \#3) we add reducible 10[555655x505x]-13 from group ~\ref{articulated 505} to limit the cross edge flow to zero on the 55 edge between the Major vertices.  Again, the 55655 group has a limit of 6 units and the remaining 55 edge may contribute 2 more units, for a total of 8 units and this neighbourhood also cannot overcharge the central 10-valent vertex.
\item 
 10[\{$\text5^6$,$\ast^4$\}] This is the most complicated case.
Several 14-ring reducible configurations are added to the set to be avoided. The maximum via radial inflows is 2$\cdot$6+4$\cdot$5+4$\cdot$3=44, so there are many subcases and much work is needed  to show that in every case, either there is  a reducible configuration or the inflow is at most 40 units.
	\begin{itemize} 
	\item Two 60506 structures: Each neighbourhood is from 10[\{5,5,5,5,60506,60506\}] and contains a reducible configuration 
from group ~\ref{two60506s}
	\item One 60506 structure: The neighbourhood is from 10[\{5,5,5,5,5,60506,$\ast$,$\ast$\}].  Reducible configurations from
group ~\ref{xxx=mm6,m,60506} eliminate the cases with the $\ast$'s consecutive or one 5 apart.  The $\ast$' s three, four or five 5's apart are eliminated by group ~\ref{articulated60506}. This leaves the * 's two 5's apart and Major by  group ~\ref{xxx=mm6,m,60506}.
The 55 between these Majors has zero flow across it by group ~\ref{articulated 505}, so the only neighbourhoods without a reducible configuration are 10[M5(*)5M55560506]   and 10[M5(*)5M55605065]. (Figure ~\ref{fig-10Targets}, row 1, \#4 and row 2, \#1).
The 2,3,4  contributions from 5-valent vertices are 21 and 20 respectively, and the cross edge maximum totals are 18 and 20 respectively and these 10-valent vertices are not overcharged.

	\item No 60506 structure: The Maximum radial inflow is 6$\cdot$5+4$\cdot$3=42. For overflow, at least three of the * 's must be 6-valent, so the neighbourhood is 10[\{5,5,5,5,5,5,6,6,6,*\}].	
The \{6,6,6,*\} could be grouped 4-0-, 3-1-, 2-2-, 2-1-1-, or 1-1-1-1- with the dashes holding one or more 5's. 
For 3-1-, the lone entry has all six 5-valent vertices around it and must be Major with zero units radial inflow and a 10-valent vertex with this neighbourhood is not overcharged.  For 2-2-, 
at least one pair is 66 and the neighbourhood is reducible from group ~\ref{xx,**not66}.

For 2-1-1-, if the 2 contains 6*, then the neighbourhood is 
10[5\{5,5,5,5,6,6\}56*], containing a reducible configuration from group ~\ref{xx,**not66}.
Otherwise the * has 5-valent vertices on both sides, as does the * in 
1-1-1-1-,  i.e. the neighbourhood is one of nineteen in the set 10[ 5\{5,5,5,5,6,6,6\}5*]-14, included in Remark ~\ref{Rmk10}, group ~\ref {-14ring}. 

 For 4-0, there are three subcases: 10[5555556666], 10[555555666M], and 10[55555566M6] (Figure ~\ref{fig-10Targets}). The first is A- or directly reducible
 ~\cite{Franklin}, the middle neighbourhood has a Mayer inflow of 23 and across the rim edges, a maximum of 17 more. For the last case,
note that in Figure ~\ref{fig-AllSources2} none of the minor hub trajectories (rules 21 to 28) has \{6,6\} shoulders, so the maximum cross Major-minor edge flow is 1 unit for the M6 edges of 6M6. The 
 total inflow is 6x4=24 plus a possible 15.
 	\end{itemize} 
\end{itemize}
\end{proof}

\section{The Discharging Rules - Part 2: Transfers across (Major,minor) edges}
\begin{figure}[ht]
    \begin{center}
    \includegraphics[
   width=6in
    ]{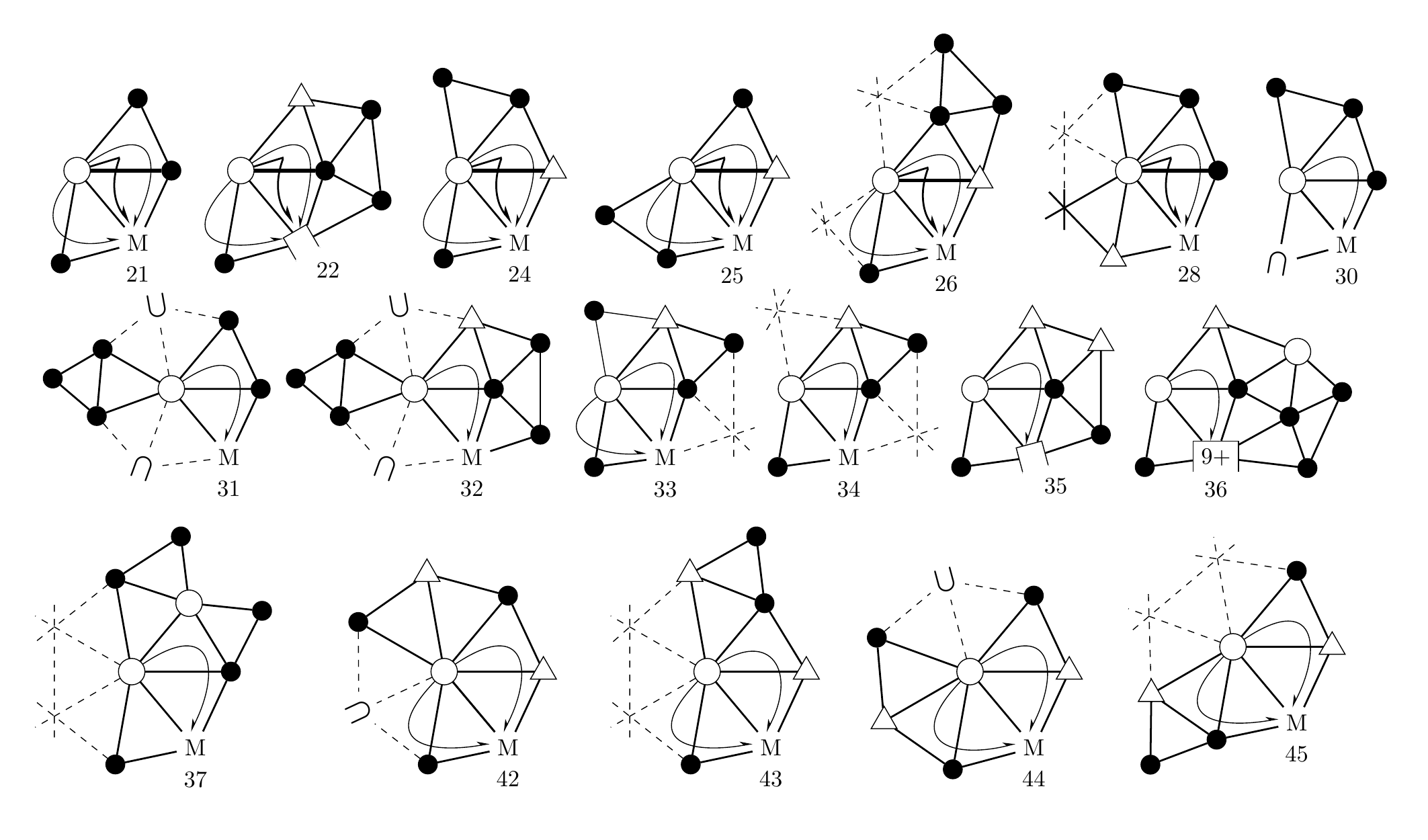}
  \end{center}
  \caption{Discharging Rules 21-45 }
  \label{fig-AllSources2}
\end{figure}

\begin{figure}[ht]
    \begin{center}
    \includegraphics[
   width=6in 
   ]{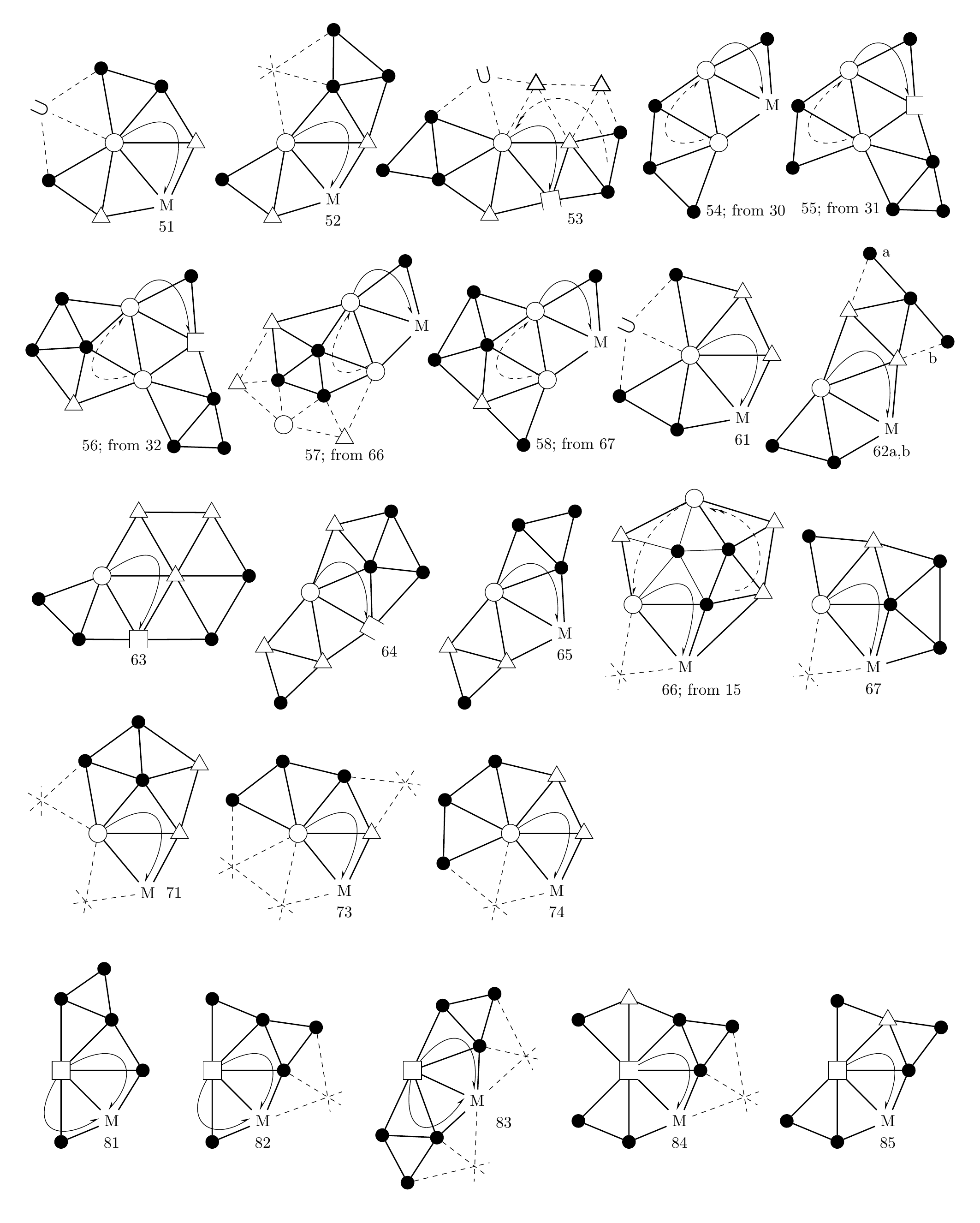}
  \end{center}
  \caption{Discharging rules 51 to 85}
  \label{fig-AllSources3}
\end{figure}


The remaining discharging rules are numbered for reference and listed in Figures ~\ref{fig-AllSources2} and ~\ref{fig-AllSources3}.

Rules 21 to 28 describe the six  
neighbourhoods of a (7,m) edge which support a transfer of 2 units from the 7-valent source across the (7,m) edge (in bold)  to the Major sink adjacent to the 7-source . All other (M,m) edges support a transfer of at most one unit, that unit is with the Major vertex as the source and hub.
\begin{remark}
These discharging rules are not additive or summable, as in RSST, but rather ORable: more than one rule can "activate" the same trajectory for the same one unit of flow. For example 7[55M55] has two occurrences of Rule 21, 7[55M5], which generates 3 units of flow by itself, but the total flow from the 7 to M in 7[55M55] is 4, not 6 units. It is the vertices attached by solid edges that define the discharging rule and activates the trajectories.
\end{remark}

For a cross (Major,minor)-edge transfer rule to apply, the target vertex must be Major. This is not a consequence of the minor maxim and reducible configurations, as was the case in the cross (m,m)-edge rules, but simply an a priori neighbourhood requirement for the discharging rule to be activated. 
However, whenever the degree of a target vertex is listed as greater than 7, that condition is a consequence of the configuration with a 7- and sometimes 8-valent target being reducible. (Rules 22, 35, 36, 53, 55, 56, 63, and 64). Similar to this, several vertices drawn with dashed edges have a minimum degree that is a consequence of a lesser degree vertex leading to a reducible configuration or subconfiguration, as in Figure ~\ref{fig-Rules11to17}.
In particular, if the structure given by the solid edges is in a m5CPG, these vertices attached by dashed edges will have the indicated minimum degrees. 

Still more vertices joined by dashed edges
 have a specified minimum degree. In these cases, if the vertex was of a lesser degree, then the structure could appear in a m5CPG, but it would contain an additional source structure with at least the same flow trajectories. For example,
Rule 43 would contain Rules 25 or 42, and Rule 67 would contain Rule 22. 
Another collection is Rules 22, 33 and 34. 
In particular, if the source structure given by the solid edges occurs in a m5CPG, the outflow from the source vertex to the target vertex is guaranteed. Because the rules are ORable, nothing is lost if the rule is limited to those occurrences with the dashed attached vertices having the indicated minimum degrees. 

Many rules are a consequence of a previous rule possibly overcharging its target of degree 7 which becomes the source vertex in such a rule. This previous rule is indicated by the dashed arcs in Rules 53,  54 to 58 and 66. 
In Rules 53 and 57, the vertices connected exclusively by dashed edges show the complete neighbourhood of that previous rule, but are not necessary for the discharging rule given by the solid edges to be activated. In fact, these are the only neighbourhoods where these discharging rules are needed. Those extra vertices and edges could be part of the rule and any discharging claims would still be true since these rules are not used anywhere else. The smaller structure is needed to satisfy a desirable feature described later.

Rule 62 requires a 55cap on the (6,6) edge (and its transfer of 3 units to the 7-valent  source). One of the "dangling"  5-valent vertices with one solid and one dashed edge is required, and, if necessary, the rules can be distinguished as 62a and 62b. 

\subsection {Features of Rules 28 and 30.}
\label{ssec-Rule28}
The pair of Rules 28 and 30 specify the transfers from 7 to M in 7[555M*x]. The flow is 1 unit from 7 to M iff * is Major (Rule 30), and  2 units iff * is 6-valent AND x is at least 6-valent. In Rule 28, this vertex x is drawn with solid edges indicating that this is a requirement of the rule. Other rules apply for 7[55M5] and 7[55M65]. These requirements on Rule 28 avoid the generation of too much outflow due to the ORable nature of overlapping rules. In particular, 
there is no additional charge transfer across the other shoulder, the (7,6) edge, to M in the structure of Rule 28.
Any flow across this (7,6) edge would need to come from a rule shaped like 7[*6M555] and only rules 62 and 63 support this pattern. 
The combinations result in reducible configurations 7[555x606x]-11 and 7[555x61(5)6x]-12.

\subsection {Transfers of 3 or more units from a Major source to the same Major sink.} 
Because these transfers are ORable, the maximum flow across a (M,m) edge is 2 units, and any cluster transferring 3 or more units from a Major source to the same Major target must use Rules 21 to 28. For Rules 24 to 26, a fourth unit would require \{5,6\} shoulders i.e. Rule 28, but that rule requires \{5,6*\} shoulders.

Only three clusters deliver a total of 4 units from a Major source to the same Major sink: 7[55M55] mentioned earlier, Rule 21 overlaid on Rule 22, and Rule 22 overlaid on Rule 22 for which the target will have degree at least 9 because 8[5571505]-11 is reducible.

A transfer of exactly 3 units from a Major source to the same adjacent Major target occurs only from Rules 24 to 26 and 21 and 22 with no overlay.

\subsection {The grouping of Cross (M,m)-edge transfers, Rules 21 to 85 }
\begin{itemize}
\item Source structures with an 8-valent source: Rules 81 to 85.\\ Otherwise the source is 7-valent.
\item Source structures with a 2-unit flow across a (M,m)-edge: Rules 21 to 28.
\item Source structures with \{55,M\} shoulders: Rules 30 and 31 . 
\item Source structures with \{65,M\} shoulders: Rule 32 . 
\item Other source structures with \{5,5\} shoulders: Rules 33 to 37.
\item A 65 cap added on either (7,5) edge of 7[56T5], yields a 1+1 unit source: Rules 42 to 45.
\item Source structures with \{6,65\} shoulders: rules 51, 52 and 53. 
\item Source structures with \{5,7\} shoulders: rules 54 to 58. See Lemma ~\ref{lem-7[5MM5]}
\item Source structures with \{66,5\} shoulders: 61 to 63 and 64, 65. 
\item Source structures with \{5,*\} shoulders:  66 and 67.
\item Source structures with \{6,*\} shoulders: 71 to 74.
\end {itemize}

\section {Counterflow and Limited Flow Targets}

\begin{figure}[ht] 
  \begin{center}
    \includegraphics[
    width=6in
    ]{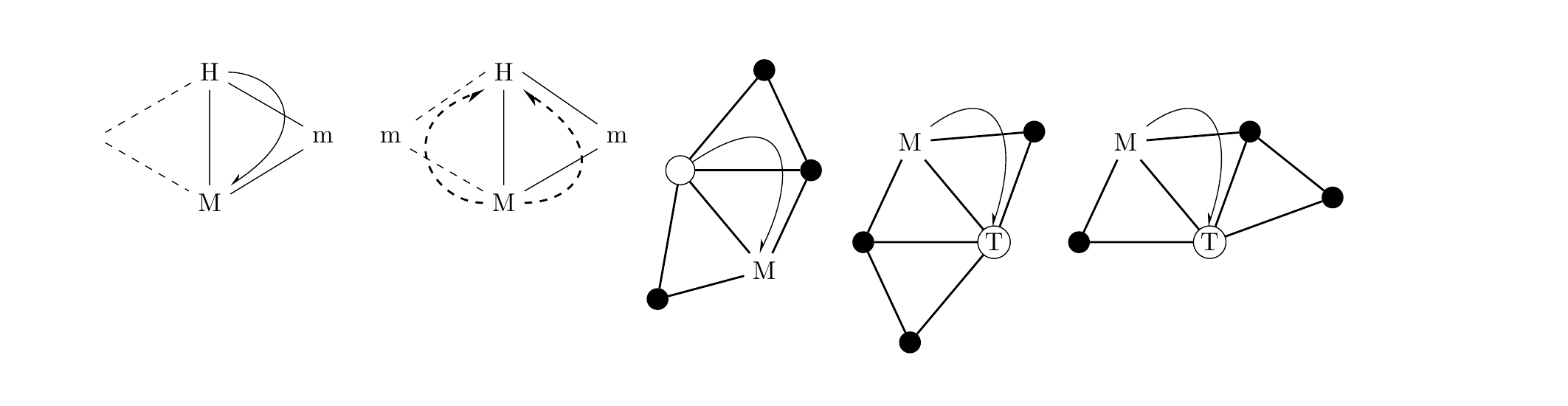}
  \end{center}
  \caption{Flows from H to M and 7[55M5] and their counterflows}
  \label{fig-Counterflow}
\end{figure}

\subsection{Counterflow for flow across (Major, minor) edges} Similar to flow across (minor, minor) edges, we now show that when there is flow from a  Major source vertex H across a (H, minor shoulder vertex m) edge to  a target Major vertex M,  there is no flow from the target vertex M back to the source vertex H. 
Since the minor hub unit of a cross-Mm flow is always in conjunction with a unit with the Major vertex as hub, it suffices to consider only the Major hub trajectories here. 



Each shoulder vertex of the (H,M) edge gives a possible trajectory for a counter flow from M to H (Figure ~\ref{fig-Counterflow}). As long as neither of these flows occurs, then the full specified flow from H to M is a guaranteed net outflow. 
For example, rotate rule 21, 7[55M5] 180 degrees around the midpoint of the (7,M) edge and the source cluster becomes a target cluster T=7[55\xedge{M5}] 
where the 7-valent vertex is the target, the crossing edge is overlined,  and M is now the Major source and hub of a possible counterflow. Mirror that along the (M,7) edge to get  the other possible counterflow target cluster, T=7[5\xedge{M5}5]. As long as each source rule cluster overlaid on these target structures,
aligned  at the source, target, and flow crossing edge, and 
compatible everywhere, produces a configuration that
is or contains a reducible configuration, the counterflow will not occur. The  full 3 (or 4) units of flow is  then guaranteed as the net charge transfer from the 7-valent vertex to M in rule 21, 7[55M5].


A target structure such that every source structure, overlaid and consistent with that target structure produces a configuration that is or contains a reducible configuration is called a zero flow target (ZFT). The supporting set of reducible configurations is represented by this target structure.
   T=7[55\xedge{M5}] and T=7[5\xedge{M5}5] are zero flow targets. While 7[56T5] is not a source, 
 T=7[56\xedge{M5}] and T=7[5\xedge{M6}5], are zero flow targets, assuring that sources 24 to 26 and 42 to 45 have their full specified  outflows and no counterflow. \\
 

Every source structure converts, as in Figure ~\ref{fig-Counterflow}, into one or two  counterflow targets. Each of these targets is a ZFT, like those from Source 21, or contains a ZFT, like the counterflow structures from Sources 24-26 and 42-45. These and other ZFTs for 7-valent targets are listed in Appendix A.  8- and 9-valent targets are described in appendices C and E.

\begin{remark} In RSST \cite{RSST}, their rules 7[5mT] and 7[566T] do not have this net outflow property. 77[5-6-5]-11 and 77[5-6-65]-12 are not reducible, and 77[56-6-65]-13 is SDI, not reducible by simple Kempe chaining.
 This is where their discharging and this one diverge irreconcilably.
  Also, as in  RSST, this proof does not invoke "block count reducibility", so any SDI configuration that is not A-reducible may appear in the presumptive m5CPG.
\end{remark}

\subsection{ Generalized Zero Flow Targets} 
The main use of a zero flow target is to establish a limit on the inflow across a rim edge when determining the final charge on
 the central vertex of a 7-, 8- or 9- wheel with a given neighbourhood.
Two generalizations are possible here. First, the crossing edge could also be (minorHub, minor), and secondly, for a non-symmetric target cluster, the Hub and non-hub could exchange functions. 
The previous target structures generalize to T=7[5\xedge{xx}5]  and T=7[5m\xedge{xx}]  and are zero flow targets and  each target represents the set of reducible configurations required to meet the definition.
The degrees of the end vertices of the crossing edge are now unspecified. All the discharging rules are considered and each compatible source structure must result in a reducible configuration assuring no flow crosses the edge and continues to the target. If both ends of the crossing edge are Major, no rule applies. If a zero flow target structure (ZFT) appears in a m5CPG, there is no flow from this discharging (Rules 1 to 85) across the crossing edge to the Major target vertex.


The overscore  was necessary to distinguish the vertices of the crossing edge from other rim vertices of a Major wheel.
Instead we now use L and R to identify the left and right vertices of the crossing edge. 
The degree of T is always specified, but the degrees of L and R are usually unrestricted but occasionally limited.

Every source structure except Rule 28 and 30 has the property that apart from its Major Target vertex, the structure of vertices of specified degrees (attached by solid edges) is Geographically Good (GG: each vertex adjacent to at most 3 consecutive  or 2 non-consecutive boundary vertices, no hanging 55 pair or other 2-splice impediment to reducibility, see ~\cite{FAPart1}). For Rules 28 and 30, the property holds for their common substructure, 7[555T]. In the target structures above and to come, the structure of vertices of specified degrees together with the vertices of the crossing edge,
is GG apart from the vertices of that crossing edge. 
With these caveats, source and target structures will nevertheless be called GG. 
When a GG source structure is overlaid and combined with a GG target structure, the resulting configuration will be GG and hence likely to be reducible. At the least, it will have no obvious impediment to reducibility.

\begin{remark}  H=7[56T5]  is GG, and with many targets, produces a reducible configuration. However, 
if it was a source then the very useful T=7[55\xedge{xx}] or T=7[55LR], would not be a zero flow target because 7[56705]-11 is SDI, so in this proof, it is not forbidden from occurring in a m5CPG. Similarly, if H=7[506T5] was a source, then some other very useful targets are not zero flow. 
The extra 5, 55 or 65 caps in Rules 24 to 26 and 42 to 45  are sufficient to avoid invalidating the targets used in this discharging. 
If this general approach to proving 4CT, i.e. discharging rules and ZFTs,  is used, then other sets of discharging rules would lead to different sets of ZFTs. The choice basically comes down to how many exceptions the designer wants to cope with because these must all be handled as special cases, as was the 60606 structure in the radial flow limits. 
\end{remark}

\label{ssec-T7impC12}
A second feature of the source structures is that
 if the target vertex is a lone 7-valent vertex, then the structure is bounded by a separating circuit of at most 12 vertices. 
When combined with the target structures T=7[5LR5] or T=7[55LR], with only 5caps on this target 7-valent vertex, the resulting GG configuration is again bounded by a separating circuit of at most 12 vertices and can be checked for reducibility by a software program that handles candidate configurations bounded by at most a 12-ring. Target structures T=7[56LR] and T=7[565LR] combined with a source structure will produce a candidate GG configuration bounded by at most a 13-ring. Showing T=7[566LR] is a zero flow target requires reducibility software able to handle configurations bounded by a 14-ring, and T=7[5666LR]  requires software able to handle configurations bounded by a 15-ring.

Target structures used in this discharging are listed in Appendices ~\ref{app-ZF7}, ~\ref{app-ZF8} and  ~\ref{app-ZF9}. Target 7402, pictured as 7[5mLR], refers to both ZFTs T=7[55LR] and T=7[56LR] and represents all the configurations composed from these targets combined with each compatible source structure. The zero pointing across the LR edge indicates that every compatible combination is reducible. Target 7408 represents the set of configurations from all compatible sources with T=7[605LR] and T=7[606LR]. All these configurations are reducible except one, 7[6066(*)0(*)6]. 

Targets 7411, 7413 and several others have a *cap on the (R,T) edge. If that cap were of degree 5, 
then the target would contain target 7401 and that would suffice to limit the flow across the LR edge to zero. Thus, we can limit the use of  target 7411 et al.  to cases with that cap of degree at least 6. In 7435, if the *cap is instead 5-valent, then ZFT 7406 shows the flow is zero. In 7545, this applies only when m=5.
 Space did not permit recording this conditional * in 7402, 7503, 7604  but it was applied when those targets with m=5 were matched with all possible sources. 
Similarly, in the Articulated targets A1 through A9, the cap vertex on RT opposite L  or on LT opposite R can often be assumed to be Major, because otherwise the limit is imposed by these "smaller" ZFTs: 7402, 7503,  7509 and 7545. 
Requiring these vertices to be at least 6-valent or Major merely avoids testing  superconfigurations of an otherwise confirmed reducible configuration.


 \begin{remark}When  Part 1 - Reducibility ~\cite{FAPart1} was produced, this would have required serious computing resources, but today's computers like my 2013 MacBook Pro have the speed and memory to determine D-reducibility or confirm Symmetric D-Irreducibility of even 16-ring configurations in a matter of seconds. E-reducibility of 15-ring configurations likewise takes seconds while C-reducibility instead takes human interaction to suggest a seed and a separate program which can confirm it is a reducer, specify how it leads to a reducer, or, failing those results, supply information on what may or may not lead to a reducer. More human interaction, skill, and luck leads to an alternate or improved seed and eventually a reducer. In particular, I have yet to find a counterexample to my Conjecture 2  of Part 1, that every Asymmetrically D-Irreducible configuration has a reducer and hence is reducible.
\end{remark} 
 
 A personal challenge of this discharging is that it should require the computed reducibility of as few as possible configurations bounded by a 16-ring or larger. In the development of this set of discharging rules together with the set of  targets which would be sufficient to show overcharging a major vertex is avoided, the choices of rules and targets was influenced by this desire to produce and require the reducibility of only geographically good configurations bounded by at most a 15-ring. Relaxing this self-imposed limitation may lead to a simpler proof, just as adding more 14- and possibly 15- ring reducible configurations could simplify the demonstration that 10-valent vertices are not overcharged.
To this end, all zero flow targets must have a "delta" of at most 3, i.e. adding a target to a source must increase the size of the surrounding separating circuit by at most 3. Around the target, 5-valent vertices are "free", 6-valent vertices have a "cost" of one, and the target vertex of degree k "costs" k-7. Any articulating vertex imposes and extra "cost" of one to the size of the surrounding separating circuit.
 To stay within this 15-ring limit, a zero flow 8-target simple structure can have two 6's or a 7-valent addition, and a zero flow 9-target simple structure can have at most one 6-valent addition. See Appendices ~\ref{app-ZF8} and ~\ref{app-ZF9}.
 
For Rules 53 and 57, keeping the extra vertices of these sources with a target vertex of degree 7 would result in a boundary circuit of size 13 or 14. Dropping the dashed vertices reduces this to 12 while maintaining the GG property and these smaller sources do not invalidate the claims made by the targets.

Configurations that are D- or E-reducible do not need a reducer in the classic sense, their reducers are a relaxation of the original configuration and do not introduce loops. Typically, a classic reducer joins or merges two boundary vertices of the separating circuit. Such an edge or merger may produce a loop making the result not colourable so care must be taken in choosing a reducer to ensure that this does not happen. The almost 6-connected property of a m5CPG provides the required assurance for vertex split relaxations and for the reducers of C1- and  C2-reducible configurations ~\cite{FAPart1}. A second but minor consideration in the selection of source rules and target structures is my desire to limit the number of reducible configurations that need a custom reducer, i.e. are not D-, E-,  C1- or C2-reducible, mainly because my software, described in Part I, determines these reducibility conditions directly. To this end, some sources and targets may be larger than necessary for reducibility. 
For example, using the smaller underlying sources described with Sources 11 to 17 does not impede their discharging of the source 7-valent vertex, but when such a smaller source is combined with a target cluster, the combination is less likely to be D- or E- reducible, or, in a worst case, may be a non-reducible configuration,  leading to the failure of the cluster as a zero flow target and undesirable consequences.

\subsection{ Not Quite Zero Flow Targets} 

For a target structure with a "delta" of +4, which would require reducible configurations bounded by a 16-ring, 
one way to meet this preference 
is to limit the source structures  to those that deliver 2 or more units across the crossing edge: Rules 21 to 28 and structures from Lemma ~\ref{lem-cross-mm} items 1, 2, 3 and 4 where a+b is at least 2 units. In every case, when the target is 7-valent, these sources produce a structure bounded by at most an 11-ring. 
Combined with these sources, a target with a delta of +4 will yield a configuration bounded by a separating circuit of at most 15 vertices. 
When  each of these configurations is or contains a reducible configuration, it is called a 1-flow target, and if it 
 appears in a m5CPG,  there is a maximum flow of 1 unit across the crossing edge to the target vertex.
 For example, in Appendix ~\ref{app-ZF9}, targets 9721 to 9726, with simple first neighbourhood structures and only one 6-valent vertex are ZFTs, while targets 9830 to 9859, simple first neighbourhood structures with two 6-valent vertices are 1-flow targets.
The 1 pointing across the LR edge indicates that this target represents only those configurations formed by using a 2 or 3 unit source across the LR edge and that all these configurations are reducible. 

1-flow targets are rarer in Appendices ~\ref{app-ZF7} and ~\ref{app-ZF8} where \_1 is also appended to the target's label for example, 7A5\_1 and 8532\_1.
In Appendix ~\ref{app-7Nhbds}  where 7-wheels  are examined for possible flow across rim edges, this is turned around as 1A5 to indicate that the inflow across such an edge is no more than 1 unit due to 1-flow Target 7A5\_1. 
Target 7A9 has a "delta" of +4 when mm=66, but for that case, it is not a ZFT, but rather a 1-flow target. Similar to source 62, a 55 cap on 66 and its 3 units of flow is required for these 1-flow fully specified targets.

Another problem occurs when a target that would be very useful as a 0- or 1-flow target has a small number of sources with which it produces a configuration that is not reducible. One example of this is the Articulated target T=7[MLRM565] (7A1 in Appendix \ref{app-ZF7}) which is reducible when combined with every source except L=6[5R=6T], the lone 5cap on LR=66,  in which case it produces the SDI configuration 7[M606M565]-12. The examination of the 7-wheel neighbourhood  7[*xx*565] must entertain this possibility of 2 units across the (x,x)=(6,6) edge due to this lone 5cap but that is the only source that can produce a positive transfer across that (x,x) edge to the 7-valent vertex. (See 7[M565M**] in Appendix ~\ref{app-7Nhbds}). 
Such exceptions are explicitly noted with an otherwise 0- or 1-flow target.


\section{
Analysis Support for, 
and Enumeration of, Neighbourhoods of a 7-Wheel.
}

\begin{figure}[ht] 
  \begin{center}
    \includegraphics[
    width=6in
    ]{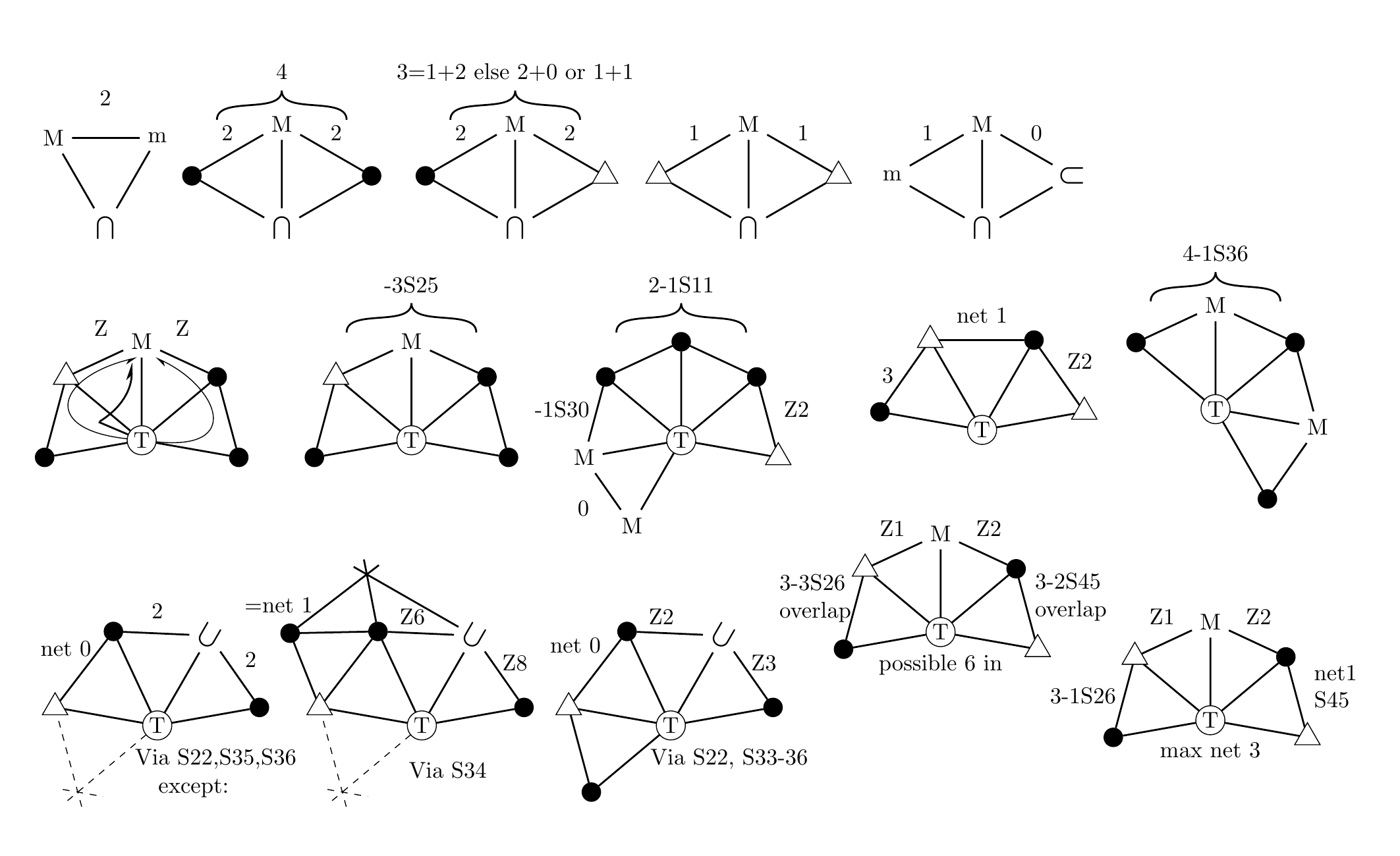}
  \end{center}
  \caption{Limits, negative inflow, and recording of consequential flows
  \label{fig-MmLimits_NegInflow}						}
\end{figure}

\subsection{Limits}
Lemmas ~\ref{lem-cross-mm} and ~\ref{lem-55655}  established maximum limits on the transfers across (minor,minor) rim edges.
A similar limitation is needed for  cross Major-minor edge transfers.
  In the absence of specification of other vertices on the rim of the wheel, the maximum transfer across a (Major,minor) edge is 2 units, across a mMm pair of consecutive edges is 2 units for 6M6 (none of sources 21 to 26 has {6,6} shoulders), and 4 units for 5M5 (Figure ~\ref{fig-MmLimits_NegInflow}, top row). The limit of  3  units for 5M6 is an improvement from that used for the 10-wheel (Figure ~\ref{fig-10Targets}), but was delayed to this point because it is a  consequence of Section ~\ref{ssec-Rule28}. None of sources 21 to 28 has a Major shoulder, so the maximum flow across the mM edge of mMM is 1 unit.

\subsection{Outflow as Negative Inflow} 
A central 7- or 8-valent vertex may have outflows to adjacent Major vertices on its rim.
Since each source structure contains a zero flow target assuring one-way flow, the edge or edges of that zero inflow could be labelled with a Z for Zero inflow. A better use of that space is to indicate the outflow to the Major rim vertex from the central Major vertex as 
a negative inflow (Figure ~\ref{fig-MmLimits_NegInflow}). 
The outflow from Rule 25, 7[56M55] is 2+1 units from 7 to M and there is no counterflow from M to 7. Combining the 3 units of outflow with the zero units of inflow, the maximum inflow across the 5M6 pair of edges i.e. between the rim vertex M and the central 7, across all four shoulder edges is --3 units and labelled --3S25.

\subsection{Recording of Consequential Flows}
\label{ssec-Consequential Flows}
When a central vertex of a wheel has inflow due to vertices beyond the rim, other rules may become active and their outflow may be included to give a lower value as a maximum net inflow.

For 7[M$_1$M$_2$5556], the maximum inflow across MM is zero. The 56 edge is labelled Z2 indicating it is the crossing edge of Zero Flow Target 7402.  Source rule 30 generates a flow of 1 unit to M$_2$ along the 5M$_2$ shoulder edge so that maximum inflow is labelled -1S30. A 5cap on 555 activates exactly one of rules 11 a, b, or c giving a net inflow of at most 1 unit. This maximum net inflow across the 555 pair of edges is labelled as 2-1S11. 
The result of an alternate possibility would be recorded after a comma.  In this case, there is only one alternate, no 5cap, producing zero inflow. The maximum net inflow is the maximum of the possibilities. The zero option from no flow, in this case across 555, can be omitted when the result of any other option is $\geq$0.

T=7[5656] has a maximum inflow of zero units across the last 56 edge and the only positive inflow across the middle 65 edge is two units from a lone 5cap activating the outflow from Source Rule 14.
This is recorded as Z2 on the last 56 edge and either 2-1S14 or net 1 for the middle 65 edge.
n1 or n0 is also used when space does not permit the full word "net", signalling to the reader that an inflow and its consequential outflow rule are being applied.

The maximum inflow to T across a 5M5 pair of edges is 4 units but for  T=7[5M5M5]  the net inflow can be indicated as 4-1S37 because the maximum 4 unit transfer across the 5M5 pair to a 7-valent target is only from M=7[55T55], precisely the condition of Rule 37 and its one unit of outflow from T to the other Major vertex. More will be said in Corollary ~\ref{cor-7[5M5M5]}. 

The maximum flow across a 66 edge is 3 units, but Rules 62a and b give a net inflow limit of  at most 2 units across the 66 edge of 7[66M55], and would be indicated on that edge as 3-1S62 or net 2.

Rules 22, 35 and 36, corresponding to 3, 1 and 1 units flowing across the 65 edge of 7[65M5] produce a net inflow of 0 across this edge. 
If the inflow across the 65 edge is in fact zero, then there is no impediment to a possible inflow of 2+2 units across the 5M5 edges as shown.
Otherwise,  the lone 5cap and rule 34 give net inflow of exactly 1 unit while targets 7406 and 7408 limit the flow across 5M5. This limits can be restricted to 7[*65M5] due to an improved result:
 T=7[565M5] has
zero flow across the 5M5 edges by ZFTs 7402 and 7503 and a net inflow of zero across the adjacent 65 edge 
because rule 33 applies.

Because the units of flow are ORable, care must be taken to avoid double counting an outflow, in particular when other sources are combined with Sources 22 or 26 to the same target.
For example, 7[56M56] could have 3 units inflow across each 56 edge for which S26 describes 3 units of outflow to M and S45 two units to M, suggesting 6 units in, 5 units out, a wrong result. By describing the contribution of S26 as 3-1S26,   this indicates that only the minor hub trajectory of S26 is being counted as output to M, the other two units from the 55 cap are counting as net inflow to its 7-target. This avoids an overlap with outflow trajectories from 7=T to M of Source 45, due to a 5cap on the other 56 edge. That net possible inflow is one of 3-2 for a 55 cap, 2-2 for a lone 5cap, or at most 1-0 for any other flow across that 56 edge. The worst case is therefore listed as net 1 by S45. Similarly, 3-1S22 would indicate that of the 3 units from the 55cap, only the minor hub trajectory is being output, the other 2 units contribute to the possible net inflow to the target/source 7-valent vertex. Although the 2 units flowing across the (m,M) edge is not exclusively from Source 22 or 26, the minor hub unit can only be from those sources, so they must be part of the structure whenever 2 units flows across that (m,M) edge.



Limited flow 7-targets are given in Appendix ~\ref{app-ZF7}. 
The 4 digit label is 7knn where k is the number of 
first neighbours of the target, k$\epsilon$\{4,5,6\}, and nn is a two-digit sequence number, with possible gaps and that is their order of presentation.

Articulated 7-targets are numbered 7A1 to 7A9. Note that 7A5, 7A7, and the mm=66 versions of 7A9 are 1-flow targets.

The crossing edge of a target structure could be unspecified as (L,R) but all other vertices must be fully specified to produce a configuration to be tested for reducibility. If a vertex of a target diagram is 
specified as minor, that diagram represents a set of two fully specified (apart from L and R) GG targets, one with the minor vertex 5-valent and a second with it 6-valent. Not only does a target diagram represent a set of reducible configurations, it could also represent a set of fully specified targets and also the union of the sets of their reducible configurations. Another notation to describe a set of targets is introduced now.

Target 7441 of Appendix ~\ref{app-ZF7} looks fully specified and contains 7[56T5], the common structure of sources 24 to 26 (and also 42 to 45). The little '3' inside the 7-valent vertex means that this target diagram represents not T=7[607L=5R], but rather the three fully specified targets where the 7[56T5] is expanded into the three sources 24 to 26, producing the fully specified targets shown as 7441a, 7441b, and 7441c. Each fully specified target is a zero flow target so there is zero flow across the crossing edge of 7441whenever a source 24 to 26 and its 2 units of flow across the (7,6) edge is present on the 7[56T5]. 
Both ways can be used to represent sets. To include source 21 in this arrangement, the common structure is 7[5mT5] and a '4' would be coded inside the circle representing the 7-valent vertex. In this way four fully specified targets 7549 and 7550 a, b, and c could be combined into one target diagram representing all the configurations formable from those fully specified targets when combined with each of  the possible sources that are compatible and consistent with the limited flow value. When this limited flow value is zero, one target diagram may represent over a hunderd reducible configurations. This coding is used frequently in Appendix ~\ref{app-ZF8}: Limited Flows for an 8-valent Target.

As mentioned earlier, the set of reducible configurations used here is not minimal. Rules and targets were chosen to minimize the number of configurations that were not D- or E-reducible. 
For example,  the common structure target drawn as 7441 is a Zero Flow Target, all the configurations are reducible, but  5 of these configurations need the reducer finding and verification step. Using the three targets 7441a, b, and c instead takes  extra  compute time to resolve three times as many configurations, but these are all D- or E-reducible.  Using 7441a, b, and c
avoids needing to find reducers for those five configurations. 
Similarly, the common structure target of 7550 is a ZFT, but requires finding a reducer for another 5 resulting configurations. Using the full sources 24 to 26, all the resulting configurations are D- or E-reducible.
This tactic is not needed for target 7545 because only T=7[5706606] needs a reducer, and for this case, the set of configurations represented by target 7408, all D- or E-reducible, confirms that should a full Source 24 to 26 be used, the resulting configuration is, in fact, D- or E-reducible.

The targets 7421 to 7550 support limits on the flow crossing an edge in the following more general target sets labelled L1 to L6.

\begin{figure}[ht] 
  \begin{center}
    \includegraphics[
    width=6in
    ]{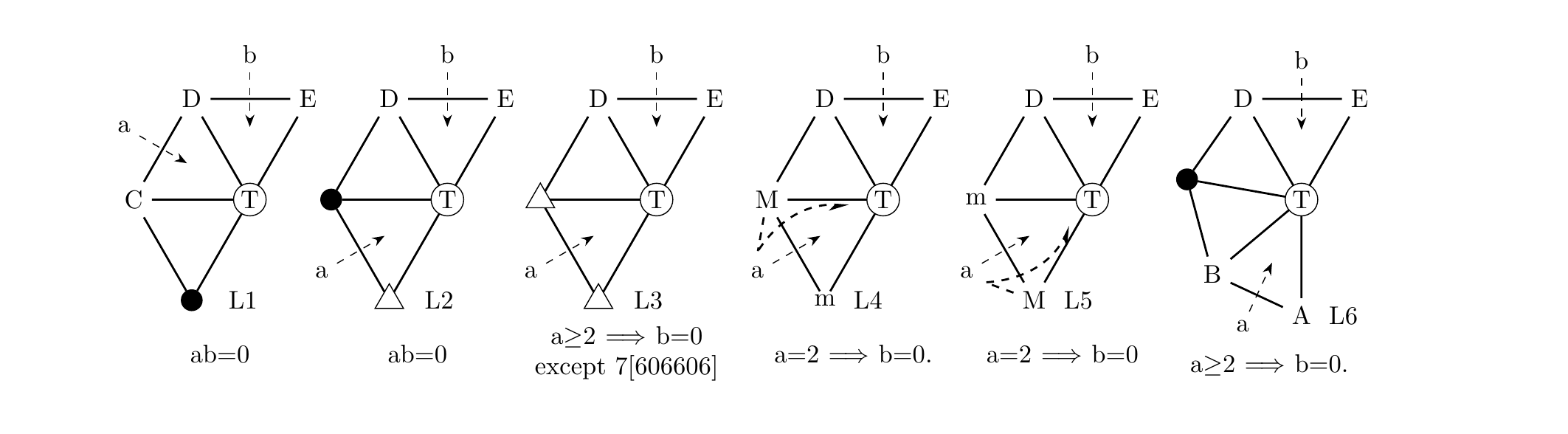}
  \end{center}
  \caption{Combinations for Zero Flow to a 7-valent Target}
  \label{fig-SixV7Lemma}
\end{figure}

\begin{lemma} Refer to Figure ~\ref{fig-SixV7Lemma}. 
It is assumed and required that the flows a and b are non-negative inflows to the target 7-valent vertex.
\label{lem-sixV7Lemmas}
\begin{enumerate}
\item For T=7[5CDE], if there is a positive transfer across (C,D) to T=7, then zero units cross (D,E) to T=7.

\item For T=7[65DE], if there is a positive transfer across 65 to T=7, then zero units cross (D,E) to T=7.

\item For T=7[66DE], if there is a 2- or 3-unit transfer across the 66 edge, then zero units cross (D,E), except a=b=2 units for T=7[6(*)0(*)66(*)0(*)6].

\item For T=7[mMDE], if there is a 2-unit transfer across the (m,M) edge to T=7, 
then zero units cross (D,E) to T=7.

\item For T=7[MmDE], if there is  a 2-unit transfer across the (M,m) edge to T=7, then zero units cross (D,E) to T=7.

\item For T=7[AB5DE], if there is a 2- or 3-unit transfer across the (A,B) edge to T=7, 
then zero units cross (D,E) to T=7.
\end{enumerate}
\end{lemma}

\begin{proof} 
\begin{enumerate}
\item L1, 7[5CDE]: If C is minor then b=0 by 7402. Otherwise, C is Major and a$>$0 requires D to be minor.
If D is 5-valent, then possible sources with \{5,5\} shoulders and a 7 target are Rules 21 (two orientations), 33, 34, 37 and 81 to 85. 7-Targets 7421 to 7429 assert zero flow across DE.

\noindent Otherwise, D is 6-valent and a$>$0 requires C=7, a source with \{6,5\} shoulders, allowing a 7-valent target, and flow across the (7,6) edge. Rules  24 to 26 and 42 to 45 all have 7[56T5] as a substructure for which 7425 asserts zero flow across DE. 
Rules 61 and 62 also support a transfer across a (C,D)=(7,6) edge to a 7-valent target for which 7426 asserts zero flow across DE.

\item L2, 7[65DE]:  7408  and 7431 to 7433 cover a$>$0 transfers from  rules 3, 4, 5 and 11a, 11b,  14 and 15. 

\item L3, 7[66DE]: ZF7Ts 7408, 7434 to 7436, and 7503 cover  2-unit transfers from rules  3, 4,  6 and 7 with either 6-valent vertex as the hub for the inflow to T=7. Only 7[606606]-12 is SDI. Since the 55 caps on LR are also part of the 7408 target and these are not an exception, they are reducible. Including symmetry, the exception is 7[6(*)0(*)66(*)0(*)6]-12.

\item L4, 7[mMDE], a=2: $\implies$ minor hub transfer $\implies$ 
sources 21 to 28 but not 22. 21 and 28 are covered by target 7426 and 7441 covers sources 24 to 26.

\item L5, 7[MmDE]:  a=2 $\implies$ sources 21, 24-26 and 28, and are covered by targets 7545 and 7446.
\item L6, 7[AB5DE]: B=5 is covered by 7402. B=7, a=2$\implies$ sources 21 and 24-26 which are  covered by targets 7549 and 7550.
BA=65 is covered by 7503. BA=67, a=2$\implies$ sources 24 to 26 producing reducible 7[56075]-11. Lastly we have BA=66 for which a$\geq$2 is possible with a 5cap covered by 7509, and the
6-hub neighbour sequences, (hub=6)[T=76655] and (hub=6)[T=766655], containing a reducible configuration   6[555] or a reducible configuration from zero flow targets 7408, 7548 or 7604.
\end{enumerate}
\end{proof}

\begin{corollary}Let 7[BCDExxx] be a 7-wheel with consecutive neighbour vertices B, C, D and E, unspecified vertices x, and an inflow of a units across BC and b units across DE. Then:
\label{cor-7[BCDE]} 
\begin{enumerate}
\item a+b $\le$4 and  a+b=4 only for 7[606606]-12 SDI.
\item a+b=3 only for either a=3, b=0 or a=0, b=3.
\item Otherwise a + b $\leq$ 2.
\end{enumerate}
\end{corollary}
\begin{proof}
 2 or more units across BC is an instance of ZFT 7402 or cases L2 through L5  all with zero flow across DE, except for 7[606606]-12. By symmetry, a similar conclusion applies to 2 or more units across DE. Otherwise, the maximum flows across BC and DE are 1 unit each.
\end{proof}

\begin{corollary}Let 7[AB5DExx] be a 7-wheel with 
a possible non-negative inflow of a units across AB and b units across DE. Then:
\label{cor-7[AB5DE]} 
\begin{enumerate}
\item a+b$\leq$3 and a+b=3 only for either a=3, b=0 or a=0, b=3.
\item Otherwise a + b $\le$ 2.
\end{enumerate}
\end{corollary}
\begin{proof}
2 or more units across AB is an instance of case L6 giving zero flow across DE. By symmetry, a similar conclusion applies to 2 or more units across DE. Otherwise, the maximum flows across AB and DE are 1 unit each.
\end{proof}

\begin{corollary}For 7[5M5M5], the net total inflow across the four 5M edges to 7 is at most 3 units.
\label{cor-7[5M5M5]}
\end{corollary}
\begin{proof} 
7[5M5M5] is shown in Figure
 ~\ref{fig-MmLimits_NegInflow}. By L1, 7[5CDE], at most one of the middle two rim edges has positive flow crossing it to the 7-valent vertex. By Corollary ~\ref{cor-7[AB5DE]}, the total across the two end edges, which are not (m,m), is at most 2 units. To exceed 3 units total, the middle transfers must be 2+0 and now the transfer across at least one of the two end edges must be 0 units, 
making the four inflows 2+2+0+0 by Corollary ~\ref{cor-7[BCDE]}.  Now this 2+2 must be from the source structure for Rule 37 and its outflow of 1 unit gives a net inflow of at most 3 units across all four edges.
\end{proof}

In Appendix ~\ref{app-7Nhbds}, this result is indicated by labelling the possible inflows across these edges as 3a, 3b, 3c and 3d.

\begin{figure}[ht] 
  \begin{center}
    \includegraphics[
    width=3.4in
    ]{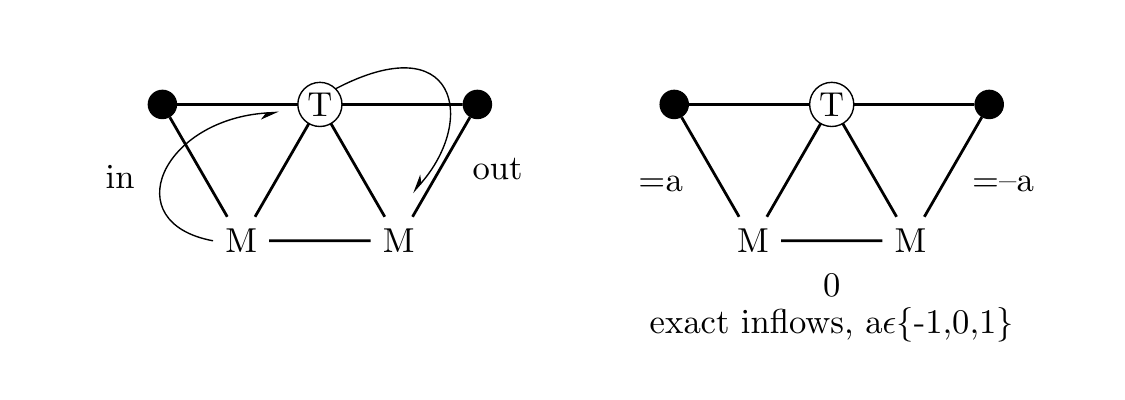}
  \end{center}
  \caption{Flows across 5M and M5 of 7[5MM5]}
  \label{fig-7 5MM5}
\end{figure}

\begin{lemma}The net inflow across the two (5,M) edges of T=7[5MM5] is zero.
\label{lem-7[5MM5]}
\end{lemma}
\begin{proof} Refer to Figure ~\ref{fig-7 5MM5}. All cross (Major, 5) edge flows with \{5, Major\} shoulders are only one unit, from Rules 30, 31, 32, 54 to 58, 66 and 67, but 55 and 56 do not allow a 7-valent target. If the inflow is from Rule 30, the result is Rule 54 and an outflow of 1 unit. If the inflow is from Rule 31, the result is Rule 55 and its outflow of 1 unit. If the inflow is from Rule 32, the result is Rule 56 and its outflow of 1 unit. Inflow from Rule 66 yields Rule 57 and inflow from Rule 67 yields Rule 58. An inflow to 7[5MM5] from any of Rules 54, 57, or 58 yields the reducible configuration 7[557075]-12. 
\end{proof}
In Appendix B, this result is recorded by labelling the 5M edges with exact flow values of a and -a where it is understood that a $\in$ \{--1,0,1\}.

\subsection {Enumeration and evaluation of neighbourhoods of 7-valent wheels }

Grouped by the size of the separating ring, GG first neighbourhoods of a 7-wheel and some other configurations are categorized in Table ~\ref{tabl-7[5mm...5]}.

\begin{table}[ht]
\caption{ GG first neighbourhoods of a 7-wheel  
and some other configurations.
}\label{tabl-7[5mm...5]}
\begin{tabular}{r | l | l }
\label{tabl-7-nhbdsRed-SDI}
Ring Size		&Reducible		&Symmetrically D-Irreducible (SDI)\\
\toprule
8: &7[5555xxx].					&\\
9:&7[5565xxx], 7[555x505x].		&\\
10:&7[56565xx], 7[55665xx], 7[555x506x], 		&7[5665xxx].\\
&7[565x505x].		&\\
11:&7[56665xx], 7[565665x], 7[565x506x],		&7[55706xxx].\\
&  7[555x606x].		&\\
12: &7[566665x], 7[5566666], 7[557075]. 	&7[565x606x], 7[56706xxx], 7[57075xxx]. \\
13:&7[5666666].				&\\
14:&7[6666666].				&\\
\hline
\hline
\end{tabular}
\end{table}


Similar to the 10-valent vetex, first neighbourhoods are enumerated in terms of the number of 5-valent and at least 6-valent vertices, which may be required to be Major by a reducible configuration or subcased as degree 6 and then as Major.

The exhaustive enumeration of symmetric neighbourhoods exploits symmetry to enumerate a pair of unspecified vertices of degree at least 6. In particular, a symmetric structure is extended by 66, MM, and lastly 6M, recursing whenever the neighbourhood has maintained the symmetry and adding M6 if not. For example, neighbourhoods 7[M5M ****] are examined in the sequence:

7[6M5M6 **]: 7[6M5M6 66], 7[6M5M6 MM], 7[6M5M6 6M], 

7[MM5MM **]: 7[MM5MM 66], 7[MM5MM MM], 7[MM5MM 6M], 

7[6M5MM **]: 7[6M5MM 66], 7[6M5MM MM], 7[6M5MM 6M], 7[6M5MM M6].



Given a central Major vertex and its  full first neighbourhood described in terms of 5-valent, 6-valent, *, and Major vertices, there will be the known exact or maximums on the Mayer 2,3,4 inflow transfers from the 5-valent vertices. If the central Major vertex is 7-valent,   there  may also be known outflow from discharging rules consistent with this first neighbourhood: Rules 21, 24, 25, 28,  30, 42, 44, 51, 61, 73 and 74.
The initial charge together with these guaranteed transfers are combined to give a preliminary charge or precharge. Additional charge transfers may be specified, for example as the result of a division into subcases, and these are also incorporated into this precharge. 

The remaining charge transfers, unspecified inflows across the rim edges, will be dependent on the vertices in the second or further neighbourhoods. Limits on these possible inflows from Zero Flow Targets and the Lemmata above contribute upper bounds on these possible inflows or combinations of possible inflows. If the total of these upper bounds is not enough to overcharge the Major vertex, that neighbourhood is disposed of. As in the case of an 11 or 10-valent vertex, if only one unit of overcharging is possible, then all maximums must be active to achieve overcharging. If  two units of overcharging is possible including 2 units from a 5cap on a 55 rim edge or on a 555 pair of rim edges, that cap would be forced because otherwise the inflow across those edges is zero units. In these cases the forced additional structure may reduce a previous upper bound enough to eliminate the possibility of overcharging, or develop a reducible configuration. 

The flow across a 656 pair of rim edges is 6 units for 60506, 3 units for the 56 cap, 2 units for a 5M cap, and otherwise zero units. If the maximum possible overcharge, including a possible 6 units from a 656 neighbour group, is at most 3 units,  then the 60506 structure is required for  overcharging.


Another frequently used subcasing is the flow across the 56 edge of  M[*56].  This will be 3 units for the 55 cap, 2 units for a lone 5 cap, or 1 unit from sources 11b or the 56 cap, and that is all. After the subcase of a 5cap with its 2 or possibly 3 units of inflow, the cap will be at least 6-valent  and the maximum flow across 56 will be at most 1 unit. 

The possible sources for a positive flow across a (Major, minor) edge are sometimes limited and these  can be enumerated.
For example, only one source structure delivers 2 units across the (5,M) edge of  T=7[5L=5R=M6xxx], namely Rule 28: 7[555T6*]. Once neighbourhoods with this source are disposed of, the maximum is only 1 unit. If that still allows overcharging,
that one unit must come from Rule 67: 7[56050T*]. Once that subcase is handled satisfactorily, the inflow across the (5,M) edge is zero for any neighbourhood that in some other way overcharges the central 7-valent vertex. 

A worst case would be the discovery of a neighbourhood of a vertex that is not reducible but which overcharges the central vertex, in which case the proof as is has failed. 
An example of this was 7[5060(5)5M55M], not GG ($\implies$ not reducible), so allowable in a m5CPG. The precharge using Mayer and rule 21 is --10+12-3-3=--4, and 5 units cross the 565 edges giving a possible final charge of +1. New rules were tried to send a fourth unit to one of the Major vertices before Rule 22 was found to be (in my mind), most appropriate for this case, while also being functional in other cases. (Other candidates were 7[565T55] and 7[6050T55].)
Another is the neighbourhood of Rule 53 for which the final charge was found to be +1. Rule 53 with its outflow of 1 unit was introduced and checked against all limited flow targets to confirm that their claims were sustained. 

Similar worst cases happened several times over the past four decades before the current set of rules and their intricacies were developed.

As an example of this enumeration task and the notations used to record it for verification, consider the configurations T=7[606606xxx], those with the exception to 7408 of Appendix A and hence possible to appear in a m5CPG.
By avoiding 7[566606]-12 of 7503, allowable neighbourhoods are T=7[606606 $*$x$*$] and  we subcase the $*$x$*$. Since the structure is symmetric, the subcase sequence is $**$=66, then MM, and finally 6M. 

\begin{figure}[ht] 
  \centering
    \includegraphics[
    width=4.6in
    ]{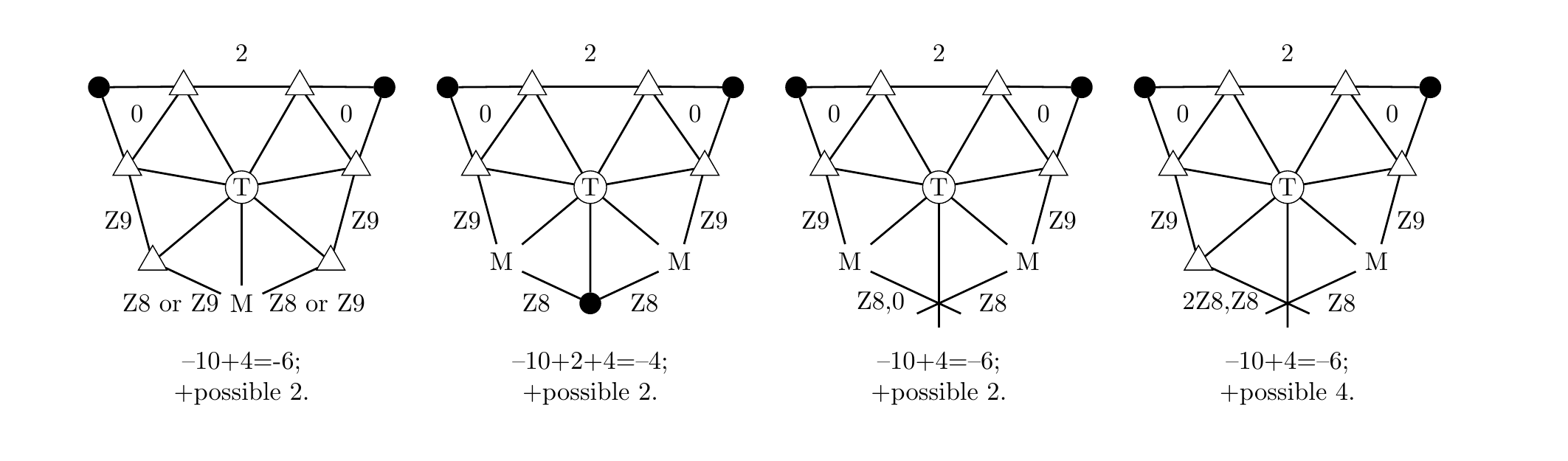}
  \caption{Possible neighbourhoods of 7[606606xxx] and maximum inflows.
  \label{fig-7 606606}		}
\end{figure}

7[606606 656] contains reducible 7[56606]-11 of 7402, 7[606606 666]
contains reducible 7[6666 666]-14,
leaving only 7[606606 6M6], the leftmost neighbourhood of Figure ~\ref{fig-7 606606}.  *x* as MxM is  subcased as the next two neighbourhoods. Lastly, *x* as 6xM drops x=5 by  reducible 7[56606]-11 leaving the rightmost neighbourhood.

The first line below a cluster evaluates the preliminary charge (precharge): the initial charge plus the transfers from first neighbourhood 5-valent vertices (there is only one here), and the known across rim edge transfers. 

7408 has only the one exception to zero flow. This means a 55cap produces a reducible configuration so the additional inflow across the capped 66 edges 
is zero as indicated inside the caps. The (6,6) edge between the caps now delivers at most 2 units by a lone 5cap. 

Most of the remaining rim edges are crossing edges of Zero Flow to 7 Targets 7[6066LR] and 7[606LR], labelled 7509 and 7408 in Appendix ~\ref{app-ZF7}. Those zero limits are indicated in the final charge analysis diagrams by Z9 and Z8 respectively. In the first diagram, the (6,M) edges match (L,R) for both Z9 and Z8-not-the-exception, so either indicator could be used.
The commas in the last two pictures separate the subcases of the * at the bottom being 6 or M. *=6 produces zero flow across the (M,6) edges because that is not the exception to 7408, and   *=M  gives 0 flow since there is no charge transfer across a (M,M) edge. This latter case is vacuously part of the definition of a Zero Flow Target, so Z8 suffices as in the right sides of the last two diagrams.  

For the final diagram, "2Z8, Z8" is also the result of subcasing * into 6 or M. This reports that the maximum possible inflow across the (6,*) edge is 2 units if *=6, by the exception to 7408, and zero units by the same 7408  if *=M since this would not be the exception. The largest possible inflow across this edge  is the maximum of these comma-separated coded values.

In every case the maximum inflow is not enough to overcharge the central 7-valent vertex.

An exhaustive enumeration of neighbourhoods 7[606606xxx], where each vertex x could be 5, 6, or Major would have around $3^3/2$ cases with possible subcases. Using * for degree $>$5,  
generalized Zero Flow Targets, and comma subcasing reduces this to four diagrams.

\begin{figure}[ht] 
  \centering
    \includegraphics[
    width=5in
    ]{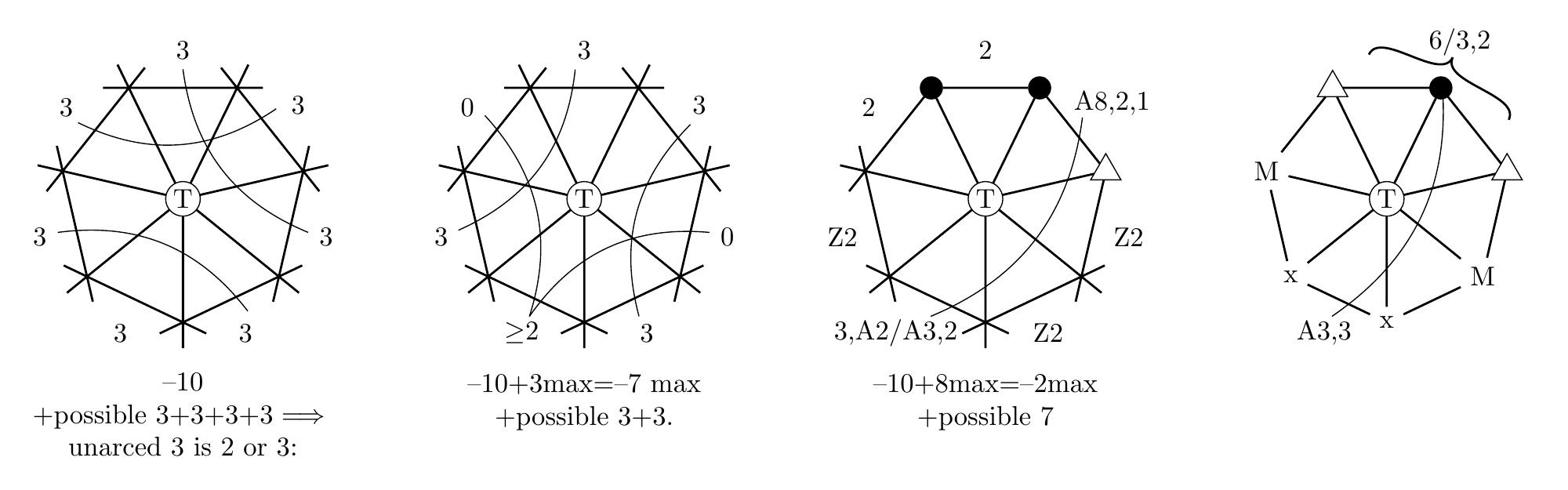}
  \caption{Maximum inflow to 7[**** ***] and arced pairs of transfers.
  \label{fig-7[**** ***]}		}
\end{figure}

Having shown that all the occurrences of the exception to 7408 are either reducible or not overcharged, case 1 of Corollary ~\ref{cor-7[BCDE]} is no longer a consideration for the remaining neighbourhoods of a 7-valent vertex in a m5CPG.  In this light, consider 7[**** *** ]: zero 5-valent neighbours, in Figure ~\ref{fig-7[**** ***]} where each rim edge may deliver a maximum of 3 units. 
Using Corollary ~\ref{cor-7[BCDE]}  the three arcs indicate the now maximum inflow of 3 units across the PAIR of rim edges connected by each arc. The unconnected edge contributes a maximum of 3 units for a possible inflow of 12 units. Overflow requires that unconnected edge contributes 2 or 3 units and a reevaluation of the flows shows that in this case, the total inflow is at most 9 units.   A 7-valent vertex whose immediate neighbourhood has no 5-valent vertices will not be overcharged.

In the log of enumeration of neighbourhoods, a limit on the maximum flow across the rim edges is reported.  
When the sum of flows across a pair of rim edges is constrained,
  this is indicated by joining them with an arc.  
When constrained by Corollary ~\ref{cor-7[BCDE]} or ~\ref{cor-7[AB5DE]},
    the maximum of the sum of the flows is the same as the maximum of the individual maximums and this is easily 
remembered from Figures ~\ref{fig-Cross mm edge summary} and ~\ref{fig-MmLimits_NegInflow},
but that result depends on two conditions. 
First, the conclusion values in  Lemma ~\ref{lem-sixV7Lemmas} are b=0;
if they weren't then the sum could exceed this maximum of maximums.
Secondly, the maximum flow across a single edge is 3 units, and the lemmas are still active if the flow across an edge is only 2 units.  Were this not the case, the sum might be 2+2 units. 

The sum of pairs of flows across rim edges can be constrained in other ways, and is shown by comma subcasing. Consider 7[556****] in Figure ~\ref{fig-7[**** ***]}. If the flow across the middle ** edge is 3 units, i.e. a 55cap on **=66, then the flow across the 56 edge is zero by Articulating  target 7A8. In this case, there is a reciprocal target limiting the flow on the same two edges.  If the flow across that 56 edge is its maximum of 2 units, and that happens in only two ways, then the flow across the middle ** edge is zero by  7A2 or 7A3. A maximum on their sum is max\{ 3+0, 0+2,  2+1\} =3, where the 2+1 is the case where neither edge has its flow equal to its "trigger" value, in these cases, their maximums.
In the discharging log,  the options are separated by a comma or slash.
Since the maximum total flow across both edges is equal to the maximum of the maximums of the individual flows, these edges are joined by a solid arc.
A dashed arc indicates that care should be taken as the total flow across the two  or more rim edges may exceed the maximum of maximums. 

Target 7A3 may also limit the sum of the flows on 3 edges as in 7[xxM656M] where the xx edge is joined by an arc to the middle 5-valent vertex of the 656 pair of edges. The 656 pair  delivering 6 or 3 units, limits the flow across the xx edge to 0 units by 7A3. Otherwise, the flow across the three edges is at most 2+3 units. 
The xx edge 
 can be 
 labelled with A3,3 and the 656 group tagged 6/3,2. There is no reciprocal target here. 
Only one of the arced edges has a lettered limit code.
The ",3" is the maximum flow across the ** edge and the ",2" is the maximum allowable value that does not trigger the A3 target.
The 656 group could also  be labelled merely with 6, its maximum, and the xx edge by A3, with a  solid arc indicating that any possible flow across the xx edge is subsumed by at least a corresponding reduction in the possible 6 units across the 656 pair and the maximum of maximums property holds.
This limited sum of flows across 3 edges is more common with neighbourhoods of 8- and 9-valent vertices.

A reference to a target or a structure of Lemma  ~\ref{lem-sixV7Lemmas} (later ~\ref{lem-ZF8Ts1} or  ~\ref{lem-ZF8Ts2}) is given by its sequence number preceded by a letter. Prefixes "Z", "A", and "L" indicate zero flow, while those starting "1Z", "1A"  or "1L"  indicate an inflow of 1 unit is possible. This limit is typically across the LR edge of the target or L structure,  but if the flow across the LR edge is already known and exceeds the limit specified in
 an L structure, then the limit will be on the edge of the side source, i.e. the contrapositive is being used. As in Figure ~\ref{fig-MmLimits_NegInflow}, -1Sk refers to Source rule k and indicates an outflow of 1 unit (maximum inflow of -1 unit).

\begin{theorem}Vertices of degree 7 are not overcharged by this discharging. 
\end{theorem}

\begin{proof} 7[ $\text5^\text{7}$] and 7[ $\text5^\text{6}\ast$] contain reducible 7[5555].  
The remaining cases with five down to one 5-valent first neighbour are enumerated in Appendix ~\ref{app-7Nhbds}. 
Extra reducible configurations are given to dispose of a neighbourhood or to support a claim in the diagram.

A new first neighbourhood, to be followed by its collection of neighbourhoods produced by subcasing, is indicated  by the initial evaluation of the precharge as "--10 +...=...". All the remaining possible inflows are summed and listed  as "+possible...". If the resulting neighbourhood cannot be overcharged, we move on to the next subcase or next first neighbourhood. If it can be overcharged, then this is the assumed condition and its consequences are followed, indicated by "$\implies$" or sometimes ":". The record for the next neighbourhood will 
 continue with the previous precharge plus any consequential flows due to the first of a set of possible subcases. The next subcase of that set can usually be recognized because it will  start with the same prechage. For example, see 7[5M56**6] where the sequence of starts of the precharge calculations is
 --10, --4, --4, --4, --4, --4 for the log of an initial neighbourhood and 4 subcases of which the fourth has a forced subcase. 

A subcasing will often be described as "=3" meaning a 55 cap is on a (m,m) rim edge, or "=3/2" meaning a 5cap is on the 56 edge of *56, since this is the only cap that guarantees at least 2 units. This is  easier to express compared to drawing the caps and shows the inflow amount explicitly.
After "=3", further subcasing may be indicated by subsequent 
 "=2" meaning any possible 2 unit source, but even without knowing which one, could  still have consequences by Lemma~\ref{lem-sixV7Lemmas}, and then "$\leq$1". The "$\leq$" is a reminder that larger flows were possible but those cases have been examined and eliminated, whereas a "1" indicates the maximum inflow is 1 unit from the neighbourhood, for example T[5MM] has a flow of at most 1 unit across the 5M edge to T. 
\end{proof}

\begin{remark}
As a example of the trade-off between adding reducible configurations and simplifying the discharging, target cluster 7A5\_1 can be established as a Zero flow target by adding an extra 24 D-reducible configurations to the set S, of which 6 are bounded by a 16-ring. As  a consequence, Source 66 is no longer needed to show that the 7[M1((6))50526***] subcase of 7[M556***] is not overcharged, and neither is Source 57 needed, shrinking by 2, or 4 including their vertical mirror images, the number of possible sources that are used to establish the zero flow targets.
\end{remark}

\section{Analysis Support for, and Enumeration of, Neighbourhoods of an 8-Wheel}
8[55555xxx]-9 is reducible so an 8-wheel with eight or seven 5-valent neighbours will not occur in an m5CPG. Since the maximum transfer across a *** pair of adjacent edges is 5 units, the final charge on an 8-valent vertex with no 5-valent neighbours is at most -20+4*5=0.  
It remains to examine first neighbourhoods with six down to one 5-valent neighbour.  



Grouped by the size of the separating ring, the GG first neighbourhoods of an 8-wheel and some other configurations are categorized 
in Table ~\ref{tabl-8[5mm...5]}. 

\begin{table}[ht]
\caption{ GG first neighbourhoods of a 8-valent vertex.}
\label{tabl-8[5mm...5]}
\begin{tabular}{r | l | l  }
Ring Size		&Reducible		&Symmetrically D-Irreducible (SDI)\\
\toprule
9: &8[55555xxx].	&\\
10:&8[55565xxx], 8[55655xxx], 8[555x555x], 8[5555x505x].	& \\
11:&8[56565xxx], 8[565565xx], 8[555665xx],&8[55665xxx], 8[556655xx].\\
&8[555x565x], 8[5655x505x], 8[5555x506x].	&\\
&8[555705xxx], 8[557055xxx],  8[570555xxx]	&\\
12:&8[565665xx], 8[5\{5,5,6,6,6\}5x], 8[56655665]. &8[56665xxx], 8[556665xx].\\
&8[5655x506x], 8[5655x605x], 8[5555x606x].&8[565x565x], 8[5665x505x].\\
13: &8[5566665x], 8[5656665x], 8[5665665x], 8[56656665].&8[566665xx].\\
&8[5655x606].&8[5665x506x].\\
14:&8[5666665x], 8[56665666], 8[56656666], 8[55666666].&8[5665x606x].\\
15:&8[56666666].	&\\
16:&8[66666666].&8[66666666] but A-Reducible.\cite{Franklin}\\
\hline
\hline
\end{tabular}\\
\end{table}


Several limited flow 8-valent target structures are labelled and listed in Appendix ~\ref{app-ZF8}. 
The 4 digit label is 8knn where k is the number of 
first neighbours of the target, k$\epsilon$\{5,6,7\}, and nn is a two-digit sequence number, with possible gaps and that is their order of presentation. 
Articulated targets are labelled 8A1 to 8A9.

Some targets have a suffix or qualifier. The Articulating targets 8A6a and 8A6b are the two source structures with 2 units crossing the 56 edge of 8[xLR5x655]. Together this pair of zero flow targets limit the total across the LR edge plus that across the 56 edge.
There are similar a and b pairings for 8503, 8504, 8A2 and 8A5.
When these pairs are used in Appendix ~\ref{app-ZF8}, they are referenced   as A6, Z3, Z4, A2 and A5. 
Target 8A2c has a 60506 cluster. Target 8A5c similarly has a 60506 cluster, while 8A5b has the 50516 cluster to pair with 8A5a.
Similar to 7A3, the flow across a 656 pair of edges is 6 or at most 3 units and the maximum of 6 units subsumes any flow across the LR edge in targets 8A2 and 8A5

A similar subsuming occurs with 8A4a and b, where 5,4, or 3 units crossing the 565 pair of edges  requires  a 5cap on a 56 edge and the resulting zero flow across the LR edge. Otherwise, the flow across the 565 pair is at most 2 units and the total across all three edges is again at most 5 units. 
In Appendix ~\ref{app-8Nhbds}, the 565 cluster is tagged as having a maximum inflow of  5 units and its central vertex is  joined by a solid arc to the LR edge which is tagged A4.  Any   positive flow across the LR edge is subsumed by the possible 5 units across the 565 pair.

Articulated targets 8A5 and 8A8 restrict the m,m part of these zero flow targets to avoid needing reducible 16-ring configurations. 8A9-1, extending 8A9, avoids requiring reducibility of 16-ring configurations by being 1-flow targets. 
8A9-1 does not have a reciprocal limited flow target when the maximum of 3 units crosses LR. For example, 8[56***65*] with 3 or 2 units across a 56, i.e. a 5cap, limits the flow across the remote ** edge to 1 unit. Otherwise, the flow across the 56 edge is at most 1 unit and the flow across the ** edge is 3 units, a total of 4 units. The pair of edges is joined by a dashed arc, the 56 edge is labelled 3/2,1 and the ** edge is labelled 1A9,3. Only one edge has a letter coded limit.		

8501, and 8502 are the most useful first neighbourhood  targets, because they have xxx, three consecutive vertices of unspecified degrees. Alas, 8501 has 2 exceptions, each with LR=606. 
Similar to the last case in figure ~\ref{fig-7 606606}, 
when target 8501 is used in Appendix ~\ref{app-8Nhbds} for a neighbourhood 8[5m6LRxxx] , it
enumerates the cross LR flows as "\underline{2Z1},Z1", meaning 2 units from the exception, and otherwise 0 units. To alert the reader to this distinguished case, the 2Z1 is underlined, meaning LR=606, the exception listed with the 8501 target, and this generates the consequences on any arc-attached edges. The default case is Z1, where the inflow across the LR edge is zero by target 8501 and the arc-attached edges are not constrained by LR. 

8[565LR] is a ZFT. Together with 8610 through 8616, all the GG first neighbourhood targets with at most two 6's and two consecutive unspecified vertices are ZFTs: 8[5\{5,m,m\}LRxx] and 8[5\{m,m,LR\}5xx]. 

Similar to the 7-valent vertex, we now give limits on combinations of inflows. Apart from the first structure, the side flows to these L-structures  use only the six 2-unit cross (Major, minor)-edge transfers, sources 21 to 28, and the resulting L-structures have a "delta" of +3 when m=5 and +4 when m=6.  

\begin{figure}[ht] 
  \begin{center}
    \includegraphics[
    width=6in,page=1
    ]{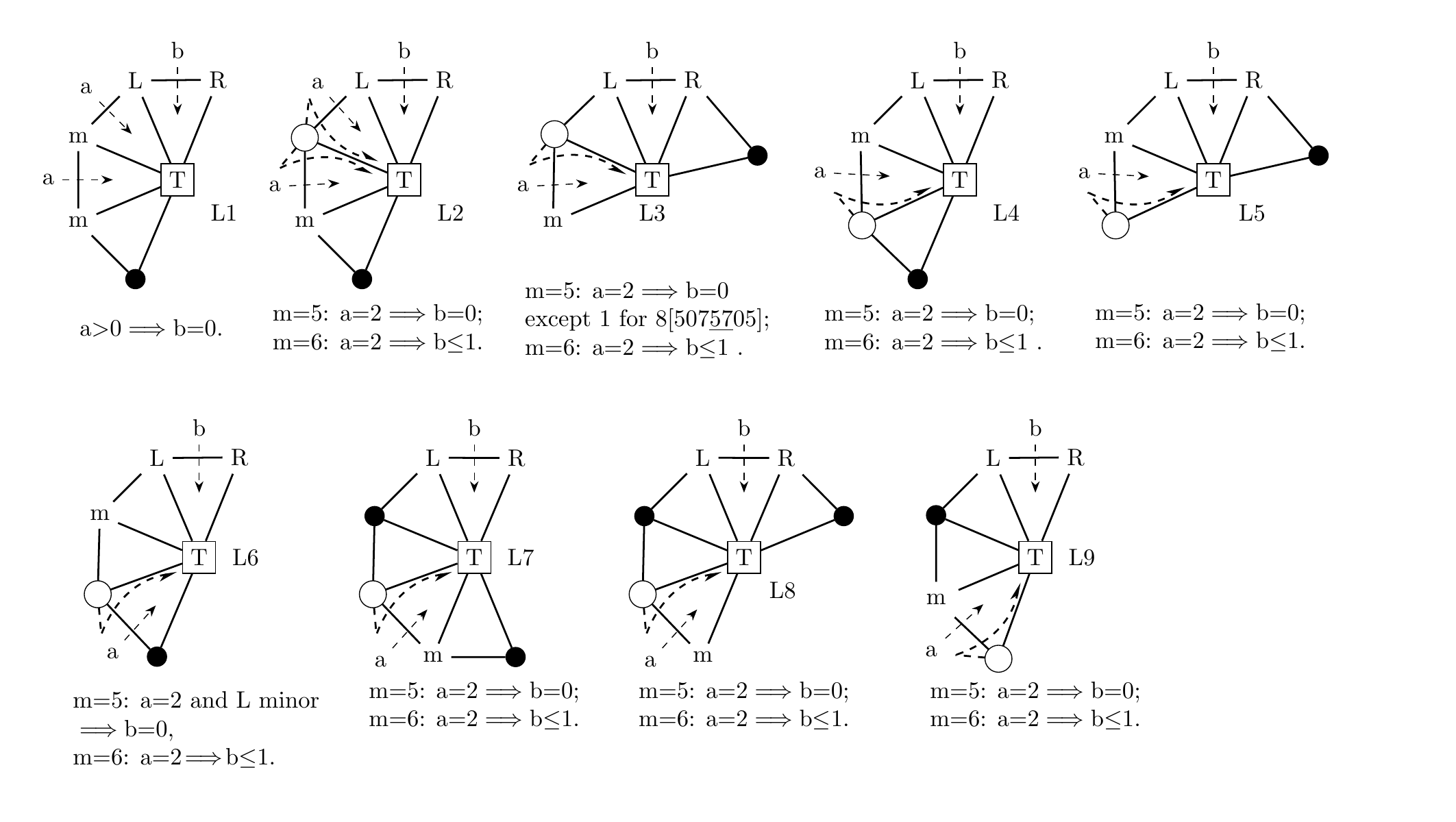}
  \end{center}
  \caption{ Combinations for limited flows to an 8-valent Target}
  \label{fig-V8Lemma1}
\end{figure}

\begin{lemma} 
\label{lem-ZF8Ts1}
Refer to Figure ~\ref{fig-V8Lemma1}. 
It is assumed and required that the flows a and b are non-negative inflows to the target 8-valent vertex. 
For each b, consider either a. 

\begin{enumerate}
\item L1: For T=8[5mmLR], if there is a positive transfer across (m,m) or (m,L) to T=8, then zero units cross (L,R) to T=8.
\item L2: For T=8[557LR], if there is a 2-unit transfer across (5,7) or (7,L) to T=8, then zero units cross (L,R) to T=8.\\
For T=8[567LR], if there is a 2-unit transfer across (6,7) or (7,L) to T=8, then at most 1 unit crosses (L,R) to T=8.
\item L3: For T=8[57LR5], if there is a 2-unit transfer across the (5,7) edge to T=8, then zero units cross LR to T=8 except 1 unit for 8[5075705]-13.\\
For T=8[67LR5], if there is a 2-unit transfer across the (6,7) edge to T=8, then at most 1 unit crosses (L,R) to T=8.
\item L4: For T=8[575LR], if there is a 2-unit transfer across the (7,5) edge to T=8, then zero units cross LR to T=8.\\
For T=8[576LR], if there is a 2-unit transfer across the (7,6) edge to T=8, then at most 1 unit crosses (L,R) to T=8.
\item L5: For T=8[75LR5], if there is  a 2-unit transfer across the (7,5) edge to T=8, then zero units cross LR to T=8.\\
For T=8[76LR5], if there is a 2-unit transfer across the (7,6) edge to T=8, then at most 1 unit crosses (L,R) to T=8.
\end{enumerate}
\setcounter{v8LemmaCounter}{\value{enumi}}
\end{lemma}

\begin{proof}Refer to Appendix ~\ref{app-ZF8} for 8-targets with Zero- and 1-flow limits.
Recall that a target or target diagram represents all the reducible configurations required to meet the limits given in the target or set of targets represented by a target diagram. Here two more styles of target diagrams are given. 8521 represents two zero flow targets, 8[50517L=5R] and  8[505175LR] and this style is also used for target 8523. Secondly,  target diagrams with a minor vertex, thus representing two targets, may have different limits on the cross LR-flow when m=5 and m=6.

\begin{enumerate}
\item  L1, 8[5mmLR]:  The exceptions to Z8501 are 8[556606]-12 and 8[566606]-13, for which ZFT's 8504a and 8505 show that a=0 for each exception and the contrapositive, b$>$0$\implies$a=0, holds.

\item - (\arabic{v8LemmaCounter}):  The  flows a are 2 units produced by sources 21 to 28. The claims are demonstrated for each source:\\
Source 22 as a in L2 upper is 8[5571505]-11, reducible from 8502.  Source 22 in L2 lower and L3  to L5 are covered by ZFTs 8521 and 8523. \\
Source 28 providing 2 units in combination structures L2 upper a, L2 lower a, L3, and L5 have the indicated results by targets 8532-1, 8536, 8539 and 8543 respectively.\\ 
This leaves sources 21 and 24 to 26 having 7[5mT=85] in common, covered by targets 8530, 8534, 8537 and 8540 for L2 upper, L2 lower, L3 and L4 respectively. These sources in L5 give superconfigurations of those in L4.

\end{enumerate}
\end{proof}

The exception for target 8537, 8[5075705] with b=1 also delivers 2 units across the adjacent R5 edge and 1 unit across the 7L edge, a total of 6 units to T. Adding a 5cap on (m,T) gives a reducible configuration from 8647 and adding a 5cap on LR drives a second unit across LR, so that is a reducible configuration from 8537. Employing symmetry, the exception for L3 with m=5 is 1 unit and occurs for only 8[*507(*)5(*)705*x].

In Lemma structures L2 lower to L5, the crossing edges for the a and b flows are one rim edge apart. The next Lemma has a and b two rim edges apart.
\begin{lemma} 
\label{lem-ZF8Ts2}
Refer to Figure ~\ref{fig-V8Lemma1}. \begin{enumerate}
\setcounter{enumi}{\value{v8LemmaCounter}}
\item L6: For T=8[575LR], with L minor, if there is a 2-unit transfer across the (5,7) edge to T=8, then zero units cross (L=m,R) to T=8.     \\ \noindent 
   For T=8[576LR], if there is a 2-unit transfer across the (5,7) edge to T=8, then at most one unit crosses (L,R) to T=8.
\item L7: For T=8[5575LR], if there is a 2-unit transfer across the (5,7) edge to T=8, then zero units cross (L,R) to T=8.\\ \noindent For T=8[5675LR], if there is a transfer of 2 units across (6,7) to T=8, then at most one unit crosses (L,R) to T=8.
\item L8: For T=8[575LR5], if there is a 2-unit transfer across the (5,7) edge to T=8, then zero units cross (L,R) to T=8.\\ \noindent For T=8[675LR5], if there is a transfer of 2 units across (6,7) to T=8, then at most one unit crosses (L,R) to T=8.
\item L9: For T=8[755LR], if there is a 2-unit transfer across the (7,5) edge to T=8, then zero units cross (L,R) to T=8.\\ \noindent For T=8[765LR], if there is a transfer of 2 units across (7,6) to T=8, then at most one unit crosses (L,R) to T=8.
\end{enumerate}
\end{lemma}

\begin{proof} 
Again, the side flows are limited to sources 21 to 28 and the proof starts by considering sources 22 and 28.

For source 22,  L6 to L8,  are covered by ZFT 8621, and 8623 covers L9. Source 28 is not compatible with L7 and L8, and with L6 and L9 it is covered by targets 8546-1 and 8551 respectively.

L6 with m=6 is not compatible with sources 21 or 24 to 26, so target 8545 with m=5 suffices for L6. 

Finally, sources 21 and  24 to 26 are covered by targets 8647, 8649, and 8652 for L7 to L9 respectively.
\end{proof}

\begin{corollary} For 8[5mMLR], the net total inflow across the three edges [mMLR]     to 8 is at most 4 units.\\ For 8[5MmL=minorR], the net total inflow across the three edges 5M, Mm, and LR to 8 is at most 4 units.
\label{cor-8[5mMLR]}
\end{corollary}
\begin{proof} 
Table 1 below shows
the consequences of a maximum flow (indicated by  "=" )  across an edge to the central 8-valent vertex. These consequences are the maximum values of the possible inflow across the other edges. For  M$\geq$8, only the last two rows apply.
For M=7, the maximum possible flow across 575 is 4 units and 3 units across 576. Entries L2, L4, L6 mean the maximum flow is zero by that part of the previous Lemmas while 1L2 etc. mean the maximum flow is 1 unit.  The maximum across LR is 3 units, only from a 55cap on a 6m edge in which case limited flow targets 8501, 8502, 8503a, 8504a and 8505 are used. Z2/4 means the flow is zero by 8502 when R=5 and zero by 8504a when R=6. \\
\begin{table}[ht]
\caption{ Maximum of 4 units across three edges of 8[5mMLR] and 8[5MmL=minorR].}
\label{tabl-8[5m7LR]}

\begin{tabular}{c || c | c | c | c || c | c | c | c || c | c | c | c || c | c | c | c ||}
& 
\multicolumn{4}{c||}{8[55MLR]} &
\multicolumn{4}{c||}{8[56MLR]} &
\multicolumn{4}{c||}{8[5M5L=mR]} & 
\multicolumn{4}{c||}{8[5M6LR]}\\
&5M&ML&LR&\!Sum\!&6M&ML&LR&\!Sum\!&5M&M5&mR&\!Sum\!&5M&M6&LR&\!Sum\!\\
\toprule
&=2&2&L2&4&=2&1&1L2&4&=2&2&L6&4&=2&1&1L6&4\\
&2&=2&L2&4&1&=2&1L2&4&2&=2&L4&4&1&=2&1L4&4\\
&Z2/4&Z2/5&=3&3&Z2/4&Z2/5&=3&3&Z1/3&Z2/4&=3&3&1L6&Z2/4&=3&4\\
\!Else\!&1&1&2&4&1&1&2&4&1&1&2&4&1&1&2&4\\
\hline
\hline
\end{tabular}\\
\end{table}
\end{proof}

In  Appendix ~\ref{app-8Nhbds}, the use of this Corollary in the evaluation of an 8-valent vertex with neighbourhood 5mMLR or 5MmmR  
 is reported by labelling the three edges 4a, 4b, and 4c and using the consequence that the total inflow across the three edges is at most 4 units.

\begin{theorem}Vertices of degree 8 are not overcharged by this discharging. 
\end{theorem}
\begin{proof} The remaining cases are enumerated in Appendix ~\ref{app-8Nhbds}. 
Outflows from an 8-valent wheel are included with a required inflow producing a net inflow. Rules 81 and 82 produce a net inflow of zero from the possible caps on 555 in 8[555M5],  and Rule 83 produces a net total inflow of at most 2 units from flows crossing  both 55 edges of 8[55M55].

The following arguments are used frequently in the examination of possibly overcharged 8- or 9-valent vertices:\\  
1)  If the flow across an edge is required to be its maximum of 3 units, then the edge must be (m,6) with a 55cap, and if the edge is **, then it must be a 66 edge. \\
2)  If the total flow across the two rim edges of M[x**] is required to be its maximum of 5 units or 4 units, then the middle * must be 6. 
\\
Lemma ~\ref{lem-9[55***5]} in the next section has a useful corollary there but is used occasionally in  Appendix ~\ref{app-8Nhbds}.
\end{proof}

\section{Analysis for 9-valent vertices and Enumeration of 9-neighbourhoods}

9[555555xxx]-10 is reducible so a 9-wheel with nine or eight  5-valent neighbours will not occur in an m5CPG. Since the maximum transfer across a *** or a 5** pair of adjacent edges is 5 units, and the one possible *5* can be split in a pairing of edges,  the final charge on a 9-valent vertex with zero or one 5-valent neighbour is at most \\-30 + 4  + 4$\cdot$5 + 1$\cdot$3 = --3.
It remains to examine first neighbourhoods with seven down to two 5-valent neighbours.

Grouped by the size of the separating ring,  GG first neighbourhoods of a 9-wheel and some other configurations are categorized in Table ~\ref{tabl-9[5mm...5]}. 

\begin{table}[ht]
\caption{GG neighbourhoods of a 9-valent vertex.}
\label{tabl-9[5mm...5]}

\noindent \begin{tabular}{ r | l | l } 
Ring		&Reducible		&Symmetrically D-Irreducible (SDI)\\
\toprule
10:	&9[555555xxx].&\\

11:	&9[555565xxx], 9[555655xxx]. &\\
	&9[5555x555x], 9[55555x505x]. &\\
12:	&9[555665xxx], 9[565565xxx].
	        &9[556565xxx], 9[556655xxx],9[5566505xxx].\\
	&9[5556565xx], 9[5565655xx], 9[5655565xx]. &\\ 
	& 9[5555x565x], 9[5565x555x].&\\
	&9[55565x505x], 9[55555x506x]. 
		&9[55655x505x].\\
&&\\
13:	&
		&9[556665xxx], 9[565665xxx].\\
	&9[5556665xx], 9[5565665xx], 9[5655665xx],
		&9[5566565xx], 9[5566655xx].\\
	& 9[5656565xx].&\\	
	& 9[555556666].&\\
	&9[5565x565x]. &9[5665x555x].\\
	&9[56565x505x], 9[55565x506x], 9[55565x605x], 
		&9[55665x505x], 9[55655x506x].\\
	&9[55555x606x].&\\
&&\\	
14: 	&	&9[5566665xx], 9[5656665xx], 9[5665665xx].\\
	&9[55566665x], 9[55656665x], 9[55665665x],
		&9[55666655x].\\
	& 9[55666565x], 9[56556665x], 9[56565665x], & \\
	&9[56566565x], 9[56655665x].&\\
	&9[555566666], 9[556556666], 9[556655666].&\\
	&9[56565x506x].&9[5665x565x].\\
&&\\
15:	&	&9[5666665xx].\\
	&9[56656665x].
		&9[55666665x], 9[56566665x].\\
\hline
\hline
\end{tabular}
\end{table}
In Table ~\ref{tabl-9-SDI-caps}, to selected SDI configurations are added  5- or 55- caps on a mm rim edge or a 56 cap on a 556 pair of rim edges to produce reducible configurations and a limit on the flows across those rim edges. These caps do not increase the size of the separating ring. 
 9[565x5665x]-14
 is also competed with a 50705 structure that results in a 15-ring configuration. Finally, 9[5mm656xxx] and 9[5666656xx] are not GG but when 656 is expanded to 60506, the results are reducible. 

\begin{table}
\caption{Reducible configurations from SDI first neighbourhoods with  caps added.}
\label{tabl-9-SDI-caps}
\begin{tabular}{r |  r  | l } 
SDI configuration&vertices	&Reducible Capped configurations:\\
\toprule
9[556565xxx]-12:&8&9[5056565],
 9[5560565], 9[5565065]
 \\
9[5565606xxx]-13:&9&9[55656(5)06]  9[556560(5)6]\\
9[556655xxx]-12:&8&9[5506655], 9[5560655]\\

9[556665xxx]-13:&8,9&9[5506665], 9[50516665]
\\

9[565665xxx]-13:&8&9[5605665], 9[5650665]\\

9[5566565xx]-13:&9&9[55660565], 9[55665065]\\
9[5665665xx]-14:&9&9[56650665]
\\
9[55666655x]-14:&10&9[55666655], 9[556660655], 9[556666055],9[556666505]\\
9[565x5665x]-14:& 9,11 & 9[565x56065x], 9[5650705665x]-15\\
9[5mm656xxx]-14 & 9 &9[55560506]-11, 9[55660506]-12, 9[56560506], 9[56660506]-13\\
\hline
\hline
\end{tabular}
\end{table}

To avoid reducible configurations bounded by a 16-ring, many 1-flow targets, numbered 9k30 to 9k74 are used.  Since these are not related to Zero flow targets, they are referred to as 1T30, ..., 1T74.
Similar to Targets 8521 and 8523, several  
1-flow limits are paired into one diagram. For example, 9761 to 9763 show 1-flow limits on each of two 57 edges, limiting the flow to most 2 units across the 5M5 pair of edges of  9[5\{5,6,6\}5M5]. Targets 9660 to 9674 are each demonstrated by at most 4 reducible configurations and sometimes only one. Only sources 21 and 22 deliver charge across a 75 edge with the other shoulder vertex also 5-valent. For 9660, the only compatible 2 unit source across the 75 edge is S21 and 9[5665075] is reducible.
This reducible configuration also enforces the other limit.

Targets 9830 to 9839 are all C(5,2)=10 versions of 9[5\{5$^3$,6$^2$\}LRx]. In the final charge evaluation of neighbourhoods, these are all referred to as 1T30. Targets 9840 to 9859 are all possible versions of 9[5\{5$^2$,6$^2$,LR\}5x] and are referred to as 1T40, 1T47, 1T51, and 1T56.

Recall how targets 8A4a and b combine to limit the flow to 5 units  across three rim edges: LR and the 565 pair. 
Any positive flow across the LR edge is subsumed by the decrease in the flow across the 565 pair which is then at most 2 units.
Target 9606F: T=9[565LR5xxx] combined with several sources produces an SDI configuration which may appear in a m5CPG; so it is labelled FAIL.  Instead, zero flow targets 9606a and 9606b together limit the total flow across the three edges: 565 and LR, to a total of at most 5 units.
Targets 9602a and b, each with one 5cap on a 56 edge give another 3 edge limit.\\
Similar problems occur with uncapped versions of targets 9605, T=9[556LR5xxx] and 9608, T=9[56LR55xxx]. 3 units occurs only for a 55 cap on LR=mm$\neq$55 and these are reducible, but other flows of up to 2 units are possible in a m5CPG. The capped versions, 9605 and 9608, are ZFTs, so any possible 2 units across LR is subsumed by the resulting drop to zero units across the 55 edge. The 55 edge is labelled with its maximum flow of 2 units and an arc to the LR edge labelled Z5 or Z8 indicates that any flow of 1 or 2 units across that LR edge is subsumed. If mm is a possibility for LR, then Z5R3 and Z8R3 would describe the pair of conditions necessary to conclude that the total flow across the 55 and LR edges is at most 2 units.
Uncapped 9602c also is reducible with a 55 cap on LR=mm$\neq$55, and target 9602c shows that the cap on the 55 edge can subsume any 1 or 2-unit flow across LR.\\

\begin{lemma} If the total flow across the four rim edges of M[x***x] is 
required to be its maximum of 10 units or 9 units, then *** must be 666. 
\label{lem-9[55***5]}
\end{lemma}

\begin{figure}
 \centering
 \includegraphics[ width=6in  ]{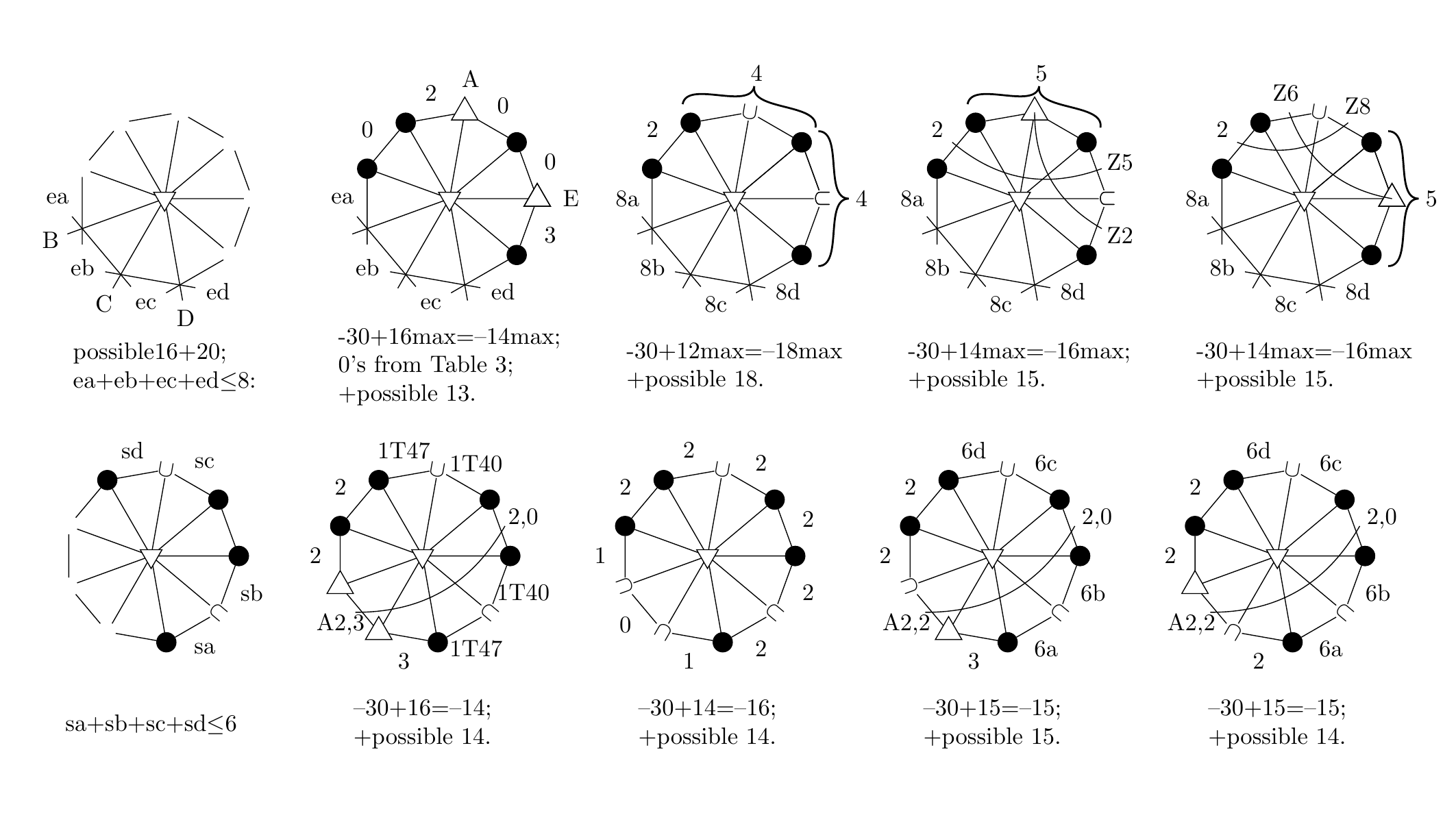}\\
  \caption{Final charge enumeration for 9[55*5*5***]}
\label{fig-9[55*5*5***]}
\end{figure}

\begin{proof} The vertices and flows are labelled in the leftmost diagram of line 1 in Figure ~\ref{fig-9[55*5*5***]}. Each of flows ea, eb, ec and ed is at most 3 units, ea+eb$\leq$5, and ec+ed$\leq$5. Assume ea+eb+ec+ed is 9 or more. Refer to  Figure ~\ref{fig-MmLimits_NegInflow}.
B Major $\implies$ea+eb$\leq$3 so B=6.
D Major $\implies$ec+ed$\leq$3 so D=6.
Now C Major $\implies$ea+eb+ec+ed$\leq$3+1+1+3=8 so C=6.
\end{proof}

\begin{corollary} In 9[55***5xxx], the flow across the four edges 5***5 is at most 8 units.
\end{corollary}

\begin{proof} Refer to the second diagram of line 1 in Figure ~\ref{fig-9[55*5*5***]}. 
If the flow is 9 or more, then *** is 666. While 9[556665]-13 is GG, it is SDI and may appear in a m5CPG. However,
9[5506665]-13 and 9[5560665]-13 are reducible; so without a 5cap, ea+eb$\leq$3 
(Lemma ~\ref{lem-cross-mm}, item 8)
and the flow across all four edges is at most 8 units.
\end{proof}
These edges can be labelled 8a, 8b, 8c, 8d to reflect this limit.  



\begin{corollary} The central 9-valent vertex in the neighbourhoods 9[55*5*5***] is not overcharged.\end{corollary}

\begin{proof}See figure  ~\ref{fig-9[55*5*5***]}, row 1.   Call the two singleton *'s A and E. 
For AE=66, the flows totalling at most 8 units are still labelled ea through ed.   While 9[556565]-12 is SDI, adding  outer 5caps yields 9[5056565]-12, 
9[5560565]-12
and 
9[5565065]-12
, all reducible and giving the limits in the second picture. AE=MM, AE=6M, and AE=M6 complete the set of neighbourhoods.
\end{proof}

Instead of 32 separate cases where the *'s are specified as 6 or Major,  only four diagrams suffice to show that the central 9-valent vertex in these neighbourhoods cannot be overcharged.

\begin{lemma} The total across the four 5M edges of 9[5M55M5] is at most 6 units.
\label{lem-9[5M55M5]}
\end{lemma}

\begin{proof}The flows are labelled in the leftmost diagram of line 2 in Figure ~\ref{fig-9[55*5*5***]}.
Assume sa+sb+sc+sd is 7 or more. Since each of sa,sb,sc and sd is at most 2 units, at least one of sb and sc must be 2 units, say sb. Targets 9672 and 9674 force sc$\leq$1 and sd$\leq$1 and the total is at most 2+2+1+1=6.
\end{proof}

Four cases from 9[5M55M5xxx] are shown. For xxx=566, the 6 unit limit is not sufficient, but the 566 makes the structure match targets from the 9840 and 9847 families to give stronger limits with a  total inflow of at most 4 units. For xxx=5MM, the maximum of 8 units for sa+sb+sc+sd does not lead to overflow, so this lemma is not needed. For xxx=5M6 and 56M, this lemma is a good fit to show the neighbourhood will not overcharge the central 9-valent vertex. All four cases of 5** needed to be examined separately.
\begin{theorem}Vertices of degree 9 are not overcharged by this discharging. 
\end{theorem}
\begin{proof}  The remaining cases are enumerated in Appendix ~\ref{app-9Nhbds}. Limited flow 9-targets are given in Appendix~\ref{app-ZF9}. 
One new trick used here is in the calculation of the precharge from Mayer 2,3,4 contributions. If a first neighbourhood containing 5**5 is reducible when **=66, this reduces the maximum possible Mayer contributions by one unit since one of the *'s must be Major. For example, 9[555**5***] would otherwise list the precharge as --30+16max=--14max, but  9[555665]-12 is reducible, so the percharge is --15max. The maximum inflow across rim edges is at most 15 so this neighbourhood cannot overcharge the central 9-valent vertex. As before, had the maximum inflow across rim edges produced a possible overcharge of at most one unit, then all the cross edge flows must be at their maximums and here the 5**5 must be 56M5 or 5M65.
\end{proof}

\begin{corollary} Every planar graph is 4-colourable.
\end{corollary}

\section{Comparisons with proofs by  RSST and Appel\&Haken}
There are three differences between this and the other two proofs.
The main difference is the discharging and my claim that after going through the text and understanding  the notations, each of the enumerations of the final charge on vertices of degrees 7 to 10 can be verified in a day or less. 

The second is the number of reducible configurations used. A computer and software proving reducibility  is an indispensible factor in providing a sufficient number of configurations to prove a discharging works. 
The Kempe chain closure algorithm is ratcheting, any small error in the input or the algorithm will immediately cascade into a large discrepancy with proven theoretical expectations on D-irreducibility and  also expectations on the openness and closure of sets of boundary colourings \cite{FAPart1}.
Since the correctness of the software can be supported by this evidence of agreement with theoretical expectations, one should be assured  that the reducible configurations claimed in RSST,  A\&H  and here, have been proven by correct programs and are in fact reducible.
The number of times this software has been invoked is then of insignificant importance relative to the accuracy of the software.

The bulk of reducible configurations used in this proof are obtained when the zero flow targets are combined with the compatible sources. The targets are in separate files and the combination with the sources in a fourth file is automated. The  data of  Table ~\ref{tabl_ConfigsTested} has been extracted from the results.
There is no attempt to avoid reducible subconfigurations or duplication of configurations tested. 
For example, if a 5cap on LR=mm$\neq$55of a target is reducible, then the 55 cap is still tested, and on both sides for LR=66. 
In Lemma ~\ref{lem-sixV7Lemmas},
whenever a source of 2 or more units crossing LR is tested, the resulting configuration is the same as one with the 2-unit sources exchanging  positions, but both are constructed and tested.

\begin{table}[ht]
\caption{Categorization of configurations tested for reducibility establishing limited flow Targets of Appendices  ~\ref{app-ZF7},~\ref{app-ZF8} and ~\ref{app-ZF9} 
\label{tabl_ConfigsTested}		}
\noindent \begin{tabular}{ r | r | r | r | r}

\multicolumn{2}{r|}{ }&7-Targets&8-Targets&9-Targets\\
\hline
\multicolumn{2}{r|}{Time (minutes)}&		9&	22&	12\\
\multicolumn{2}{r|}{Configurations Tested}&	3030&		5045&	2052\\
\toprule
Ring Size	&	8&		2&		0&		0\\
		&	9&		20&		3&	 	0\\
		&	10&		86&		40&		3\\
		&	11&		291&		245&		41\\
		&	12&		584&		780&		200\\
		&	13&		895&		1530&	545\\
		&	14&		855&		1743&	831\\
		&	15&		297&		704&		432\\
		&	16&		0&		0&		0\\
\hline
Reducibility&	D&		2634&	4107	&	1647	\\
Type&		SDI&		2&		3&		1\\
		&	E1&		79&		150&		69\\
		&	E5&		232&		 641&	287\\
		&	Ee&		42&		86&		23\\
		&	E8&		2&		7&		0\\
		&	E9&		0&		0&		6\\
		&	C1&		0&		21&		0\\
		&	C2&		34&		4&		0\\
\multicolumn{2}{r|}{needs reducer finding}	&		5&		26&		19\\
\hline
\hline
\end{tabular}
\end{table}

A third difference is the number of discharging rules. RSST's proof has 32  while A\&H used 487 (\cite{RSST}page 18). Appendix G shows  56 rules, but rules 2 to 7 and 62 are multiple rules, giving an extra 1+3+3+1+1+1=10 for a total of 66 fully specified rules where each non-sink vertex in the source structure has a specified degree. Similar counting of minor vertices (6-valent with a "--") in RSST gives an extra 1+3+3+1+1+1+1+3+1+1=16, so their list is at least 48 fully specified charge transfer rules. 

A final difference is the absence of an algorithm to produce a 4-colouring. This is deliberate. My goal was to find a discharging that can be verified by a human, completing ~\cite{FAPart1}. 


\cleardoublepage

\appendix


\section*{Supplementary material}

## Appendix B. Enumeration of Neighbourhoods of a 7-valent Vertex

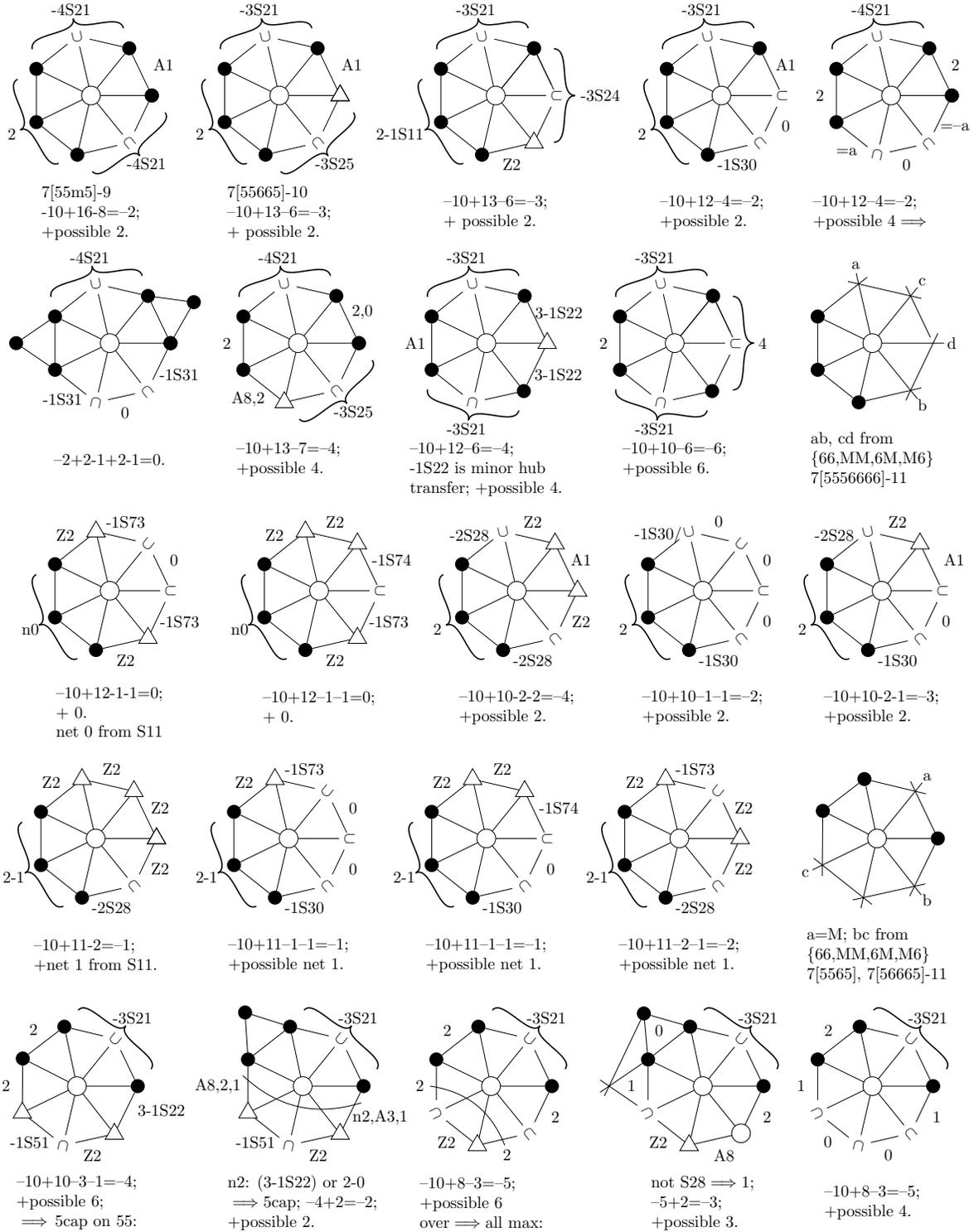

7[55m5]-9
−10+16-8=−2;
+possible 2.

7[55665]-10
−10+13-6=−3;
+ possible 2.

−10+13-6=−3;
+ possible 2.

−10+12-4=−2;
+possible 2.

−10+12-4=−2;
+possible 4 ⟹

−2+2-1+2-1=0.

−10+13-7=−4;
+possible 4.

−10+12-6=−4;
−1S22 is minor hub
transfer; +possible 4.

−10+10-6=−6;
+possible 6.

ab, cd from
{66,MM,6M,M6}
7[5556666]-11

−10+12-1-1=0;
+ 0.
net 0 from S11.

−10+12-1-1=0;
+ 0.

−10+10-2-2=−4;
+possible 2.

−10+10-1-1=−2;
+possible 2.

−10+10-2-1=−3;
+possible 2.

−10+11-2=−1;
+net 1 from S11.

−10+11-1-1=−1;
+possible net 1.

−10+11-1-1=−1;
+possible net 1.

−10+11-2-1=−2;
+possible net 1.

a=M; bc from
{66,MM,6M,M6}
7[5565], 7[56665]-11

−10+10-3-1=−4;
+possible 6;
⟹ 5cap on 55:

n2: (3-1S22) or 2-0
⟹ 5cap; −4+2=−2;
+possible 2.

−10+8-3=−5;
+possible 6
over ⟹ all max:

not S28 ⟹ 1;
−5+2=−3;
+possible 3.

−10+8-3=−5;
+possible 4.

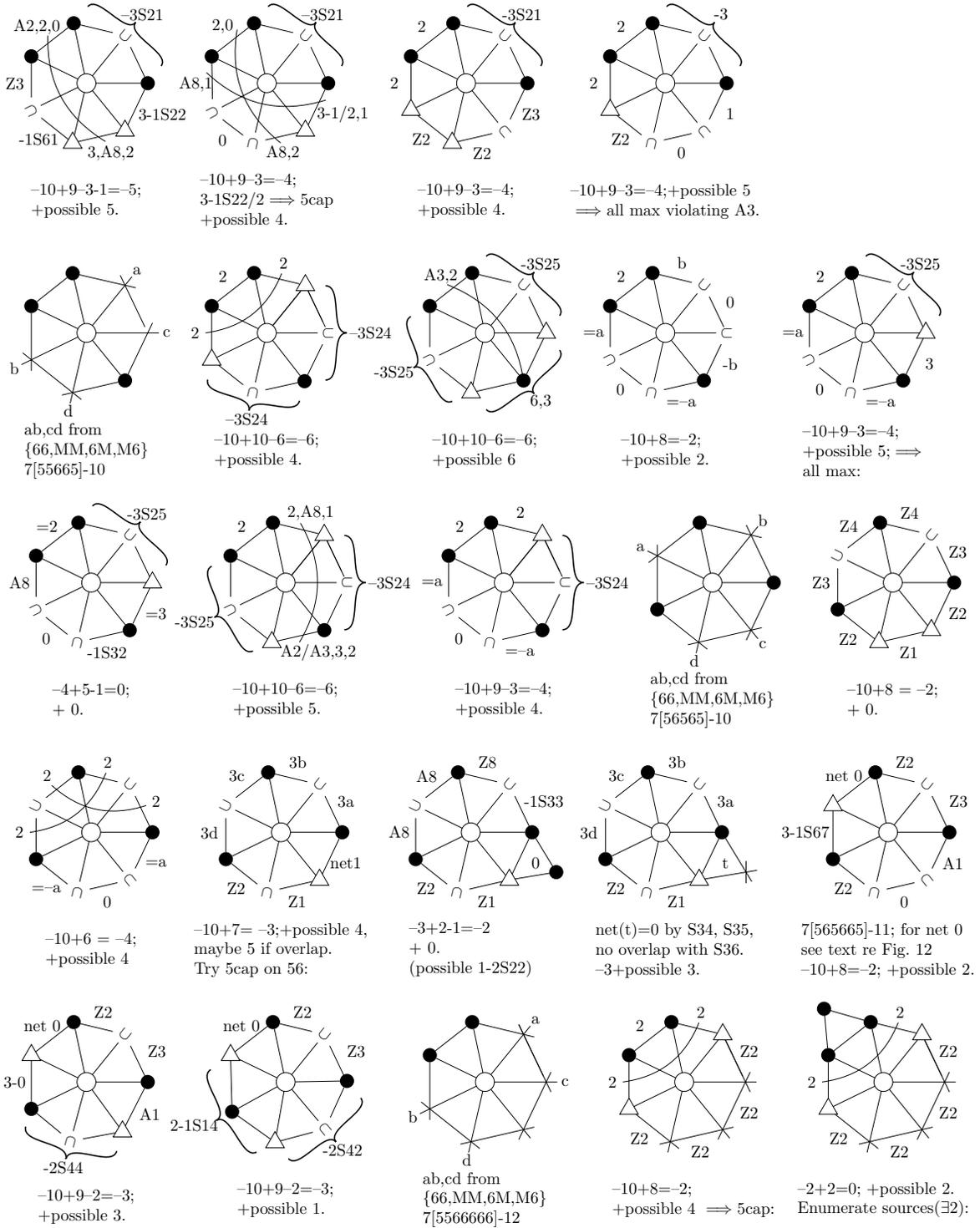

−10+9−3−1=−5;
+possible 5.

−10+9−3=−4;
3−1S22/2 ⟹ 5cap
+possible 4.

−10+9−3=−4;
+possible 4.

−10+9−3=−4;+possible 5
⟹ all max violating A3.

ab,cd from
{66,MM,6M,M6}
7[55665]−10

−10+10−6=−6;
+possible 4.

−10+10−6=−6;
+possible 6

−10+8=−2;
+possible 2.

−10+9−3=−4;
+possible 5; ⟹
all max:

−4+5−1=0;
+ 0.

−10+10−6=−6;
+possible 5.

−10+9−3=−4;
+possible 4.

ab,cd from
{66,MM,6M,M6}
7[55665]−10

−10+8 = −2;
+ 0.

−10+6 = −4;
+possible 4

−10+7=−3;+possible 4,
maybe 5 if overlap.
Try 5cap on 56:

−3+2−1=−2
+ 0.
(possible 1−2S22)

net(t)=0 by S34, S35,
no overlap with S36.
−3+possible 3.

7[565665]−11; for net 0
see text re Fig. 12
−10+8=−2; +possible 2.

−10+9−2=−3;
+possible 3.

−10+9−2=−3;
+possible 1.

ab,cd from
{66,MM,6M,M6}
7[5566666]−12

−10+8=−2;
+possible 4 ⟹ 5cap:

−2+2=0; +possible 2.
Enumerate sources(∃2):

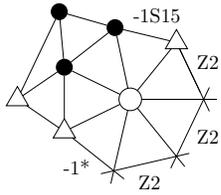

0+2−1−1=0;
−1* from S17 or S71.

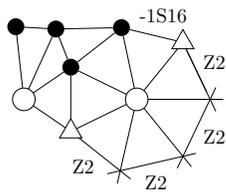

0+1−1=0; + 0.
5[71515176]−11
Else 0 across 65.

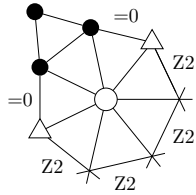

0+ 0.

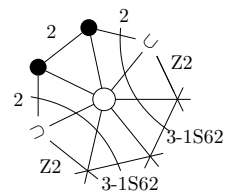

−10+6=−4;+possible 6:
⟹ 5cap on 55:

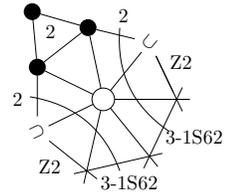

−4+2=−2; Try sources
across M55; ∃5: T=7,
Z2, shoulders{*5,*}:

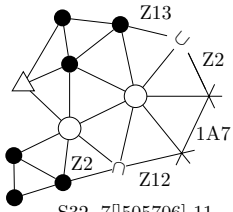

S32, 7[]505706]−11
−2+1=−1
+possible 1.

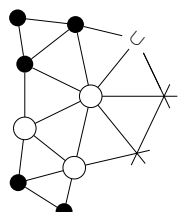

Sources 54,57
77[55-7-505]−12

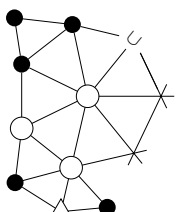

Source 58
77[565-7-55]−13

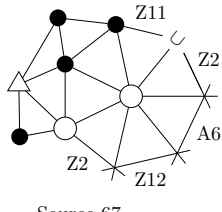

Source 67
−2+1=−1; + 0.
Else zero across 5M.

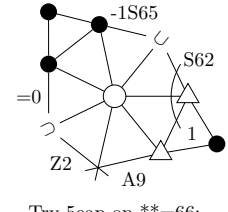

Try 5cap on **=66:
−2+2-1S65=−1;
+possible 1-1S62.

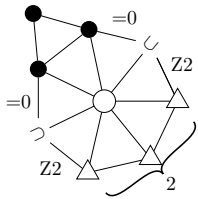

666, no 5cap⟹max=2
−2+possible 2.
Else *** not 666

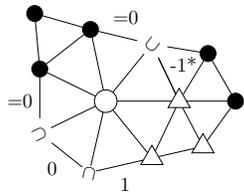

−2+2-1*S63=−1
+possible 1.

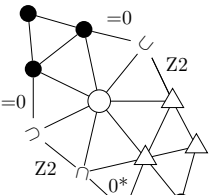

0*: No source
7[506TM]. −2+2=0
+ 0.

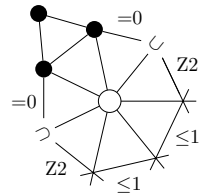

No 2 across any **
−2+possible 1+1.

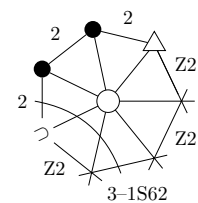

−10+7=−3;
+possible 6
Try *cap on 55:

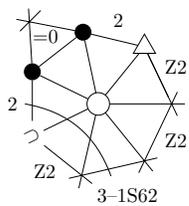

−3;+possible 4
⟹ 5cap on 56:

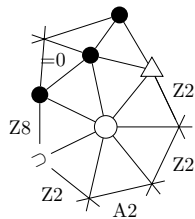

−3+2=−1;
+ 0.
Else 5cap on 55.

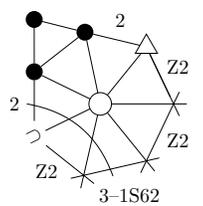

−3+2=−1;+possible 4
enumerate cross 56;∃3:

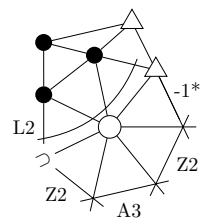

−1* from S71 or S17
−1+2-1*=0;
+ 0.

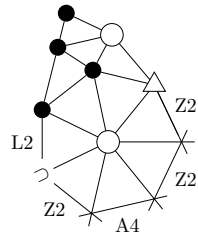

−1+1=0;
+ 0.

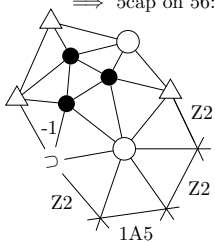

−1+1-1S66=−1;
+possible 1.
Else 0 across 56:

Similar to 7[M505M***]
above, enumerate
sources across M55.
Sources 54, 57
and 58 yield the
same reducible
subconfigurations
as above.

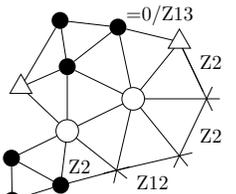

Source 32
−1+1=0
+ 0.

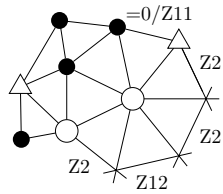

Source 67: −1+1=0.
+ 0.
Else 0 across 5M.

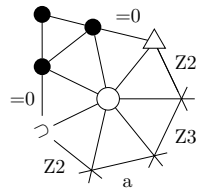

−1+possible a
⟹ a≥2, crossing 66:

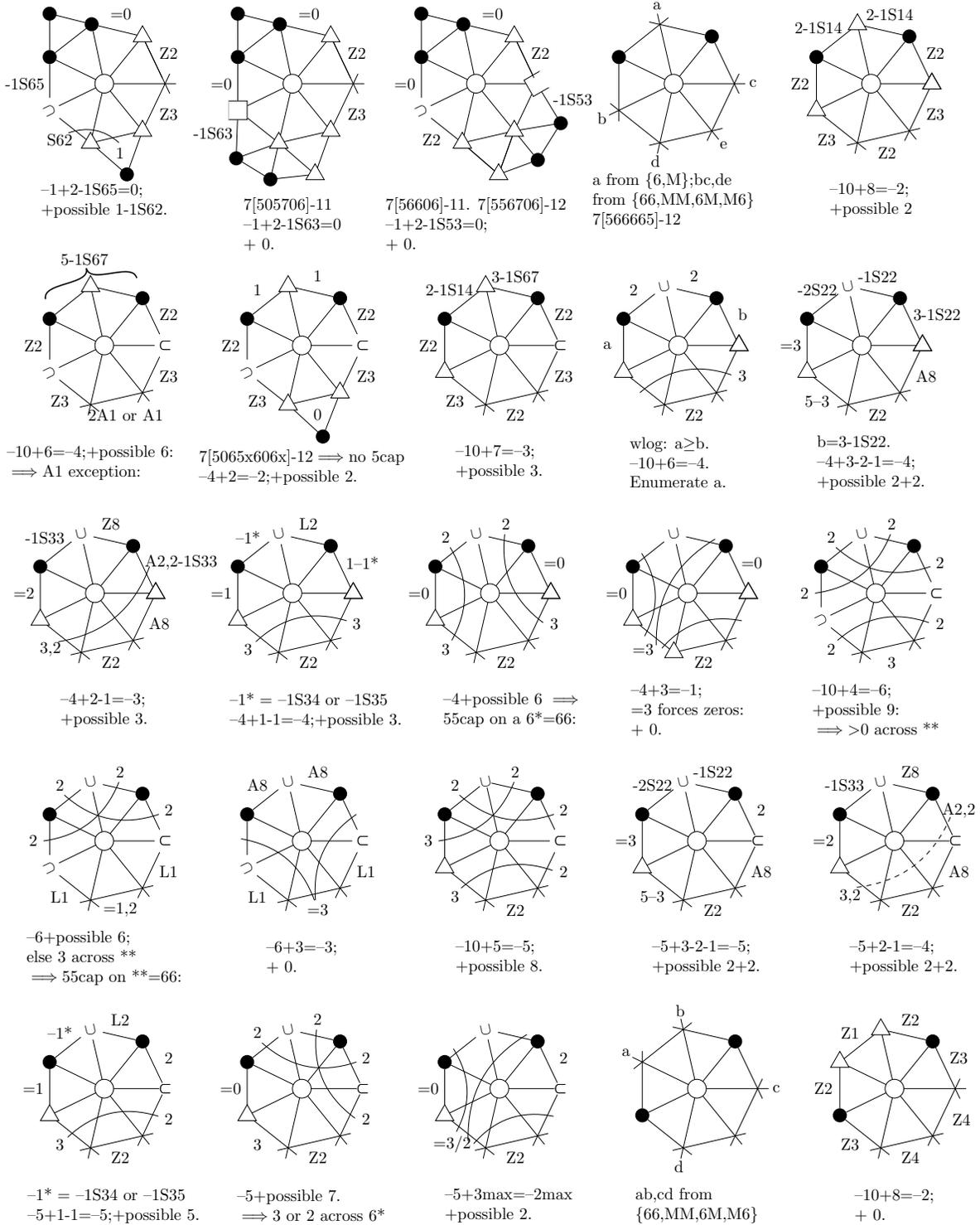

$-1+2-1S65=0$;
+possible $1-1S62$.

$7[505706]-11$
$-1+2-1S63=0$
$+0$.

$7[56606]-11.$ $7[556706]-12$
$-1+2-1S53=0$;
$+0$.

a from $\{6,M\}$;bc,de
from $\{66,MM,6M,M6\}$
$7[566665]-12$

$-10+8=-2$;
+possible 2

$-10+6=-4$;+possible 6:
$\Longrightarrow$ A1 exception:

$7[5065x606x]-12\Longrightarrow$ no 5cap
$-4+2=-2$;+possible 2.

$-10+7=-3$;
+possible 3.

wlog: a≥b.
$-10+6=-4$.
Enumerate a.

b=3-1S22.
$-4+3-2-1=-4$;
+possible 2+2.

$-4+2-1=-3$;
+possible 3.

$-1^* = -1S34$ or $-1S35$
$-4+1-1=-4$;+possible 3.

$-4+$possible 6 $\Longrightarrow$
55cap on a 6*=66:

$-4+3=-1$;
=3 forces zeros:
$+0$.

$-10+4=-6$;
+possible 9:
$\Longrightarrow$ >0 across **

$-6+$possible 6;
else 3 across **
$\Longrightarrow$ 55cap on **=66:

$-6+3=-3$;
$+0$.

$-10+5=-5$;
+possible 8.

$-5+3-2-1=-5$;
+possible 2+2.

$-5+2-1=-4$;
+possible 2+2.

$-1^* = -1S34$ or $-1S35$
$-5+1-1=-5$;+possible 5.

$-5+$possible 7.
$\Longrightarrow$ 3 or 2 across 6*

$-5+3max=-2max$
+possible 2.

ab,cd from
$\{66,MM,6M,M6\}$

$-10+8=-2$;
$+0$.

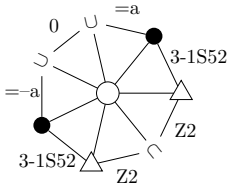

7[56665]-11
-10+6=-4;
+possible 4.

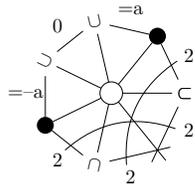

-10+4=-6;
+possible 4.

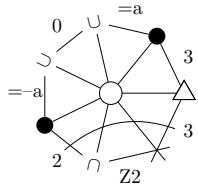

-10+5=-5;+possible 6
$\Longrightarrow$ 3+3 violating limit
of 5 units across 56*.

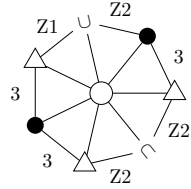

-10+7=-3;
+possible 9:

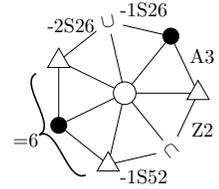

-3+6-4=-1;
+ 0.

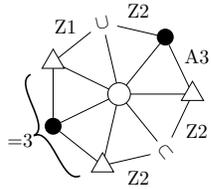

-3+3=0;
+ 0.

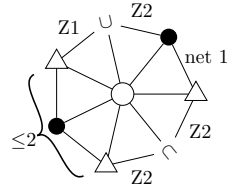

net 1 by S45
-3+possible 3.

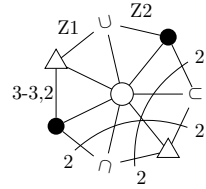

-10+5=-5;
+possible 6 $\Longrightarrow$
lone 5cap on 56;

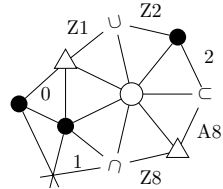

-5+2=-3;
+possible 3.

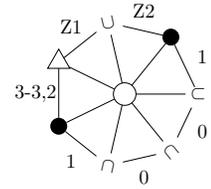

-10+5=-5;
+possible 4.

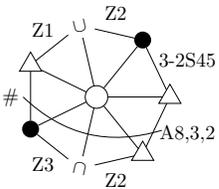

#=3-1S26/2,A2,1
-10+6=-4;
+possible 3+1.

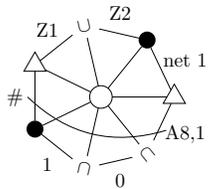

#=3-1S26/2,1
-10+6=-4;
+possible 2+1+1.

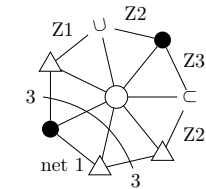

-10+6=-4;
+possible 4.
net1 uses S43.

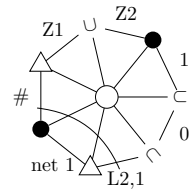

#=3-1S26/2,1
-10+6=-4;
+possible 1+2+1.

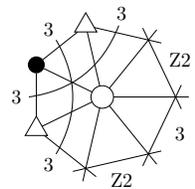

-10+4=-6;
+possible 9:

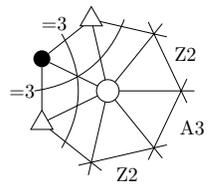

-6+6=0;
+ 0.

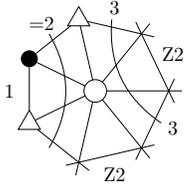

-6+2=-4;
+possible 4.

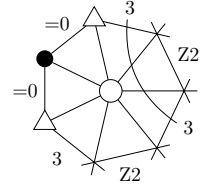

-6;+possible 6.

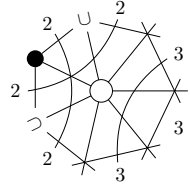

-10+2=-8;
+possible 10:
$\Longrightarrow$ at least one 3

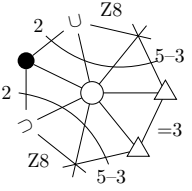

-8+3=-5;+possible 4;
else ≤2 here and
3 across another **

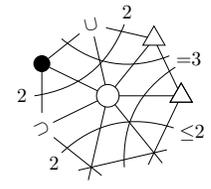

-8+3=-5;
+possible 4.

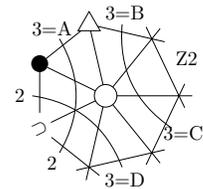

-10+3=-7;
+possible 9 $\Longrightarrow$ at least
2 of A,D,BorC equals 3:

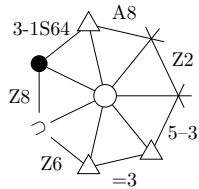

D=3:-7+3=-4;
+possible 4; else D≤2 $\Longrightarrow$
A=3 $\Longrightarrow$ B≤2 $\Longrightarrow$ C=3:

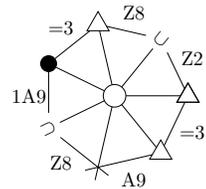

-7+3+3=-1;
+possible 1.

7[**** ***]:
see Figure 16.

## Appendix C.  8-valent Targets with Zero- or One-Flow Limits.

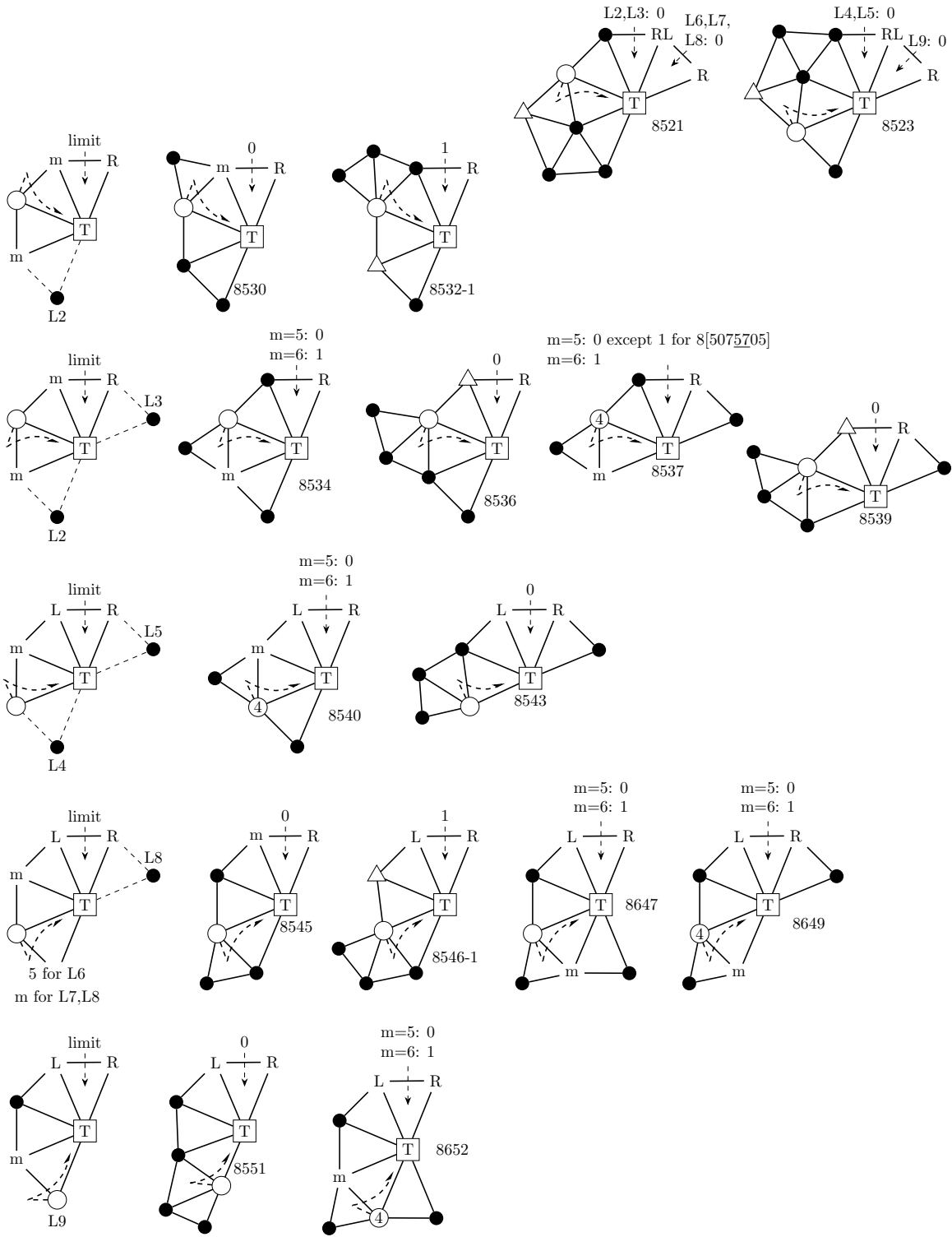

## Appendix D. Enumeration of Neighbourhoods of an 8-valent Vertex

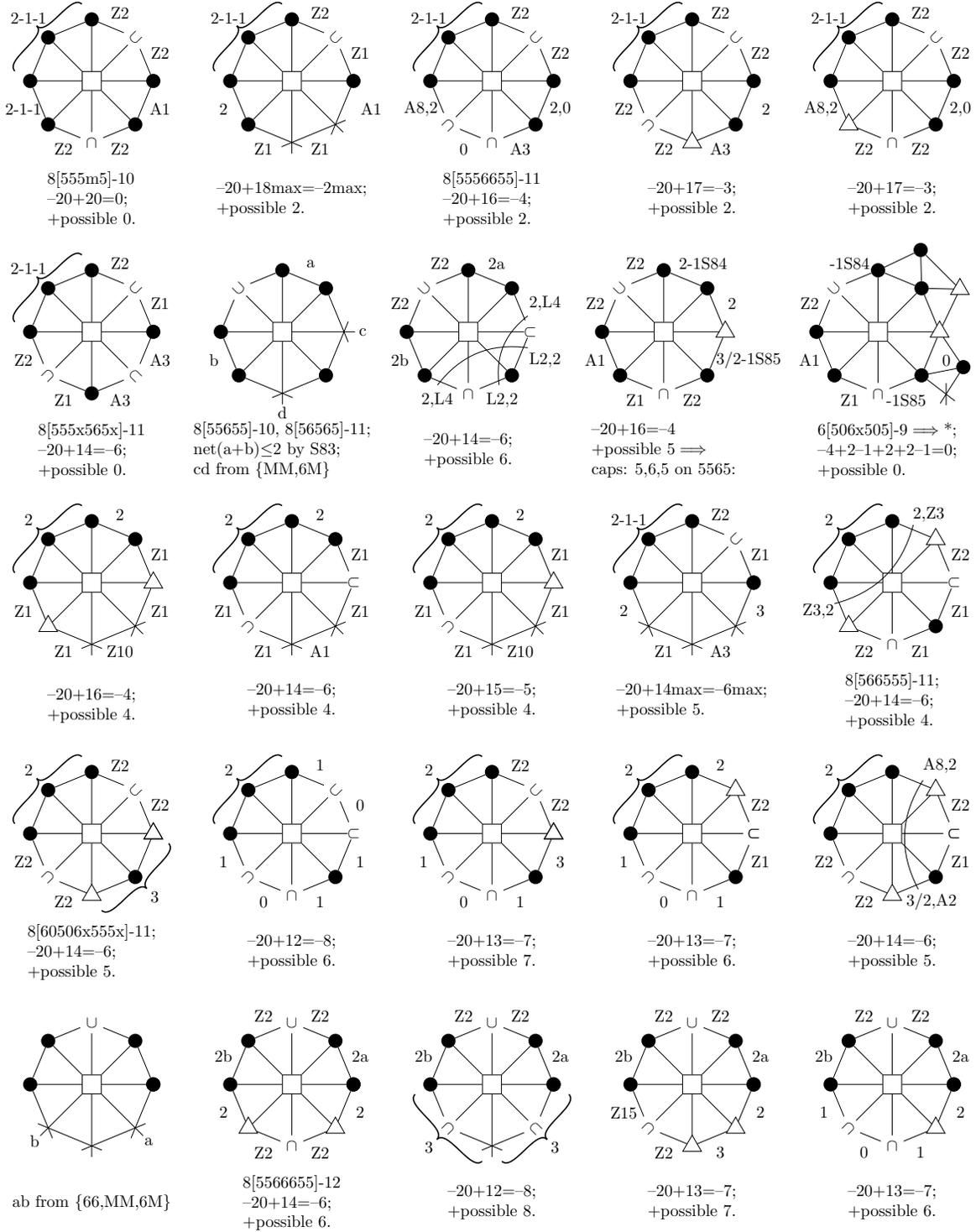

8[555m5]-10
−20+20=0;
+possible 0.

−20+18max=−2max;
+possible 2.

8[5556655]-11
−20+16=−4;
+possible 2.

−20+17=−3;
+possible 2.

−20+17=−3;
+possible 2.

8[555x565x]-11
−20+14=−6;
+possible 0.

8[55655]-10, 8[56565]-11;
net(a+b)≤2 by S83;
cd from {MM,6M}

−20+14=−6;
+possible 6.

−20+16=−4
+possible 5 ⟹
caps: 5,6,5 on 5565:

6[506x505]-9 ⟹ *;
−4+2−1+2+2−1=0;
+possible 0.

−20+16=−4;
+possible 4.

−20+14=−6;
+possible 4.

−20+15=−5;
+possible 4.

−20+14max=−6max;
+possible 5.

8[566555]-11;
−20+14=−6;
+possible 4.

8[60506x555x]-11;
−20+14=−6;
+possible 5.

−20+12=−8;
+possible 6.

−20+13=−7;
+possible 7.

−20+13=−7;
+possible 6.

−20+14=−6;
+possible 5.

ab from {66,MM,6M}

8[5566655]-12
−20+14=−6;
+possible 6.

−20+12=−8;
+possible 8.

−20+13=−7;
+possible 7.

−20+13=−7;
+possible 6.

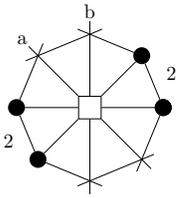

use 4-fold symmetry,
ab from {66,MM,6M}

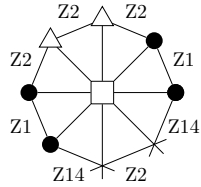

−20+16max=−4max;
+possible 0

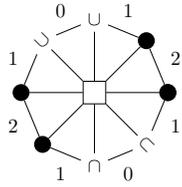

−20+12=−8
+possible 8.

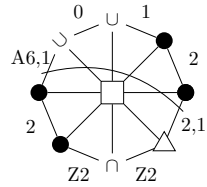

−20+13=−7;
+possible 7.

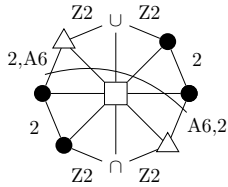

−20+14=−6;
+possible 6

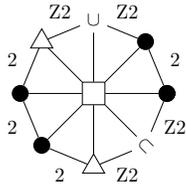

−20+14=−6;
+possible 8 ⟹ 5caps,
6cap on at least one 56

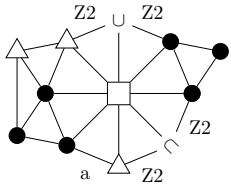

−6+6+a=a; a> 0 ⟹
8[6150516x505x]-12 or
8[62((6))50516x505]-13.

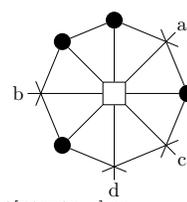

8[565565xx]-11
ab from {MM,6M}
cd from {66,MM,6M,M6}

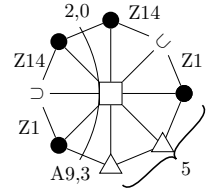

−20+12=−8;
+possible 8.

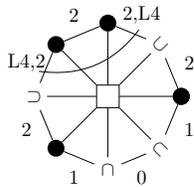

−20+10=−10;
+possible 10.

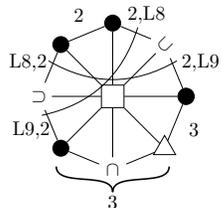

−20+11=−9;
+possible 12
try 5cap on 56:

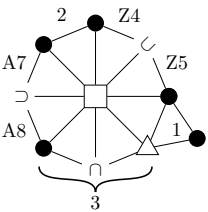

−9+2=−7;+possible 6.
Else *,−9+possible10 ⟹

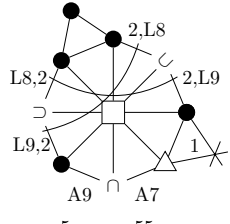

5cap on 55:
−9+2=−7;
+possible 5.

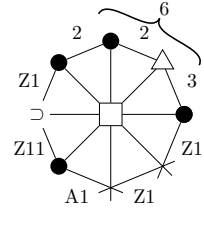

−20+14max=−6max;
+possible 6.

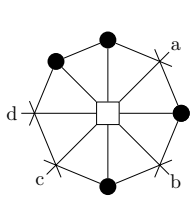

8[56565xxx]-11,
8[5655665x], 8[565565xx]-12
Try a=6,b=6,d=6

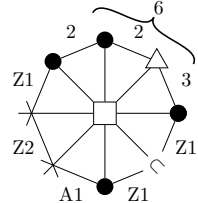

−20+14max=−6max;
+possible 6.

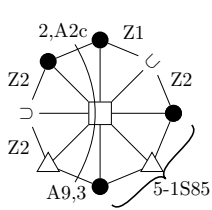

−20+13=−7
+possible 7.

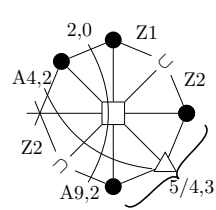

−20+13max=−7max;
5/4 ⟹ 5cap;
+possible 7.

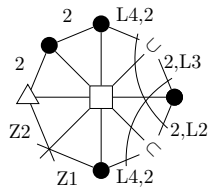

−20+12max=−8max
+possible 8.

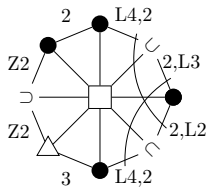

not L3 exception:
−20+11=−9
+possible 9.

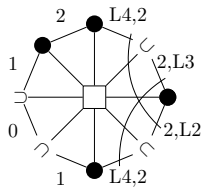

−20+10=−10
+possible 8.

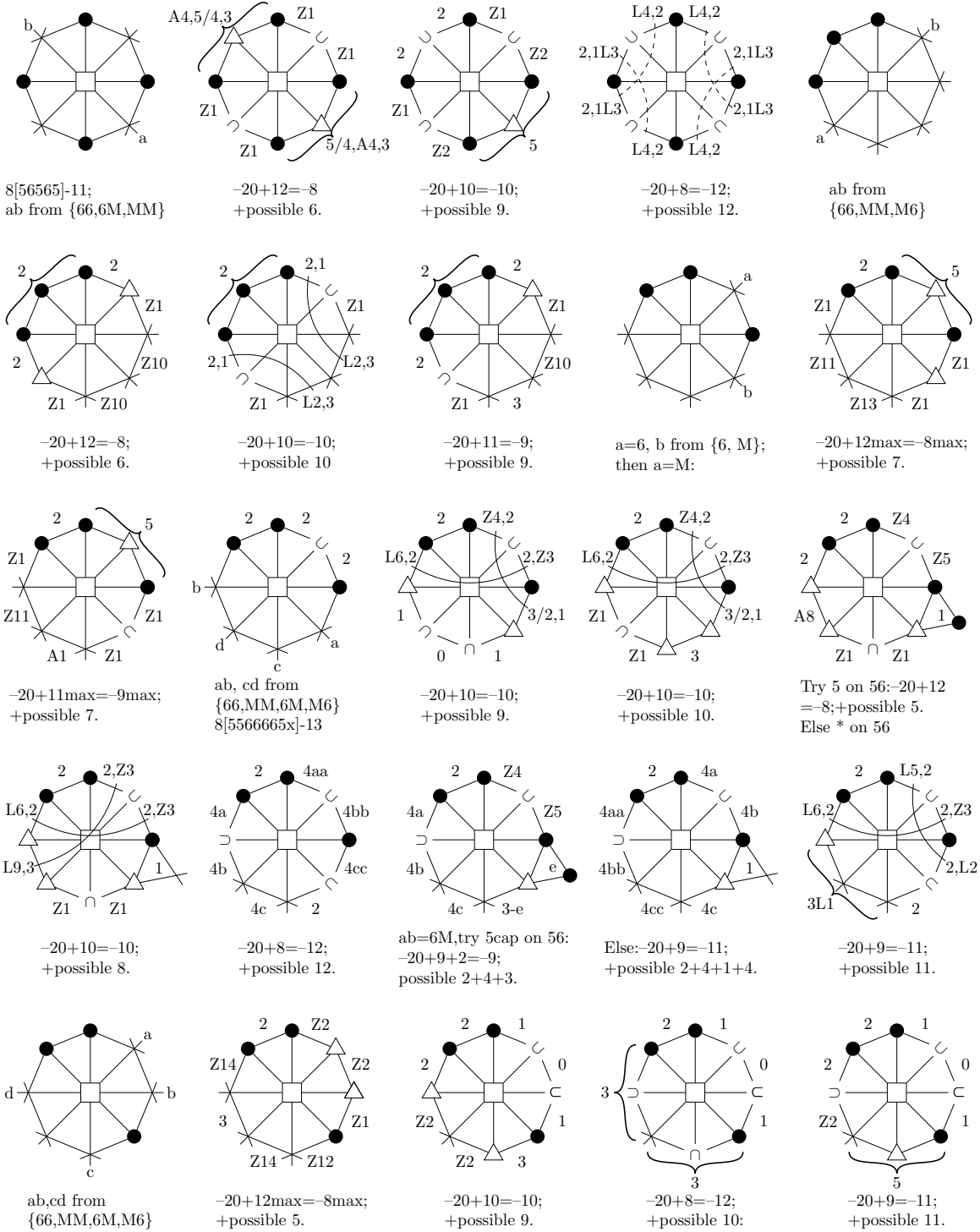

8[56565]-11;
ab from {66,6M,MM}

−20+12=−8
+possible 6.

−20+10=−10;
+possible 9.

−20+8=−12;
+possible 12.

ab from
{66,MM,M6}

−20+12=−8;
+possible 6.

−20+10=−10;
+possible 10

−20+11=−9;
+possible 9.

a=6, b from {6, M};
then a=M:

−20+12max=−8max;
+possible 7.

−20+11max=−9max;
+possible 7.

ab, cd from
{66,MM,6M,M6}
8[5566665x]-13

−20+10=−10;
+possible 9.

−20+10=−10;
+possible 10.

Try 5 on 56:−20+12
=−8;+possible 5.
Else * on 56

−20+10=−10;
+possible 8.

−20+8=−12;
+possible 12.

ab=6M,try 5cap on 56:
−20+9+2=−9;
possible 2+4+3.

Else:−20+9=−11;
+possible 2+4+1+4.

−20+9=−11;
+possible 11.

ab,cd from
{66,MM,6M,M6}

−20+12max=−8max;
+possible 5.

−20+10=−10;
+possible 9.

−20+8=−12;
+possible 10:

−20+9=−11;
+possible 11.

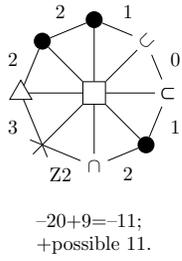

−20+9=−11;
+possible 11.

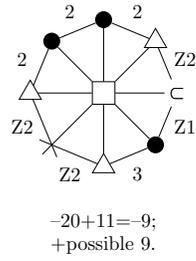

−20+11=−9;
+possible 9.

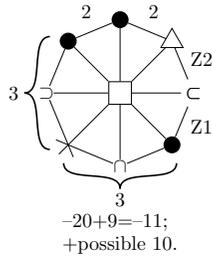

−20+9=−11;
+possible 10.

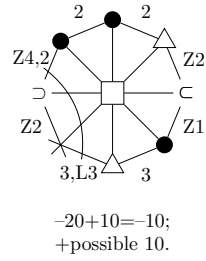

−20+10=−10;
+possible 10.

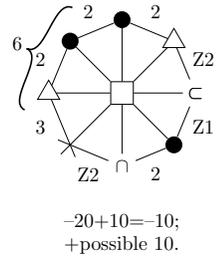

−20+10=−10;
+possible 10.

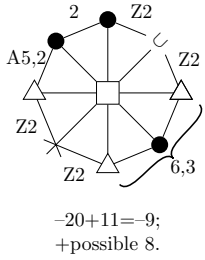

−20+11=−9;
+possible 8.

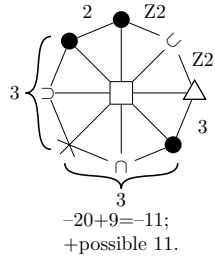

−20+9=−11;
+possible 11.

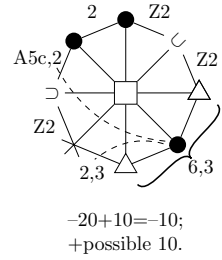

−20+10=−10;
+possible 10.

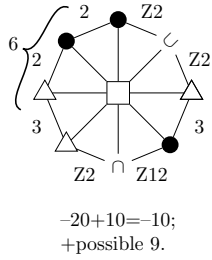

−20+10=−10;
+possible 9.

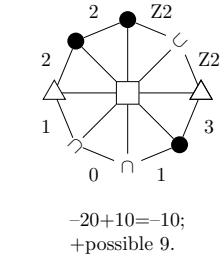

−20+10=−10;
+possible 9.

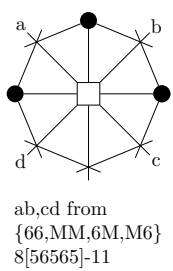

ab,cd from
{66,MM,6M,M6}
8[56565]-11

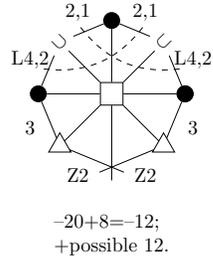

−20+8=−12;
+possible 12.

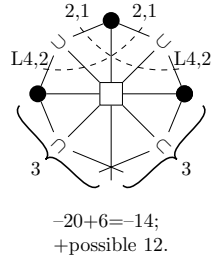

−20+6=−14;
+possible 12.

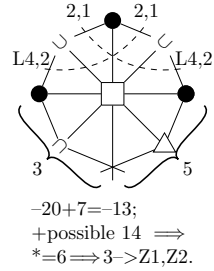

−20+7=−13;
+possible 14 ⟹
*=6 ⟹ 3−>Z1,Z2.

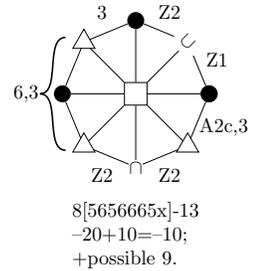

8[5656665x]-13
−20+10=−10;
+possible 9.

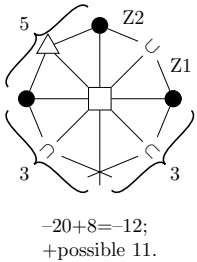

−20+8=−12;
+possible 11.

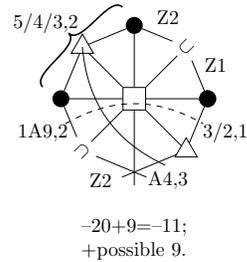

−20+9=−11;
+possible 9.

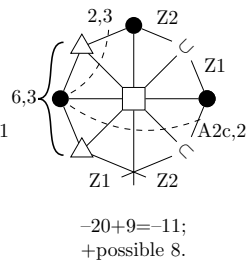

−20+9=−11;
+possible 8.

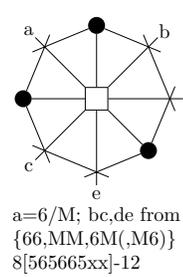

a=6/M; bc,de from
{66,MM,6M(,M6)}
8[565665xx]-12

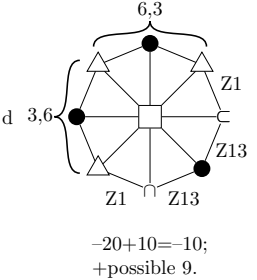

−20+10=−10;
+possible 10.

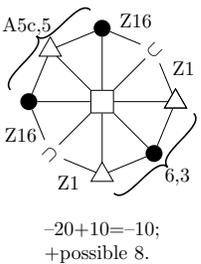

−20+10=−10;
+possible 8.

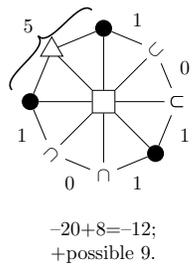

−20+8=−12;
+possible 9.

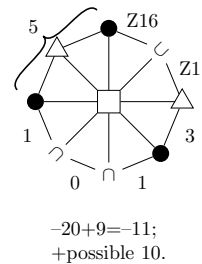

−20+9=−11;
+possible 10.

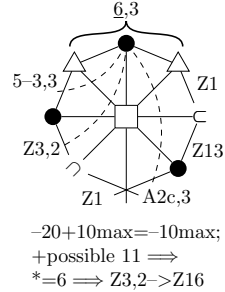

−20+10max=−10max;
+possible 11 ⟹
*=6 ⟹ Z3,2−>Z16

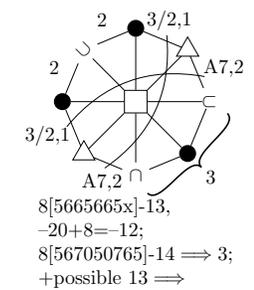

8[5665665x]-13,
−20+8=−12;
8[567050765]-14 ⟹ 3;
+possible 13 ⟹

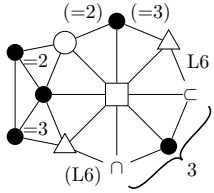

$=2 \implies$ S21 $\implies$ L6 $\implies$ =3
$\implies$ 5cap on 56;(=2) $\implies$ ...
$\implies$ 8[605070506]-11.

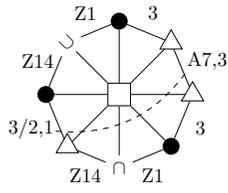

$-20+9=-11;$
+possible 3+3+4.

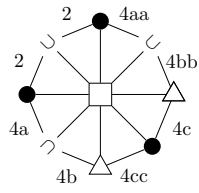

$-20+8=-12;$
+possible 12.

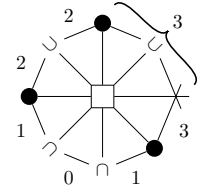

combine de=MM,6M
$-20+7max=-13max;$
+possible 12

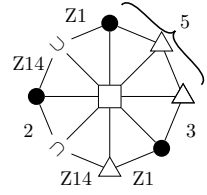

$-20+9=-11;$
+possible 10

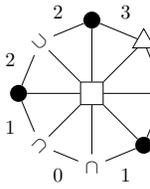

$-20+7=-13;$
+possible 13.

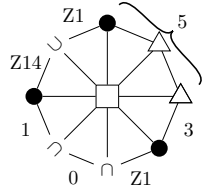

$-20+8=-12;$
+possible 9.

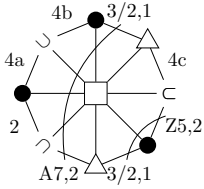

$-20+8=-12;$
+possible 12.

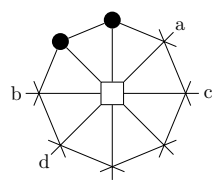

ab,cd from
{66,MM,6M,M6}

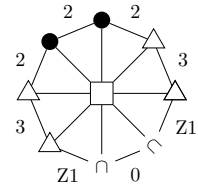

8[55666666]-14
$-20+8=-12;$
+possible 12.

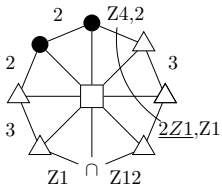

$-20+8=-12;$
+possible 12.
2Z1 from 606

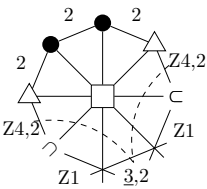

$-20+8=-12;$
+possible 12.

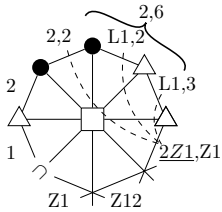

$-20+8=-12;$
+possible 9.

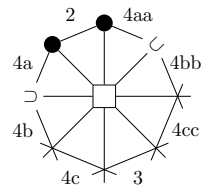

cd=**
$-20+6=-14$
+possible 13.

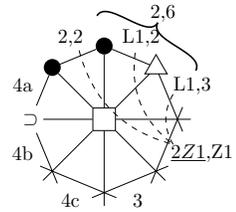

cd=**
$-20+7=-13;$
+possible 13.

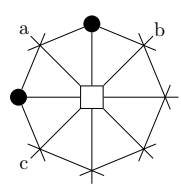

a from {6,M}, bc
from{66,MM,6M}

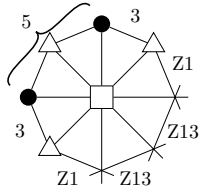

$-20+8=-12;$
+possible 11.

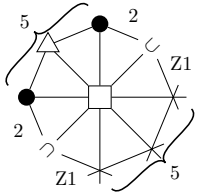

$-20+6=-14;$
+possible 14.

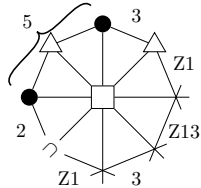

$-20+7=-13;$
+possible 13

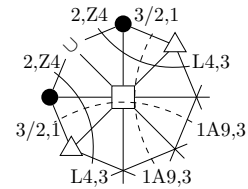

$-20+6=-14;$
+possible 3+4+4+3.

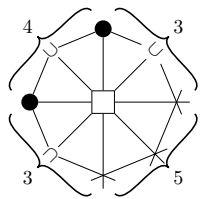

$-20+4=-16;$
+possible 15.

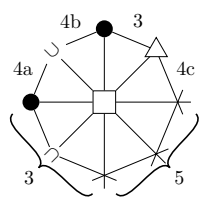

$-20+5=-15;$
+possible 15.

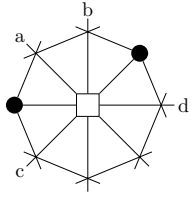

ab,cd from
{66,MM,6M,M6}

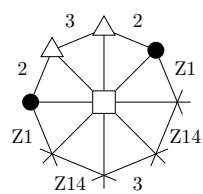

$-20+8max=-12max;$
+possible10.

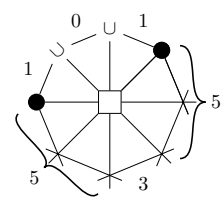

$-20+6max=-14max;$
+possible 15 $\implies$
8[5660665] from 8616

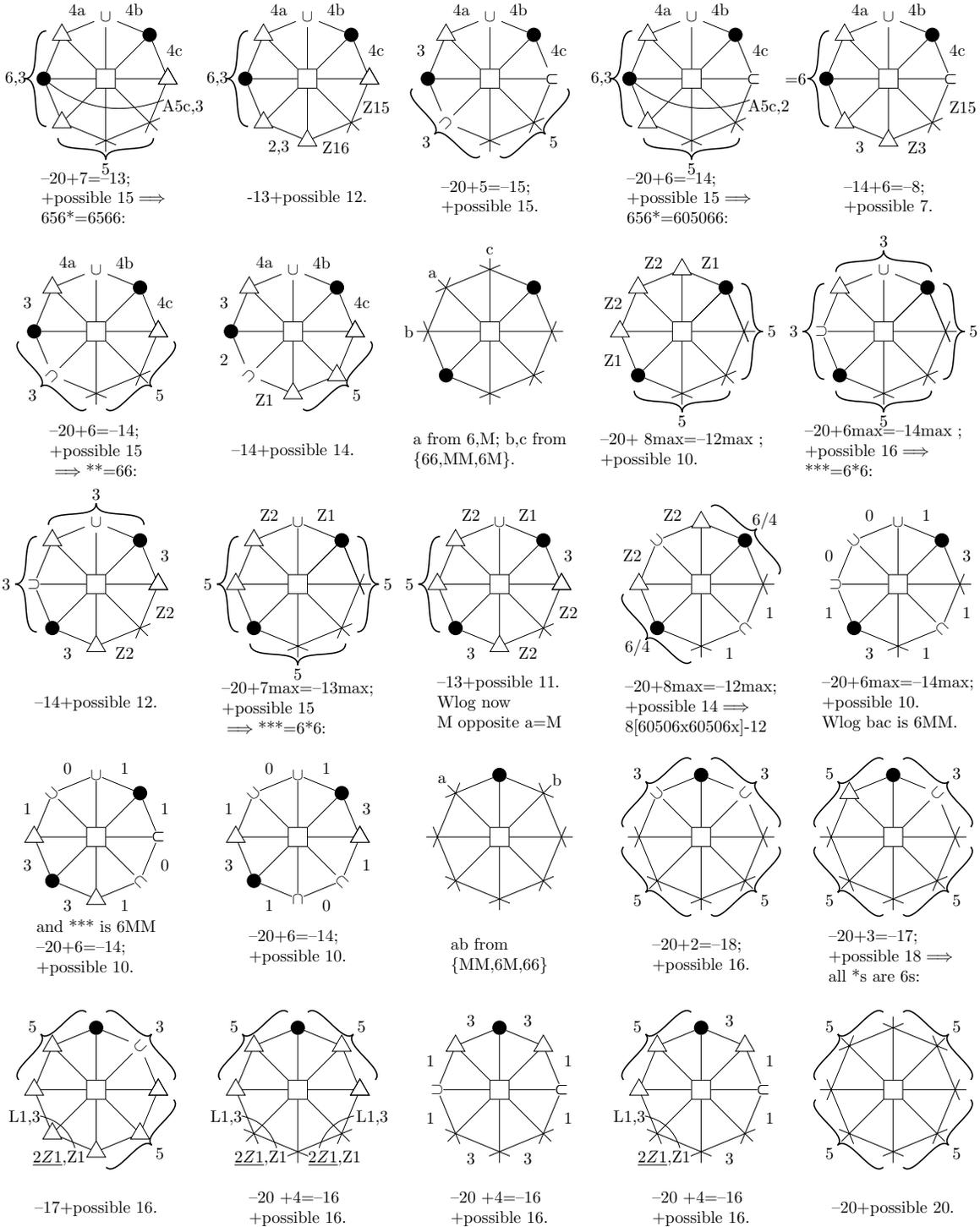

$-20+7=13$;
+possible 15 $\implies$
656*=6566:

$-13$+possible 12.

$-20+5=-15$;
+possible 15.

$-20+6=14$;
+possible 15 $\implies$
656*=605066:

$-14+6=-8$;
+possible 7.

$-20+6=-14$;
+possible 15
$\implies$ **=66:

$-14$+possible 14.

a from 6,M; b,c from
{66,MM,6M}.

$-20+$ 8max$=-12$max ;
+possible 10.

$-20+6$max$=-14$max ;
+possible 16 $\implies$
***=6*6:

$-14$+possible 12.

$-20+7$max$=-13$max;
+possible 15
$\implies$ ***=6*6:

$-13$+possible 11.
Wlog now
M opposite a=M

$-20+8$max$=-12$max;
+possible 14 $\implies$
8[60506x60506x]-12

$-20+6$max$=-14$max;
+possible 10.
Wlog bac is 6MM.

and *** is 6MM
$-20+6=-14$;
+possible 10.

$-20+6=-14$;
+possible 10.

ab from
{MM,6M,66}

$-20+2=-18$;
+possible 16.

$-20+3=-17$;
+possible 18 $\implies$
all *s are 6s:

$-17$+possible 16.

$-20$ $+4=-16$
+possible 16.

$-20$ $+4=-16$
+possible 16.

$-20$ $+4=-16$
+possible 16.

$-20$+possible 20.

## APPENDIX E. 9-VALENT TARGETS WITH ZERO- OR ONE-FLOW LIMITS.

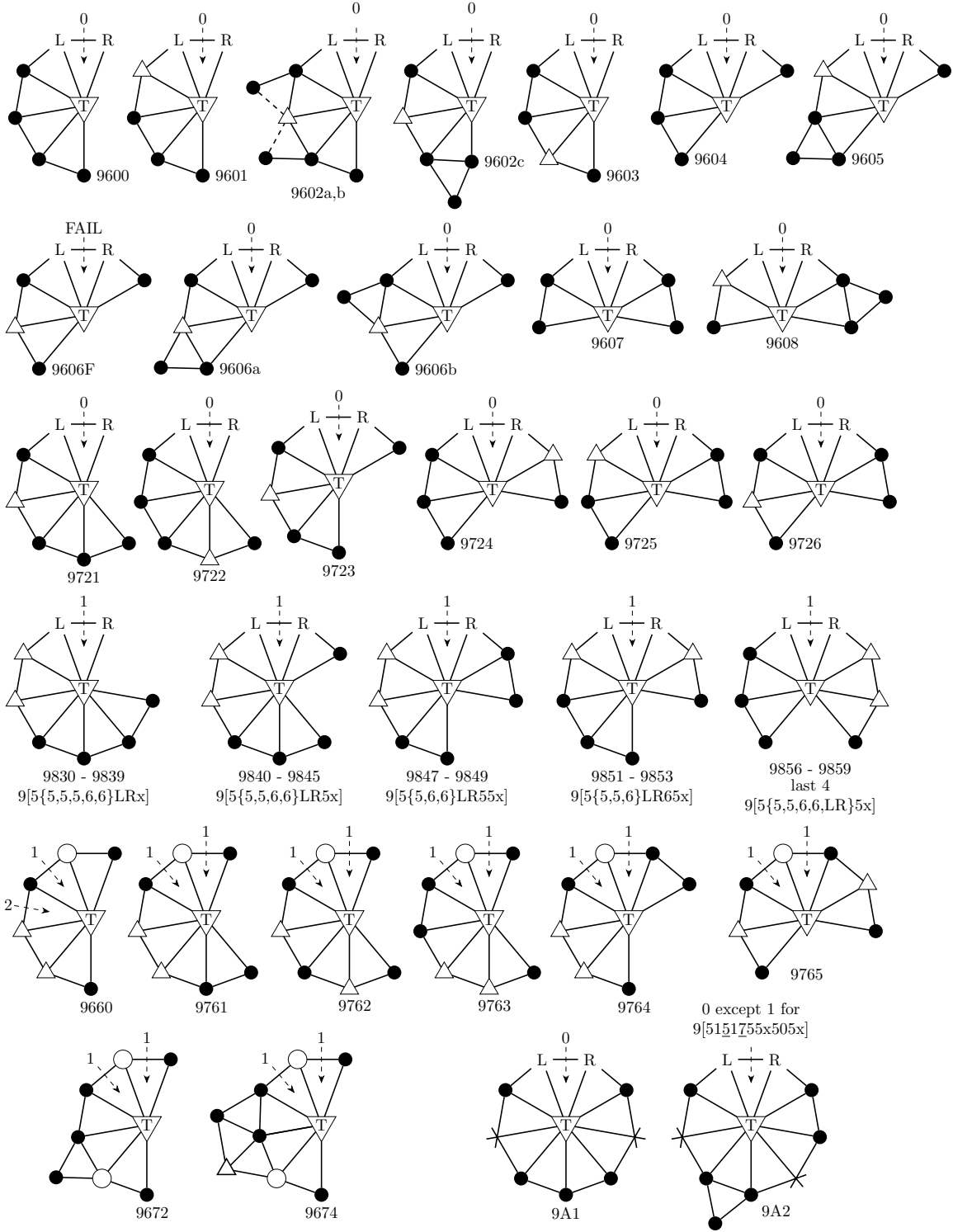

## Appendix F. Enumeration of Neighbourhoods of a 9-valent Vertex

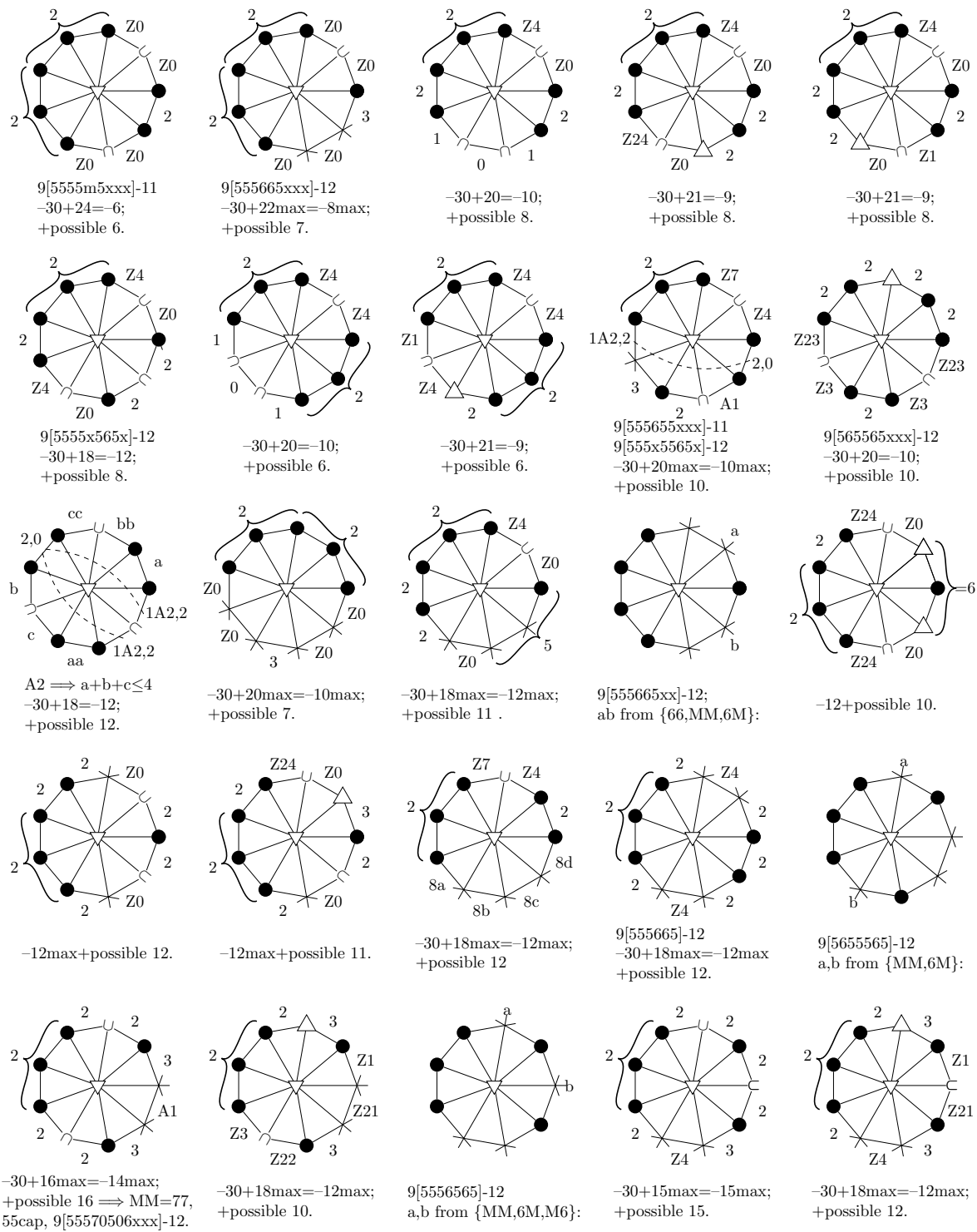

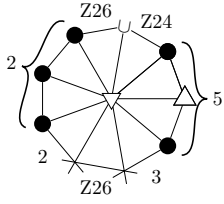

-30+18max=-12max;
+possible 12.

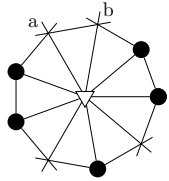

9[5565655]-12;
-30+18max=-12max;
a,b from {66,MM,6M}:

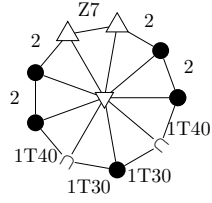

9[55565565]-13;
-30+16=-14;
+possible 12.

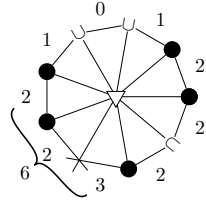

-30+16max=-14;
+possible 14

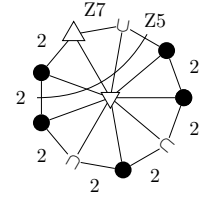

-30+15=-15;
+possible14.

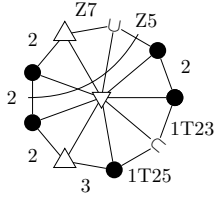

-30+17=-13;
+possible 13

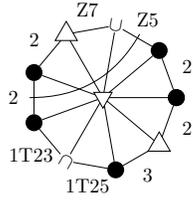

-30+17=-13;
+possible 13

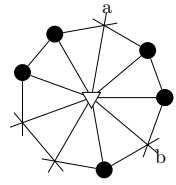

9[565565]-12
a,b from {MM,6M,M6}:

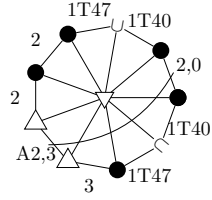

-30+16=-14;
+possible 14.

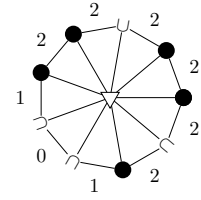

-30+14=-16;
+possible 14.

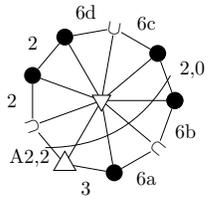

-30+15=-15;
+possible 15.

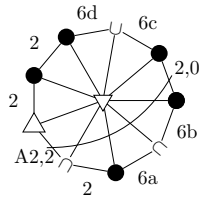

-30+15=-15;
+possible 14.

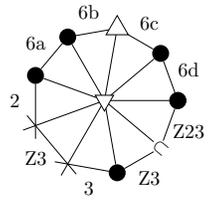

-30+18max=-12max;
+possible 11

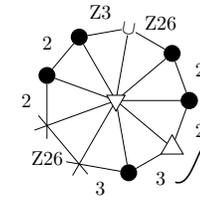

9[5565665xx]-13
-30+17max=-13max;
+possible 13.

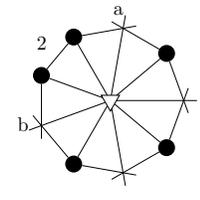

9[565565]-12,
9[5565x565x]-13;
a,b from {MM,6M}:

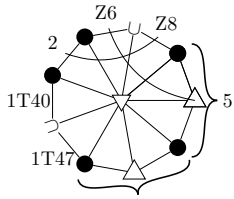

-30+16=-14;
+possible 14

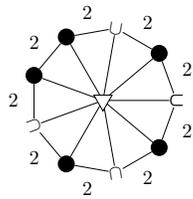

-30+12=-18;
+possible 18

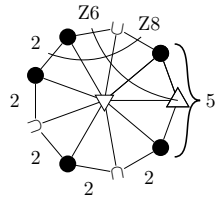

-30+14=-16;
+possible 15

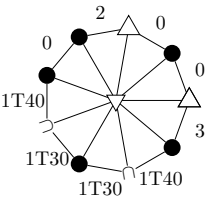

-30+16=-14;
0's from Table 5;
+possible 9

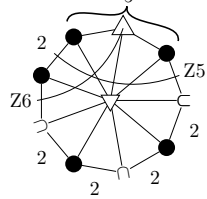

-30+14=-16;
+possible 15

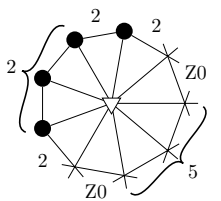

-30+16max=-14max;
+possible 13.

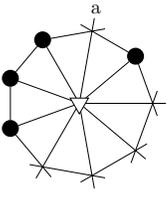

a from {6,M}

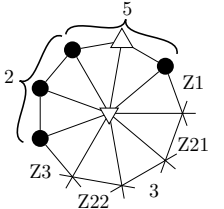

-30+16max=-14max;
+possible 10.

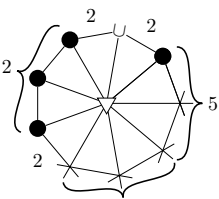

-30+14max=-16max;
+possible 18 ⟹
M5****=M5666*

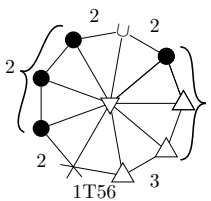

*=M⟹-30+13=-17
+possible 17,
else 9[55566665x]-14.

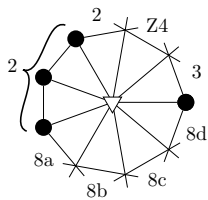
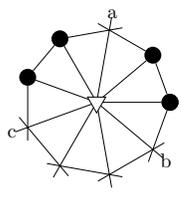
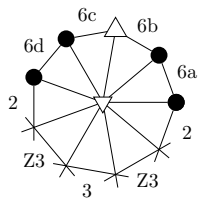
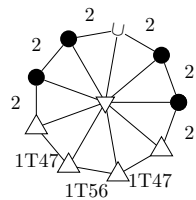
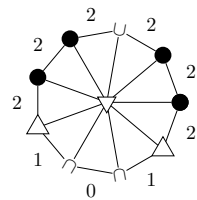

9[555665]-12;
−30+15max=−15max;
+possible 15.

a from {6,M}:
b,c from {66,MM,6M}:

−30+16max=−14max;
+possible 13.

−30+14=−16;
+possible 15;

−30+14=−16;
+possible 14.

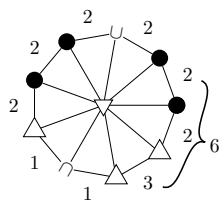
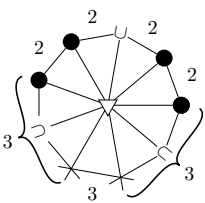
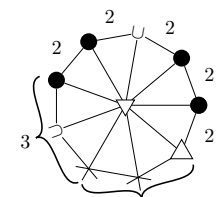
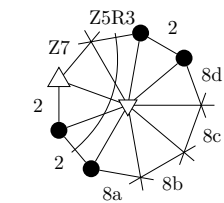
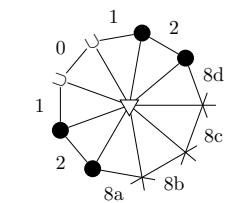

−30+14=−16;
+possible 16.

−30+12=−18;
+possible 17

−30+13=−17;+possible 18
⟹ 9[50517150516xxx]-13

−30+16max=−14max;
+possible 14.

−30+14max=−16max;
+possible 14.

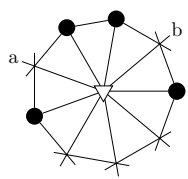
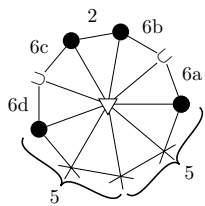
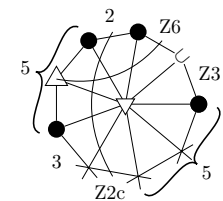
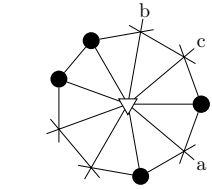
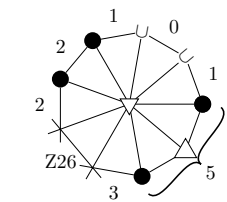

ab from {66,MM,6M};
9[565565]−12.

−30+12max=−18max
+possible 18.

−30+14max=−16max
+possible 15.

a from {6,M}
b,c from {MM,66,6M,M6}

−30+14=−16
+possible 14

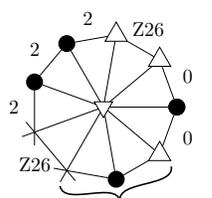
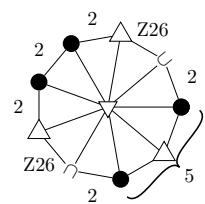
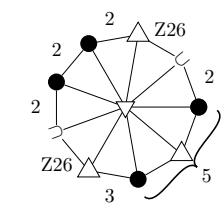
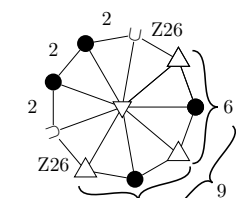
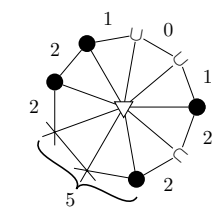

−30+16max=−14max
0's by Table 5;
+possible 12.

−30+14=−16
+possible 15

−30+14=−16
+possible 16

−30+14=−16
+possible 15

−30+12max=−18max;
+possible 15.

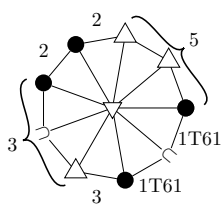
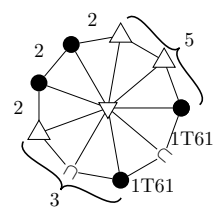
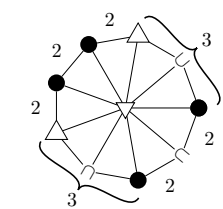
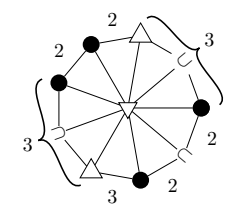
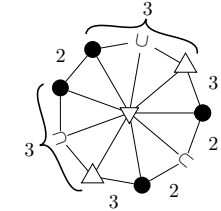

9[56655665x]-14
−30+13=−17;
+possible 17.

−30+13=−17;
+possible 16

−30+12=−18
+possible 16.

−30+12=−18
+possible 17.

−30+12=−18
+possible 18.

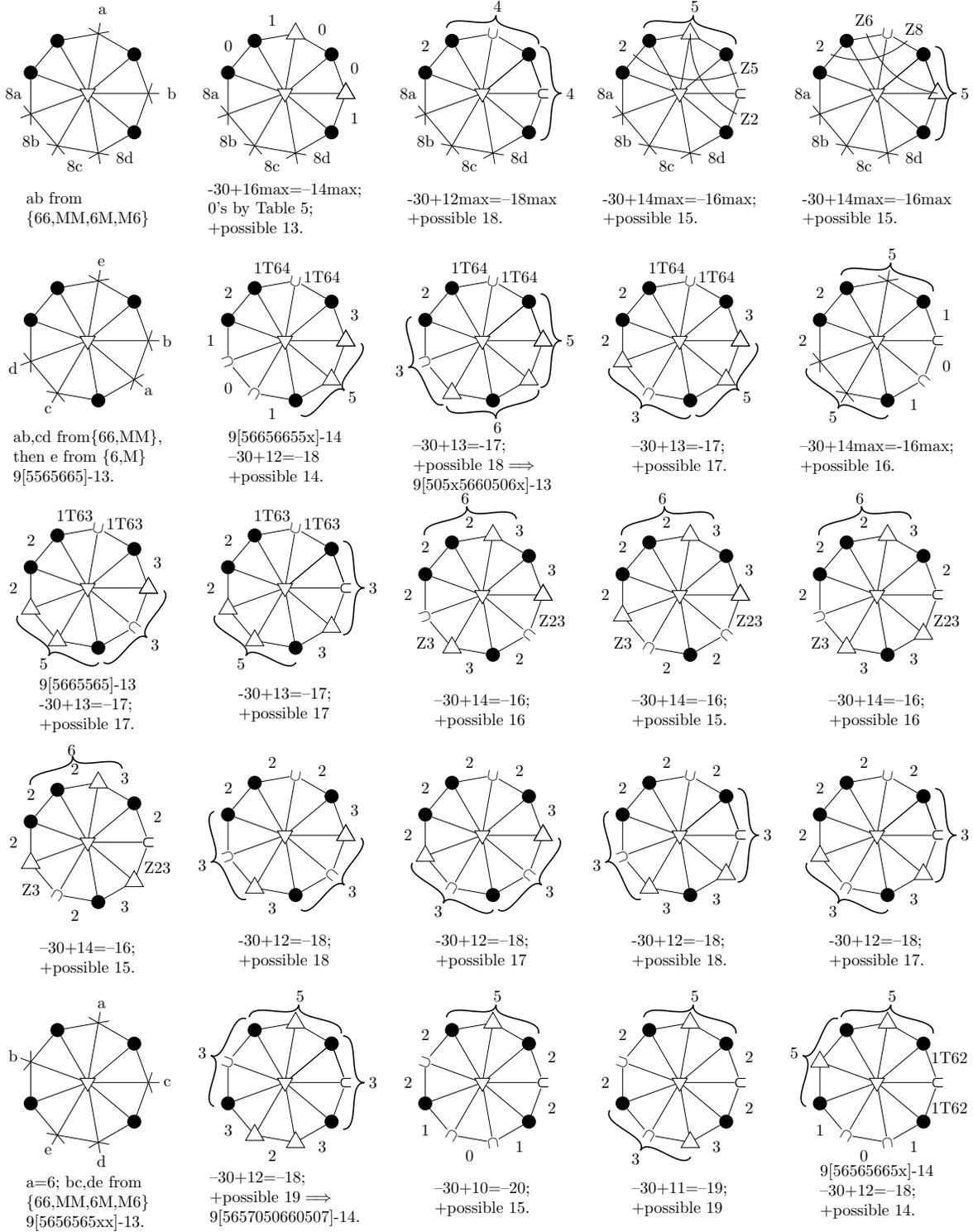

ab from
{66,MM,6M,M6}

−30+16max=−14max;
0's by Table 5;
+possible 13.

−30+12max=−18max
+possible 18.

−30+14max=−16max;
+possible 15.

−30+14max=−16max
+possible 15.

ab,cd from{66,MM},
then e from {6,M}
9[5565665]-13.

9[5665655x]-14
−30+12=−18
+possible 14.

−30+13=−17;
+possible 18 ⟹
9[505x5660506x]-13

−30+13=−17;
+possible 17.

−30+14max=−16max;
+possible 16.

9[5665565]-13
−30+13=−17;
+possible 17.

−30+13=−17;
+possible 17

−30+14=−16;
+possible 16

−30+14=−16;
+possible 15.

−30+14=−16;
+possible 16

−30+14=−16;
+possible 15.

−30+12=−18;
+possible 18

−30+12=−18;
+possible 17

−30+12=−18;
+possible 18.

−30+12=−18;
+possible 17.

a=6; bc,de from
{66,MM,6M,M6}
9[5656565xx]-13.

−30+12=−18;
+possible 19 ⟹
9[5657050660507]-14.

−30+10=−20;
+possible 15.

−30+11=−19;
+possible 19

9[56565665x]-14
−30+12=−18;
+possible 14.

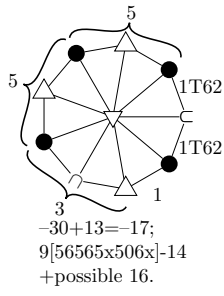

−30+13=−17;
9[56565x506x]-14
+possible 16.

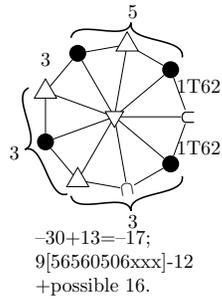

−30+13=−17;
9[56560506xxx]-12
+possible 16.

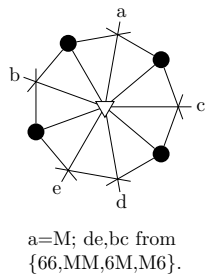

a=M; de,bc from
{66,MM,6M,M6}.

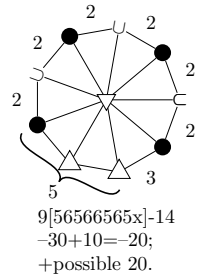

9[56566565x]-14
−30+10=−20;
+possible 20.

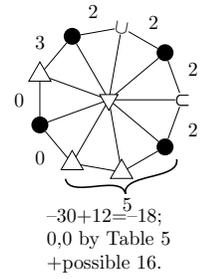

−30+12=−18;
0,0 by Table 5
+possible 16.

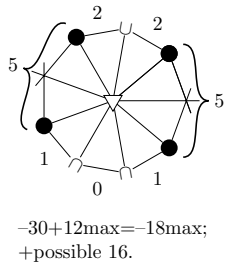

−30+12max=−18max;
+possible 16.

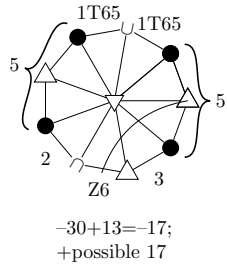

−30+13=−17;
+possible 17

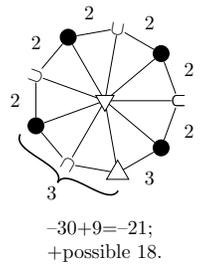

−30+9=−21;
+possible 18.

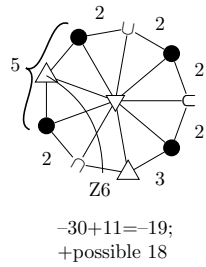

−30+11=−19;
+possible 18

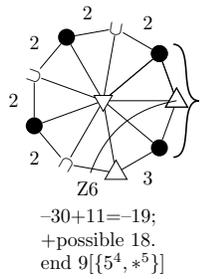

−30+11=−19;
+possible 18.
end 9[{5⁴,∗5⁵}]

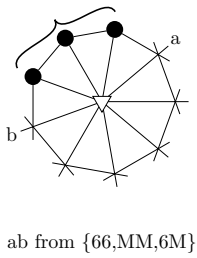

ab from {66,MM,6M}

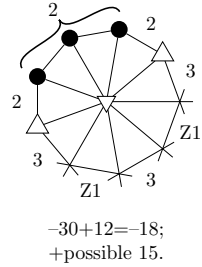

−30+12=−18;
+possible 15.

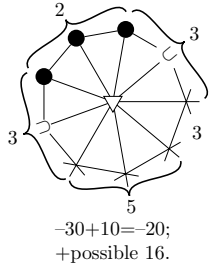

−30+10=−20;
+possible 16.

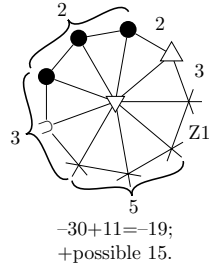

−30+11=−19;
+possible 15.

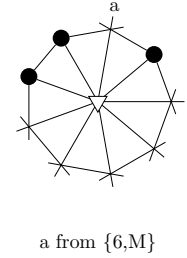

a from {6,M}

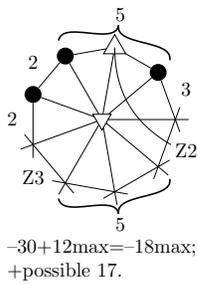

−30+12max=−18max;
+possible 17.

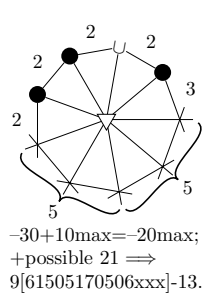

−30+10max=−20max;
+possible 21 ⟹
9[61505170506xxx]-13.

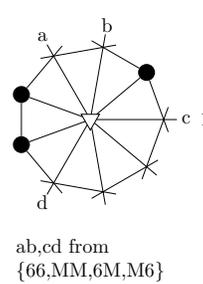

ab,cd from
{66,MM,6M,M6}

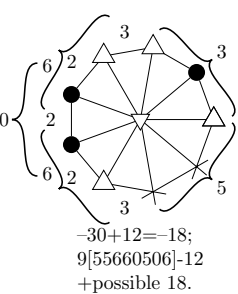

−30+12=−18;
9[55660506]-12
+possible 18.

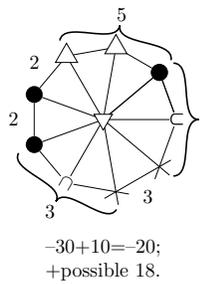

−30+10=−20;
+possible 18.

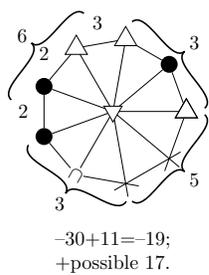

−30+11=−19;
+possible 17.

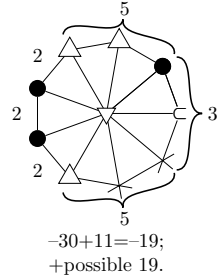

−30+11=−19;
+possible 19.

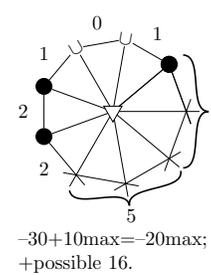

−30+10max=−20max;
+possible 16.

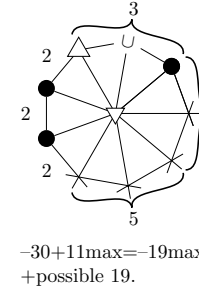

−30+11max=−19max;
+possible 19.

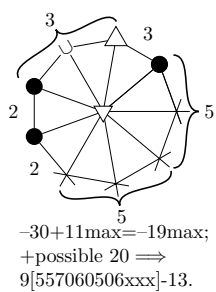

−30+11max=−19max;
+possible 20 ⟹
9[557060506xxx]-13.

ab,cd from
{66,MM,6M,M6}

$-30+12\max=-18\max$;
+possible 20 $\Longrightarrow$ *=6:

1,1,3 from Table 5;
$-18$+possible 17.

$-30+12\max=-18\max$;
+possible 16.

$-30+12\max=-18$.

$-30+10\max=-20\max$;
+possible 18.

$-30+11\max=-19\max$;
+possible 19.

ab,cd from
66, MM, 6M, M6

$-30+12=-18$;
9[56561506], 9[56560516]-13
+possible 18.

$-30+10=-20$;
+possible 19.

$-30+11=-19$;
+possible 18.

$-30+8\max=-22\max$;
+possible 21.

$-30+10\max=-20\max$;
+possible 22;xchg 5&3 $\Longrightarrow$
****=6666:

9[566660506xx]-14;
$-20$+possible 20.

a from {6,M},
bc from {66,MM,6M,M6}

$-30+12\max=-18\max$;
0,0 from Table 5;
+possible 18.

$-30+10\max=-20\max$;
+possible 17.

$-30+11\max=-19\max$;
+possible 21 $\Longrightarrow$
*** = 666

9[56660506]-13,
9[5666560506]-14;
$-19$+possible 17.

$-30+11\max=-19\max$;
+possible 21 $\Longrightarrow$
*** = 666

$-19$+possible 17.

$-30+10\max=-20$ max;
+possible 21 $\Longrightarrow$ M=7 $\Longrightarrow$
3 becomes 2T60.

$-30+8\max=-22\max$;
+possible 16

$-30+9\max=-21\max$;
+possible 20.

$-30+9\max=-21\max$;
+possible 20.

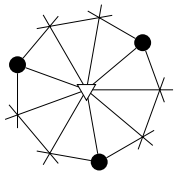

Try **=MM,
then ** **=66 66,
then at most one 66:

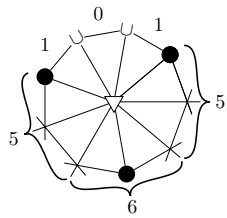

−30+10max=−20max
+possible 18.

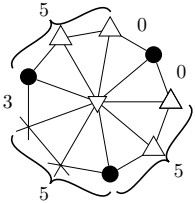

−30+12max=−18max
0,0 from Table 5;
+possible 18.

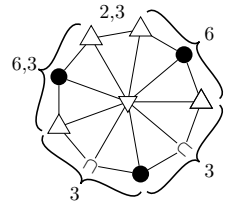

−30+10=−20;
+possible 20.

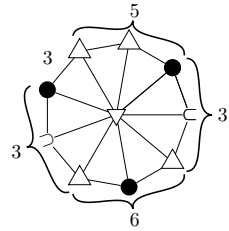

−30+10=−20;
+possible 20.

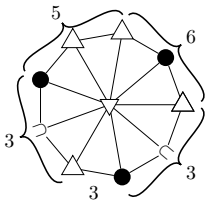

−30+10=−20;
+possible 20.
Else no MM, no 66:

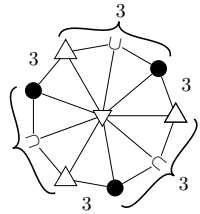

−30+9=−21;
+possible 18.

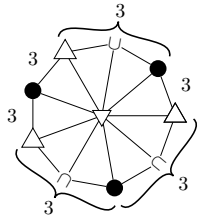

−30+9=−21;
+possible 18.
end of 9{5³,*⁶}.

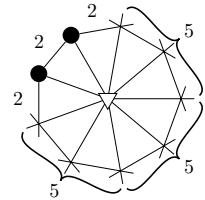

−30+8max=−22max;
+possible 21.

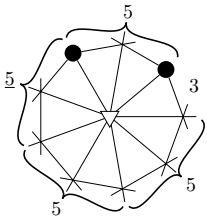

−30+8max=−22max;
+possible 23=+1:
8max ⟹ *5*=656:

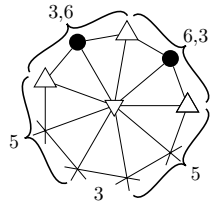

−22+possible 22.

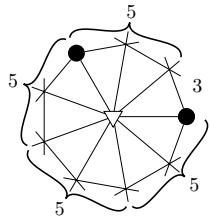

−30+8max=−22max;
+possible 23=+1:
8max ⟹ *5*=656:

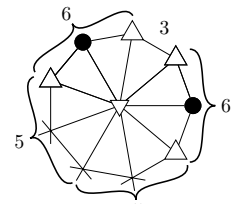

−22+possible 25 ⟹
9[6050660506]-12.

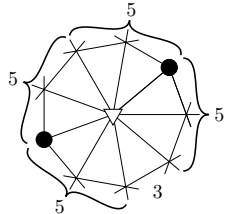

−30+8max=−22max;
+possible 23; ⟹
***=666, *5*=656:

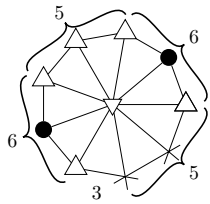

−22+possible 25 ⟹
9[60506660506]-13.

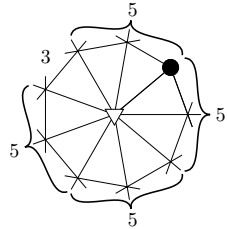

−30+4max=−26max;
+possible 23.

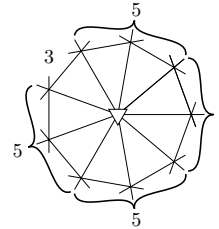

−30+possible 23.

APPENDIX G. ALL SOURCES

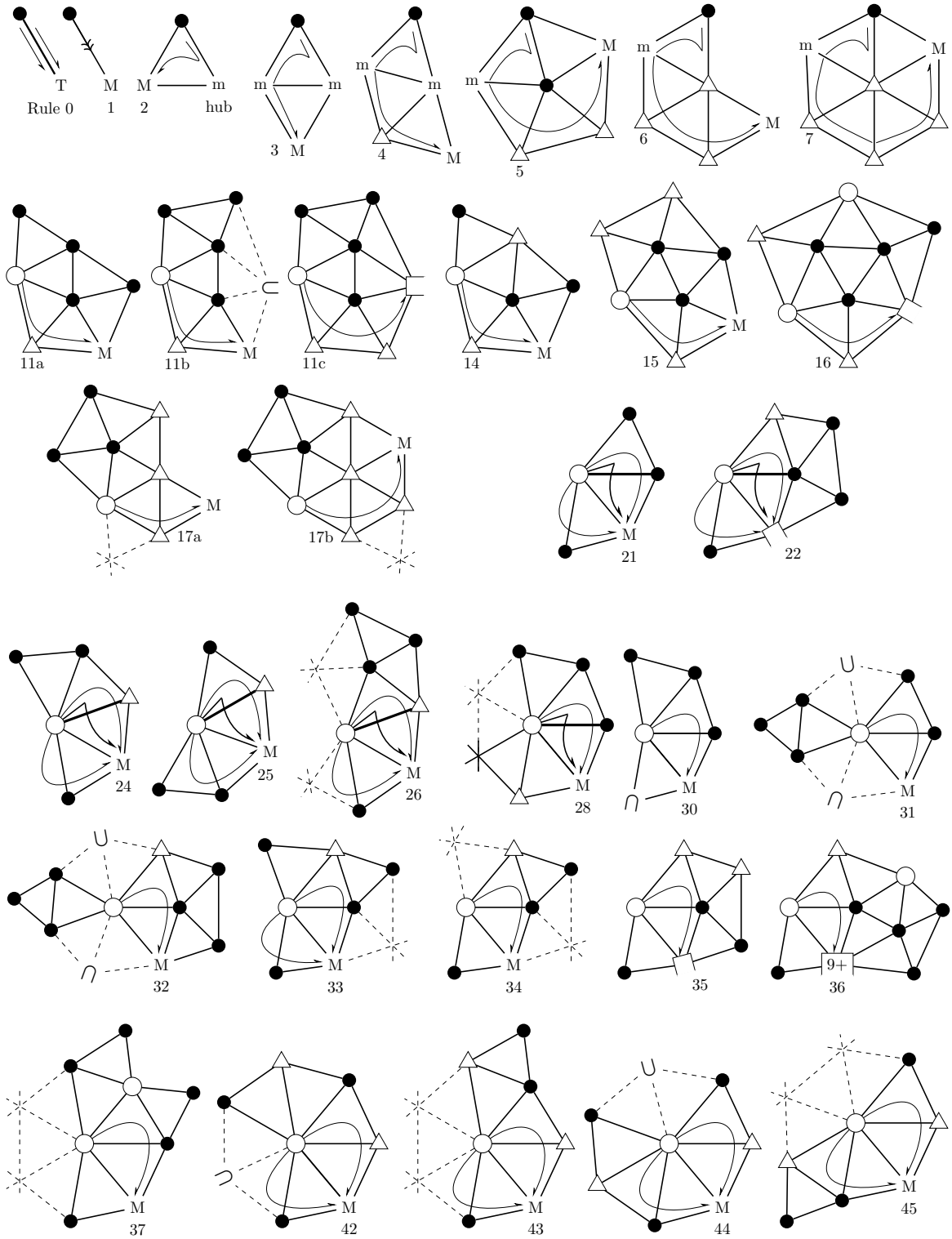

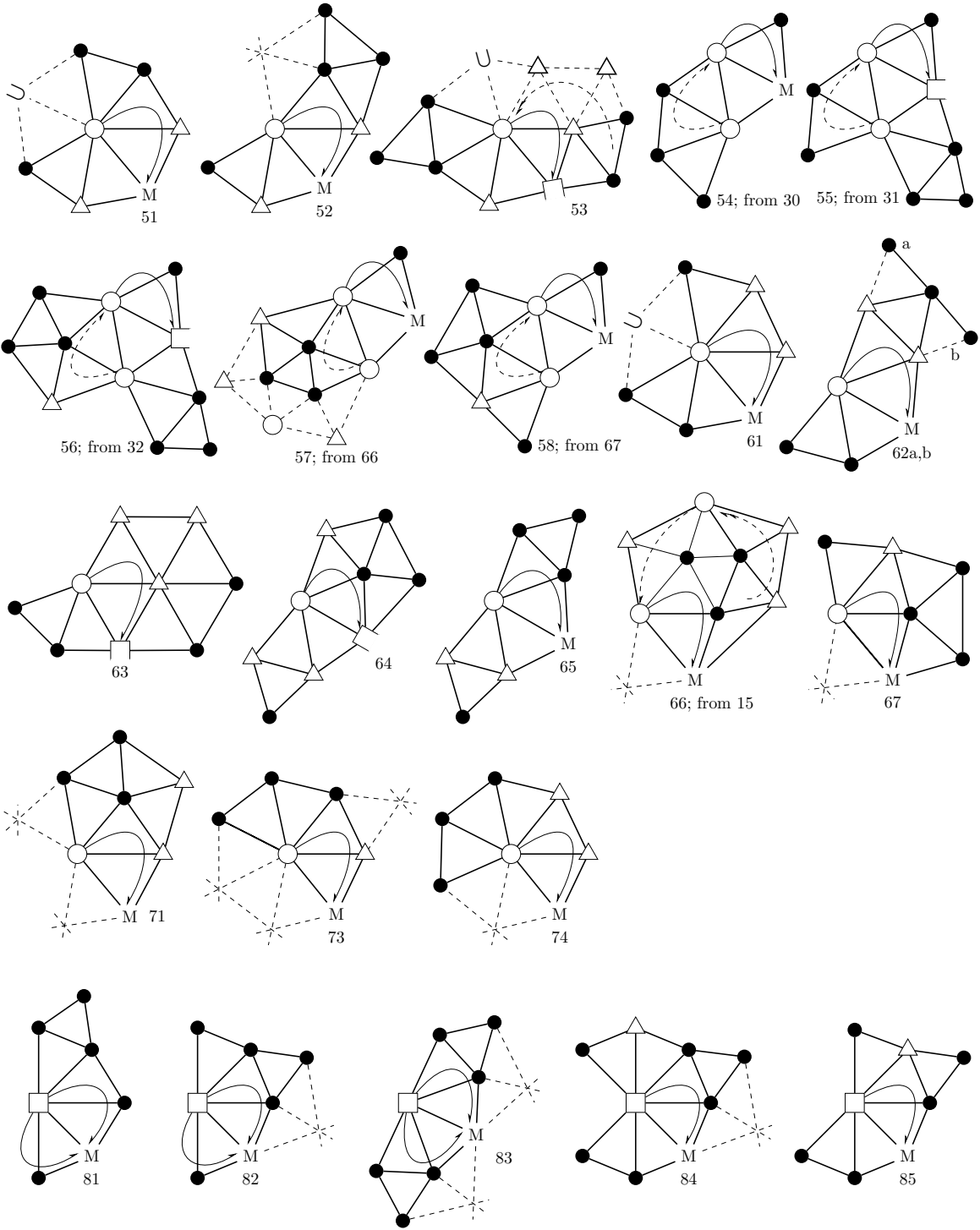

*Email address*: frank.allaire@lakeheadu.ca



\end{document}